\documentclass[12pt]{book}
 \usepackage{amsmath,amssymb,amsbsy,amstext,amscd,amsfonts}
 \usepackage{exscale,calrsfs,times,multind,theorem}
 \usepackage{color}
 \usepackage{article}
 \usepackage{pictex,epsfig}
 \makeatletter
 \renewcommand{\@evenhead}{\textsc{\thepage}\hfill {\em\scriptsize\rightmark}
        \hfill}
\renewcommand{\@oddhead}{ {\em\scriptsize \hfill A study of bounded operators on Martingale Hardy spaces} \hfill \textsc{\thepage}\,}
\renewcommand{\@evenfoot}{\scriptsize\em
\hfill  \hfill,}
\renewcommand{\@oddfoot}{\scriptsize\em
\hfill  G.Tutberidze \hfill,}

 \input article.def

\begin{document}
\def \R {{\mathbb {R}}}

\def \R {{\mathbb {R}}}

\newpage

\thispagestyle{empty}

\begin{center}
{\huge \textbf{Ph.D. Thesis}} \vspace{2cm}

{\Large \textbf{Bounded operators on Martingale Hardy spaces}} \vspace{0.5cm}

\end{center}

\vspace{3cm}

\begin{center}
\textbf{Giorgi Tutberidze}
\end{center}

\vspace{2.0cm}

\begin{center}
{\ The University of Georgia, School of Science and Technology \\[0pt]
77a Merab Kostava St, Tbilisi, 0128, Georgia \\[0pt]
and \\[0pt]
Department of Computer Science and Computational Engineering UiT - The Arctic University of Norway \\[0pt]
P.O. Box 385, N-8505, Narvik, Norway \\[0pt]
{giorgi.tutberidze1991@gmail.com}	}\\[0pt]
\vspace{8cm} \textbf{Tbilisi \\
	2021}
\end{center}

\newpage

\thispagestyle{empty}

\begin{center}
{\scriptsize Giorgi Tutberidze}
\end{center}

\newpage

\tableofcontents

\newpage

\thispagestyle{empty}


\newpage

\thispagestyle{empty}

\begin{center}
	{\scriptsize Giorgi Tutberidze}
\end{center}

\newpage

\begin{center}
{\Large \textbf{Abstract}}
\end{center}

\vspace{0.5cm}

The classical Fourier Analysis has been developed in an almost unbelieved way from the first fundamental discoveries by name Fourier. Especially a number of wonderful results have been proved and new directions of such research has been developed e.g. concerning Wavelets Theory, Gabor theory, Time-Frequency Analysis, Fast Fourier Transform, Abstract Harmonic Analysis, etc. One important reason for this is that this development is not only important for improving the "State of the art", but also for its importance in other areas of mathematics and also for several applications (e.g. theory of signal transmission, multiplexing, filtering,
image enhancement, coding theory, digital signal processing and pattern
recognition.)

The classical theory of Fourier series deals with decomposition of a
function into sinusoidal waves. Unlike these continuous waves the Vilenkin (Walsh) functions are rectangular waves. The development of the theory of Vilenkin-Fourier series has
been strongly influenced by the classical theory of trigonometric series.
Because of this it is inevitable to compare results of Vilenkin series to
those on trigonometric series. There are many similarities between these
theories, but there exist differences also. Much of these can be explained
by modern abstract harmonic analysis, which studies orthonormal systems from the point of view of the structure of a topological group.

The aim of my thesis is to discuss, develop and apply the newest developments of this fascinating theory connected to modern harmonic analysis. In particular, we investigate some strong convergence result of partial sums of Vilenkin-Fourier series. Moreover, we derive
necessary and sufficient conditions for the modulus of continuity so that norm convergence of subsequences of Fej\'er means is valid. Furthermore, we
consider Riesz and N\"orlund logarithmic means. It is also proved that these results are the best possible in a
special sense. As applications both some well-known and new results are pointed out.
In addition, we investigate some $T$ means, which are "inverse" summability methods of N\"orlund, but only in the case  when their coefficients are monotone. 

The thesis contains six chapters and two appendices, one contains some basic fact concerning classical Hardy spaces, in the other one is devoted to present about Kachmarz systems. One reason for this is that it will be more convenient for the reader to compare with the classical theory and another reason is that it give us possibility to raise new open questions. It is maybe surprising that some of these open questions concern classical situation but are motivated by the results we proved in this new situation. Totally we have explicitly pointed out 30?? open questions in this book. We hope that this can stimulate the further development of this fascinating area. We now continue by describing the main content of each of the chapters.

In Chapter 1 we first present some definitions and notations, which are
crucial for our further investigations. After that we also define some
summabilitity methods and remind about some classical facts and results. We investigate some well-known results and prove new estimates for the kernels of these summabilitity methods, which are very important to prove our main results. Moreover, we define martingale Hardy spaces and construct martingales, which help us to prove sharpness of our main results in the later chapters.

Chapter 2 is devoted to investigate some new strong convergence of partial sums and Fej\'er means with respect to Vilenkin systems. Next, we consider convergence of subsequences of Fej\'er means and prove some boundedness results for them. After that we apply these results to find necessary and sufficient conditions for the modulus of continuity for which norm convergence of Fej\'er means hold. We also prove sharpness of all our main results in this Chapter.

In Chapter 3 we consider boundedness of maximal operators of $T$ means  with respect to Vilenkin systems. We also prove that results are sharp in the special sense. After that we prove some strong convergence theorems for these summablility methods. Since Fej\'er means, Riesz means are well-know examples of $T$ means some well-known and new results are pointed out.

In Chapter 4 we consider Riesz and N\"orlund logarithmic means with respect to Vilenkin systems. In particular, we prove some strong convergence theorems of Riesz means with respect to Vilenkin systems. Moreover, we also prove sharpness of this result for only Walsh-Fourier series. Next, we  investigate boundedness of maximal operators of N\"orlund logarithmic means. We also obtain a.e. convergence of this means.

\newpage

\thispagestyle{empty}

\begin{center}
	{\scriptsize Giorgi Tutberidze}
\end{center}

\newpage

\begin{center}
{\Large \textbf{Preface}}
\end{center}

\vspace{0.70cm}

This PhD thesis is written as a monograph based on the following publications:
\vspace{0.30cm}

[1] G. Tutberidze, A note on the strong convergence of partial sums with respect to Vilenkin system, Journal of Contemporary Mathematical Analysis, 54, 6, (2019), 319–324.

[2] G. Tutberidze, Maximal operators of $T$ means with respect to the Vilenkin system, Nonlinear Studies, 27, 4, (2020), 1–11.

[3] L. E. Persson, G. Tephnadze, G. Tutberidze, On the boundedness of subsequences of Vilenkin-Fejér means on the martingale Hardy spaces, operators and matrices, 14, 1 (2020), 283–294.

[4] G. Tephnadze,  G. Tutberidze, A note on the maximal operators of the Nörlund logaritmic means of Vilenkin-Fourier series, Transactions of A. Razmadze Math. Inst., 174, 1 (2020), 107–112.

[5] D. Lukkassen, L.E. Persson, G. Tephnadze, G. Tutberidze, Some inequalities related to strong convergence of Riesz logarithmic means of Vilenkin-Fourier series, J. Inequal. Appl., 2020, DOI: https://doi.org/10.1186/s13660-020-02342-8.

[6] G. Tutberidze, Modulus of continuity and boundedness of subsequences of Vilenkin- Fej\'er
means in the martingale Hardy spaces, Georgian Mathematical Journal, (to appear).

[7] L. E. Persson, G. Tephnadze, G. Tutberidze, P. Wall, Strong summability result of Vilenkin-Fejér means on bounded Vilenkin groups, Ukr. Math. J., (to appear).

[8] G. Tutberidze, Sharp $\left( H_{p},L_{p}\right) $ type
inequalities of maximal operators of $T$ means with respect to Vilenkin systems with monotone coefficients, Journal Mathematical Analysis and Application, (to appear).

\vspace{0.70cm}

We also have some papers which are not include in this PhD thesis:
\vspace{0.30cm}

1) G. Tutberidze and V. Tsagaraeishvili, Multipliers of Absolute Convergence, Mat. Zametki,  105, 3, (2019), 433–443.

2) G. Tutberidze and V. Tsagaraeishvili, Absolute convergence factors of Lipshitz class functions for general Fourier series, Georgian Mathematical Journal, (to appear).

3) G. Tutberidze and V. Tsagaraeishvili, Multipliers of a.e convergence of general Fourier series, Ukr. Math. J., (to appear).

\vspace{0.70cm}
Remark: Also some new results which can not be found in these papers appear in this PhD thesis for the first time.

\newpage

\begin{center}
{\Large \textbf{Acknowledgement}}
\end{center}

\vspace{1.0cm}

It is a pleasure to express my warmest thanks to my supervisors Professors Lars-Erik Persson, Dag Lukkassen and George Tephnadze for the attention to my work, their valuable remarks and suggestions and for their constant support and help.

I am also grateful to Professor Natasha Samko and PhD student Singh Harpal for helping me with several practical things. Moreover, I appreciate very much the warm and friendly atmosphere at the Department of Mathematics at UiT The Arctic University of Norway, which helps me to do research more effectively.

I thank my Georgian colleagues from The University of Georgia for their interest in my research and for many fruitful discussions. 

I am also grateful to Professor Vakhtang Tsagareishvili from Tbilisi State University, for his support.

I also want to pronounce that the agreement about scientific collaboration and PhD education between The University of Georgia and UiT - The Arctic University of Norway has been very important. Especially, I express my deepest gratitude to Professor Roland Duduchava at The University of Georgia for his creative and hard work to realize this important agreement.

I thank Shota Rustaveli National Science Foundation for financial support, which was very important to do my research in the frame of this project.

Finally, I thank my family for their love, understanding, patience and long lasting constant	support.

\newpage

\thispagestyle{empty}

\begin{center}
{\scriptsize Giorgi Tutberidze}
\end{center}

\newpage

	\section{\textbf{Preliminaries}\protect\bigskip}

\vspace{0.5cm}

\subsection{Vilenkin groups and Functions}

Denote by $\mathbb{N}_+$ \index{N}{$\mathbb{N}_+$} the set of the positive integers, \index{s}{positive integers}
$\mathbb{N}:=\mathbb{N}_+\cup \{0\}.$ \index{N}{$\mathbb{N}$}
Let $m:=(m_0,m_1,\ldots)$ \index{N}{$m$} \index{N}{$m_k$} be a sequence of positive
integers not less than 2. Denote by \index{N}{$Z_{m_k}$}
\begin{equation*}
Z_{m_k}:=\{0,1,\ldots,m_k-1\}
\end{equation*}%
the additive group of integers modulo $m_k$ \index{S}{additive group of integers modulo $m_k$}.

Define the group $G_m$ \index{N}{$G_m$} as the complete direct product of the groups $Z_{m_k}$ with the product of the discrete topologies of $Z_{m_k}$.

The direct product $\mu $ of the measures
\begin{equation*}
\mu_k\left(j\right):=1/m_k\qquad(j\in Z_{m_k})
\end{equation*}%
is the Haar measure on $G_m$ with $\mu\left(G_m\right)=1.$

If $\sup_{n\in\mathbb{N}}m_n<\infty$, then we call $G_m$ a bounded Vilenkin \index{A}{Vilenkin} group.
If the generating sequence $m$ is not bounded, then $G_m$ is said to be an unbounded Vilenkin group.

In this book we discuss bounded Vilenkin groups,\index{S}{bounded Vilenkin groups} i.e. the case
when $\sup_{n\in\mathbb{N}}m_n<\infty.$

The elements of $G_m$ are represented by sequences
\begin{equation*}
x:=\left(x_0,x_1,\ldots,x_j,\ldots\right)\qquad\left(x_j\in Z_{m_j}\right).
\end{equation*}

Vilenkin group can be metrizable with the following metric:
\begin{equation*}
d\left(x,y\right):=\left\vert x-y\right\vert:=\sum_{k=0}^{\infty}
\frac{\left\vert x_k-y_k\right\vert}{M_{k+1}},\qquad\left( x\in G_m\right).
\end{equation*}

It is easy to give a base for the neighborhoods of $G_m:$  \index{N}{$I_0\left(x\right)$}  \index{N}{$I_n\left(x\right)$}
\begin{eqnarray*}
	I_0\left(x\right)&:&=G_m,\\
	I_n(x)&:&=\{y\in G_m\mid y_0=x_0,\ldots
	,y_{n-1}=x_{n-1}\} \qquad\left(x\in G_m,\ \ n\in\mathbb{N}\right).
\end{eqnarray*}%
Let  \index{N}{$ e_n$}
\begin{equation*}
e_n:=\left(0,\ldots,0,x_n=1,0,\ldots\right)\in G_m\qquad \left(n\in\mathbb{N}\right).
\end{equation*}

If we define $I_n:=I_n\left(0\right),$ \index{N}{$I_n$} for  $n\in\mathbb{N}$ and $\overline{I_n}:=G_m \ \ \backslash $ $I_n,$ \index{N}{$\overline{I_n} $} then%
\begin{equation} \label{1.1}
\overline{I_N}=\overset{N-1}{\underset{s=0}{\bigcup}}I_s\backslash
I_{s+1}=\left(\overset{N-2}{\underset{k=0}{\bigcup}}\overset{N-1}
{\underset{l=k+1}{\bigcup}}I_N^{k,l}\right)\bigcup
\left(\underset{k=1}{\bigcup\limits^{N-1}}I_N^{k,N}\right),
\end{equation}
where \index{N}{$I_N^{k,l}$}
\begin{equation*}
I_N^{k,l}:=\left\{ \begin{array}{l}I_N(0,\ldots,0,x_k\neq 0,0,...,0,x_l\neq 0,x_{l+1},\ldots ,x_{N-1},\ldots),\\
\text{for} \qquad k<l<N,\\
I_N(0,\ldots,0,x_k\neq 0,x_{k+1}=0,\ldots,x_{N-1}=0,x_N,\ldots ), \\
\text{for } \qquad l=N. \end{array}\right.
\end{equation*}

The norm (or quasi-norm when $0<p<1$) of the \index{A}{Lebesgue} \index{S}{Lebesgue spaces} Lebesgue space $L_p(G_m)$ \index{N}{$L_p(G_m)$} $\left(0<p<\infty \right)$ is defined by
\begin{equation*}
\left\Vert f\right\Vert_p:=\left(\int_{G_m}\left\vert f\right\vert
^pd\mu\right)^{1/p}.
\end{equation*}

The space $weak-L_p\left(G_m\right)$ \index{N}{$weak-L_p\left( G_m\right) $} consists of all measurable functions $f$, for which
\begin{equation*}
\left\Vert f\right\Vert_{weak-L_p}:=\underset{\lambda>0}{\sup}\lambda
\{\mu\left(f>\lambda\right)\}^{1/p}<+\infty.
\end{equation*}

The norm of the space of continuous functions $C(G_m)$ \index{N}{$C(G_m)$} is defined by
\begin{equation*}
\left\Vert f\right\Vert_C:=\underset{x\in G_m}{\sup}\left\vert
f(x)\right\vert<c<\infty.
\end{equation*}

The best approximation of $f\in L_p(G_m)$ $(1\leq p\leq\infty)$ is
defined as \index{N}{$E_n\left(f,L_p\right)$}
\begin{equation*}
E_n\left(f,L_p\right):=\inf_{\psi\in\emph{P}_n}\left\Vert f-\psi
\right\Vert_p,
\end{equation*}%
where $\emph{P}_n$ \index{N}{$\emph{P}_n$} is set of all Vilenkin polynomials \index{S}{Vilenkin polynomials} of order less than
$n\in \mathbb{N}$.

The modulus of continuity of functions of Lebesgue spaces \index{S}{The modulus of continuity of functions of Lebesgue spaces} $f\in L_p\left(G_m\right) $ and continous functions \index{S}{modulus of continuity of continous functions} $f\in C\left(G_m\right)$ are defined by \index{N}{$\omega_p\left(\frac{1}{M_n},f\right)$} \index{N}{$\omega_{C}\left(\frac{1}{M_n},f\right)$}
\begin{equation*}
\omega_p\left(\frac{1}{M_n},f\right):=\sup\limits_{h\in I_n} \left\Vert f\left(\cdot-h\right)-f\left(\cdot\right)\right\Vert_p
\end{equation*}
and
\begin{equation*}
\omega_C\left(\frac{1}{M_n},f\right):=\underset{h\in I_n}{\sup} \left\Vert f\left(\cdot-h\right)-f\left(\cdot\right)\right\Vert_C,
\end{equation*}
respectively.

If we define the so-called generalized number system based on $m$  in the
following way : \index{N}{$M_k$}
\begin{equation*}
M_0:=1, \ \  M_{k+1}:=m_kM_k \ \  (k\in\mathbb{N}),
\end{equation*}
then every $n\in\mathbb{N}$ can be uniquely expressed as
\begin{equation*}
n=\sum_{j=0}^{\infty}n_jM_j,
\end{equation*}
where $n_j\in Z_{m_j} \ \ (j\in\mathbb{N}_+)$ and only a finite number of $n_j^{^{\prime}}$`s differ from zero.

Next, we introduce on $G_m$ an orthonormal systems, which are called
Vilenkin systems.

At first, we define the complex-valued function
$r_k\left(x\right):G_m\rightarrow\mathbb{C},$
the generalized Rademacher functions, \index{S}{generalized Rademacher functions} by \index{N}{$r_k\left(x\right)$}
\begin{equation*}
r_k\left(x\right):=\exp\left(2\pi ix_{k}/m_{k}\right),\ \ \left(
i^{2}=-1, \ \ x\in G_{m},\ \ k\in\mathbb{N}\right).
\end{equation*}

Now, define Vilenkin systems \index{S}{Vilenkin systems}
$\psi:=(\psi_n:n\in\mathbb{N})$ on $G_m$ as: \index{N}{$\psi_n(x)$}
\begin{equation*}
\psi_n(x):=\prod\limits_{k=0}^{\infty}r_k^{n_k}\left(x\right), \ \ \left(n\in\mathbb{N}\right).
\end{equation*}

The Vilenkin systems are orthonormal and complete in $L_2\left(G_m\right).$

It is well-known that for all $n\in\mathbb{N},$ 
\begin{eqnarray} \label{vilenkin}
\left\vert \psi _{n}\left( x\right) \right\vert &=&1,\text{ } \\
\psi _{n}\left( x+y\right) &=&\psi _{n}\left( x\right) \text{\ }%
\psi _{n}\left( y\right) ,  \notag \\
\psi _{n}\left( -x\right) &=&\psi _{n^{\ast }}\left( x\right) =\overline{\psi }_{n}\left( x\right) ,\text{ \ \ \ }  \notag \\
\psi _{n}\left( x-y\right) &=&\psi _{n}\left( x\right) \text{\ }%
\overline{\psi} _{n}\left( y\right) ,  \notag \\
\psi _{n\widehat{+}k}\left( x\right) &=&\psi _{s}\psi _{n}\left( x\right) ,\text{\ \ }\left( s,n\in \mathbb{N},\text{ \ \ }x,y\in G_{m}\right).  \notag
\end{eqnarray}

Specifically, we call this system the Walsh-Paley system \index{S}{Walsh-Paley system} \index{A}{Walsh-Paley} when $m\equiv 2.$

\subsection{partial sums and Fejér means with respect to the Vilenkin systems}

Next, we introduce some analogues of the usual definitions in
Fourier-analysis. If $f\in L_1\left( G_m\right) $ we can define the
Fourier coefficients,  the partial sums  of Vilenkin-Fourier series, the Dirichlet
kernels,  Fejér means, 
Dirichlet and Fejér kernels with respect to the Vilenkin system with respect to Vilenkin systems 	in the usual manner:
\index{S}{Fourier coefficients with respect to Vilenkin systems}
\index{S}{partial sums  with respect to Vilenkin systems}
\index{S}{Fejér means with respect to Vilenkin systems}
\index{S}{Dirichlet kernels with respect to Vilenkin systems}
\index{S}{Fejér kernels with respect to Vilenkin systems}
\begin{eqnarray*} \index{N}{$\widehat{f}\left( n\right)$}
	\index{N}{$S_{n}f$} \index{N}{$D_{n}$}
	\widehat{f}\left(n\right)&:=&\int_{G_m}f\overline{\psi}_nd\mu,\qquad \left(n\in \mathbb{N}\right), \\
	S_nf&:=&\sum_{k=0}^{n-1}\widehat{f}\left(k\right)\psi_k, \qquad
	\left(n\in\mathbb{N}_+\right), \\
	\sigma_n f &:&=\frac{1}{n}\sum_{k=0}^{n-1}S_k f, \ \ \ \ \  \ \  \left(n\in \mathbb{N}_{+}\right) , \\
	D_n&:=&\sum_{k=0}^{n-1}\psi_{k}, \ \ \ \ \ \ \ \ \ \ \ \ \ \ \ \ \  \left(n\in\mathbb{N}_+\right), \\
	K_n &:&=\frac{1}{n}\overset{n-1}{\underset{k=0}{\sum }}D_{k}, \ \ \ \  \ \ \ \ \ \  \left(n\in \mathbb{N}_{+}\right) .
\end{eqnarray*}
respectively.

It is easy to see that%
\begin{eqnarray*}
	S_nf\left(x\right)&=&\int_{G_m}f\left(t\right)\sum_{k=0}^{n-1}\psi
	_k\left(x-t\right)d\mu \left(t\right) \\
	&=&\int_{G_m}f\left( t\right) D_n\left(x-t\right)d\mu\left(t\right)\\
	&=&\left(f\ast D_n\right)\left(x\right).
\end{eqnarray*}

It is well-known that (for details see e.g. \cite{sws}, \cite{gol} and \cite{AVD}) that for any $n\in\mathbb{N}$ and $1\leq s_n\leq m_n-1$ the following equalities holds \index{N}{$D_{s_nM_n}$}
\begin{equation} \label{dn21}
D_{j+M_n}=D_{M_n}+\psi_{M_n}D_j=D_{M_n}+r_nD_j, \qquad j\leq\left( m_n-1\right)M_n,
\end{equation}

\begin{eqnarray} \label{dn22}
D_{M_n-j}(x)&=&D_{M_n}(x)-\overline{\psi}_{M_n-1}(-x)D_j(-x)\\ \notag
&=&D_{M_n}(x)-\psi_{M_n-1}(x)\overline{D}_j(x), \qquad
j<M_n.
\end{eqnarray}

\begin{equation}  \label{3aa}
\quad \hspace*{0in}D_{M_{n}}\left( x\right) =\left\{ 
\begin{array}{l}
\text{ }M_{n}\text{ \ \ \ }x\in I_{n} \\ 
\text{ }0\text{ \qquad }x\notin I_{n}%
\end{array} 
\right.
\end{equation}

\begin{equation} \label{9dn}
D_{s_nM_n}=D_{M_n}\sum_{k=0}^{s_n-1}\psi_{kM_n}=D_{M_n}\sum_{k=0}^{s_n-1}r_n^k
\end{equation}
and
\begin{equation} \label{2dna}
D_n=\psi_n\left(\sum_{j=0}^{\infty}D_{M_j}\sum_{k=m_j-n_j}^{m_j-1}r_j^k\right).
\end{equation}

By using \eqref{3aa} we immediately get that
\begin{equation}  \label{5aa}
\Vert D_{M_n} \Vert_1=1<\infty.
\end{equation}  

It is obvious that%
\begin{eqnarray*}
	\sigma_nf\left(x\right) &=&\frac{1}{n}\overset{n-1}{\underset{k=0}{\sum}}
	\left(D_k\ast f\right)\left(x\right) \\
	&=&\int_{G_m}f\left(t\right)
	K_n\left(x-t\right)d\mu \left(t\right)\\
	&=&\left(f\ast K_n\right)\left(x\right).
\end{eqnarray*}
where $ K_n $ are the so called Fej\'er kernels.

It is well-known that (for details see e.g. \cite{gat}) for every 
$n>t,$ $t,n\in \mathbb{N}$ we have the following equality: \index{N}{$K_{M_n}$}
\begin{equation} \label{kn8}
K_{M_n}\left(x\right)=\left\{ \begin{array}{ll}
\frac{M_t}{1-r_t\left(x\right) },& x\in I_t\backslash I_{t+1},\qquad x-x_te_t\in I_n, \\ 
\frac{M_n+1}{2}, & x\in I_n, \\ 
0, & \text{otherwise. } \end{array} \right.
\end{equation}

Moreover, \index{N}{$K_{s_nM_n}$}
\begin{equation} \label{mag}
s_nM_nK_{s_nM_n}=\sum_{l=0}^{s_n-1}\left(\sum_{i=0}^{l-1}r_n^i\right)M_nD_{M_n}+\left(\sum_{l=0}^{s_n-1}r_n^{l}\right)M_nK_{M_n}
\end{equation}%
and

The next equality of Fej\'er kernels is very important for our further investigation (for details see Blahota and Tephnadze \cite{bt1}). In particular, if  $n=\sum_{i=1}^rs_{n_i}M_{n_i}$, where
$n_1>n_2>\dots>n_r\geq 0$ and $1\leq s_{n_i}<m_{n_i}$ \ for all
$1\leq i\leq r$ as well as $n^{(k)}=n-\sum_{i=1}^ks_{n_i}M_{n_i}$,
where $0<k\leq r$, then
\begin{equation} \label{kn10}
nK_n=\sum_{k=1}^r\left(\prod_{j=1}^{k-1}r_{n_j}^{s_{n_j}}\right)
s_{n_k}M_{n_k}K_{s_{n_k}M_{n_k}}+\sum_{k=1}^{r-1}\left(\prod_{j=1}^{k-1}r_{n_j}^{s_{n_j}}\right) n^{(k)}D_{s_{n_k}M_{n_k}}.
\end{equation}

It is well-known that 
\begin{equation} \label{knbounded}
\Vert K_n \Vert_1<c<\infty.
\end{equation}

Let define  maximal operators of partial sums and F\'ejer means by \index{S}{maximal operator of partial sums}
\index{S}{maximal operator of F\'ejer means}
\index{N}{${S}^{\ast}f$}
\index{N}{${\sigma}^{\ast}f$}
\begin{eqnarray*}
	{S}^{\ast}f&:=&\sup_{n\in\mathbb{N}}\left\vert S_{n}f\right\vert, \\
	{\sigma}^{\ast}f&:=&\sup_{n\in\mathbb{N}}\left\vert \sigma_{n}f\right\vert.
\end{eqnarray*}

Let define restricted maximal operators of partial sums and F\'ejer means by \index{S}{restricted maximal operator of partial sums}
\index{S}{restricted maximal operator of F\'ejer means}
\index{N}{$\widetilde{S}_{\#}^{\ast}f$}
\index{N}{$\widetilde{\sigma}_{\#}^{\ast}f$}
\begin{eqnarray*}
	\widetilde{S}_{\#}^{\ast}f&:=&\sup_{n\in\mathbb{N}}\left\vert S_{M_n}f\right\vert, \\
	\widetilde{\sigma}_{\#}^{\ast}f&:=&\sup_{n\in\mathbb{N}}\left\vert \sigma_{M_n}f\right\vert.
\end{eqnarray*}

\subsection{ Character $\rho\left( n\right)$ and Lebesgue constants with respect to Vilenkin Systems}

Let as define \index{N}{$\left\langle n\right\rangle$} \index{N}{$\left\vert n\right\vert$}
\begin{equation*}
\left\langle n\right\rangle:=\min \{j\in \mathbb{N}:n_j\neq 0\}\text{ \
	and \ }\left\vert n\right\vert :=\max \{j\in \mathbb{N}:n_{j}\neq 0\},
\end{equation*}%
that is
$M_{\left\vert n\right\vert}\leq n\leq M_{\left\vert n\right\vert+1}.$
Set
\begin{equation*}
\rho\left( n\right) :=\left\vert n\right\vert -\left\langle n\right\rangle ,\qquad \text{for all}\qquad n\in\mathbb{N}.
\end{equation*}

For the natural numbers $n=\sum_{j=1}^{\infty}n_jM_j$ and $k=\sum_{j=1}^{\infty}k_jM_j$ we define \index{N}{$\widehat{+}$}
\begin{equation*}
n\widehat{+}k:=\sum_{i=0}^{\infty}{\left(n_i\oplus k_i\right)}{M_{i+1}}
\end{equation*}
and \index{N}{$\widehat{-}$}
\begin{equation*}
n\widehat{-}k:=\sum_{i=0}^{\infty}{\left(n_i\ominus k_i\right)} {M_{i+1}},
\end{equation*}
where
\begin{equation*}
a_i\oplus b_i:=(a_i+b_i)\text{mod}m_i, \qquad a_i, b_i \in Z_{m_i}
\end{equation*}
and $\ominus $ \index{N}{$\ominus$} is the inverse operation for  \index{N}{$\oplus$} $ \oplus $.

For the natural number $n=\sum_{j=1}^{\infty}n_jM_j,$ we define
functions $v$ \index{N}{$v$} and $v^{\ast}$ \index{N}{$v^{\ast}$} by
\begin{equation*}
v\left(n\right):=\sum_{j=1}^{\infty}\left\vert\delta_{j+1}-\delta_j\right\vert+\delta_0,\qquad v^{\ast }\left(n\right)
:=\sum_{j=1}^{\infty}\delta_j^{\ast},
\end{equation*}%
where
\begin{equation*}
\delta_j=sign\left(n_j\right)=sign\left(\ominus n_j\right)\text{\ and \ }\delta_j^{\ast }=\left\vert \ominus n_j-1\right\vert\delta_j.
\end{equation*}

The $n$-th Lebesgue constant \index{S}{Lebesgue constant} is defined in the following way: \index{N}{$L_n$}
\begin{equation*}
L_n:=\left\Vert D_n\right\Vert_1.
\end{equation*}

For the trigonometric system
it is important to note that results of Fejér and Szego,  latter given in \cite{Szego} an explicit formula for Lebesgue constants.
The most properties of
Lebesgue constants with respect to the Walsh-Paley system were obtained
by Fine in \cite{fi}. In (\cite{sws}, p. 34), the two-sided estimate
is proved. In \cite{luko}, Lukomskii presented the lower estimate with sharp constant $1/4$.  Malykhin, Telyakovskii and Kholshchevnikova \cite{MTK} (see also Astashkin and Semenov \cite{AS}) improved the estimation above and proved  sharp estimate with factor 1. The first proof of \eqref{var1} with some different constants can be found in Lukomskii \cite{luko}, a new and shorter proof discribed in this book  which improved upper bound and provide similar lower bound can be found in \cite{bpt2}. In particular, for $\lambda :=\sup_{n\in \mathbb{N}}$ and for any  $n=\sum_{i=1}^{\infty }n_iM_i$ and $m_n$ we have the following two sided estimate:
\begin{equation}\label{var1}
\frac{1}{4\lambda}v\left(n\right)+\frac{1}{\lambda^2 }v^{\ast} \left(n\right)\leq L_n\leq v\left(n\right)+v^{\ast}\left(n\right).
\end{equation}

Moreover, (for details see in Memic, Simon and Tephnadze \cite{MST})
\begin{equation}\label{var2}
\frac{1}{nM_n}\underset{k=1}{\overset{M_n-1}{\sum }}v\left(k\right)
\geq\frac{2}{\lambda^2},
\end{equation}

Inequality \eqref{var1} immediately follows that for any $n\in \mathbb{N}$,  there exists an absolute constant $c,$ such that
\begin{equation} \label{Dn}
\left\Vert D_{n}\right\Vert_1\leq c\log n.
\end{equation}

For example, if we take  $q_{n_k}=M_{2n_k}+M_{2n_k-2}+M_2+M_0$, we have following two-sided inequality\index{N}{$D_{q_{n_k}}$}
\begin{equation}\label{Dnqn}
\frac{n_k}{2\lambda}\leq\left\Vert D_{q_{n_k}}\right\Vert_1\leq\lambda n_k, \ \ \ \ \ \ \lambda :=\sup_{n\in\mathbb{N}}m_n.
\end{equation}

\subsection{Definition and examples of N\"orlund and $T$ means and its maximal
	operators}

Let $\{q_k:k\in \mathbb{N}\}$ \index{N}{$q_k$} be a sequence of nonnegative numbers. The $n$-th N\"orlund means \index{S}{N\"orlund means} for the Fourier series of $f$  is defined by \index{N}{$t_nf$}
\begin{equation} \label{1.2}
t_nf:=\frac{1}{Q_n}\overset{n}{\underset{k=1}{\sum }}q_{n-k}S_kf,
\end{equation}
where   \index{N}{$Q_n$}
\begin{equation*}
Q_n:=\sum_{k=0}^{n-1}q_k.
\end{equation*}

A representation
\begin{equation*}
t_nf\left(x\right)=\underset{G}{\int}f\left(t\right)A_n\left(x-t\right) d\mu\left(t\right)
\end{equation*}%
plays a central role in the sequel, where \index{N}{$A_n$}
\begin{equation*}
A_n:=\frac{1}{Q_n}\overset{n}{\underset{k=1}{\sum }}q_{n-k}D_k
\end{equation*}
is the so-called N\"orlund kernel. \index{S}{N\"orlund kernel}

In  Moore \cite{moo} (see also Tephnadze \cite{tep12}) it was found necessary and suficient condition of  regularity of Norlund means. In particular, if $\{q_k:k\geq 0\}$ be a sequence of nonnegative numbers, $q_0>0$
and
\begin{equation*}
\lim_{n\rightarrow\infty}Q_n=\infty.
\end{equation*}
Then the summability method (\ref{1.2}) generated by $\{q_k:k\geq 0\}$ is regular if and only if
\begin{equation}
\underset{n\rightarrow \infty}{\lim}\frac{q_{n-1}}{Q_n}=0. \label{112}
\end{equation}

In addition, if sequence $\{q_k:k\in\mathbb{N}\}$ be non-increasing, then the summability method generated by
$\{q_k:k\in \mathbb{N}\}$ is regular, but if sequence $\{q_k:k\in \mathbb{N}\}$ be non-decreasing, then the summability method generated by $\{q_k:k\in\mathbb{N}\}$ is not always regular.

Let $\{q_{k}:k\geq 0\}$ be a sequence of non-negative numbers. The $n$-th $T$ means for a Fourier series of $f$  is defined by
\begin{equation} \label{nor}
T_nf:=\frac{1}{Q_n}\overset{n-1}{\underset{k=0}{\sum }}q_{k}S_kf,
\end{equation}
where $Q_{n}:=\sum_{k=0}^{n-1}q_{k}.$ 
It is obvious that 
\begin{equation*}
T_nf\left(x\right)=\underset{G_m}{\int}f\left(t\right)F_n\left(x-t\right) d\mu\left(t\right),
\end{equation*}
where $	F_n:=\frac{1}{Q_n}\overset{n}{\underset{k=1}{\sum }}q_{k}D_k$ is called the kernel of  $T$ means.

We always assume that $\{q_k:k\geq 0\}$ is a sequence of non-negative numbers and $q_0>0.$ Then the summability method (\ref{nor}) generated by $\{q_k:k\geq 0\}$ is regular if and only if $	\lim_{n\rightarrow\infty}Q_n=\infty.$

Let $t_n$ be N\"orlund means with monotone and bounded sequence
$\{q_k:k\in \mathbb{N}\}$, such that
\begin{equation*}
q:=\lim_{n\rightarrow\infty}q_n>c>0.
\end{equation*}

Then, if the sequence $\{q_k:k\in\mathbb{N}\}$ is non-decreasing, we get
that
\begin{equation*}
nq_0\leq Q_n\leq nq.
\end{equation*}

In the case when the sequence $\{q_{k}:k\in \mathbb{N}\}$ is non-increasing,
then
\begin{equation}
nq\leq Q_{n}\leq nq_{0}.  \label{monotone0}
\end{equation}

In both cases we can conclude that
\begin{equation} \label{monotone1}
\frac{q_{n-1}}{Q_n}=O\left(\frac{1}{n}\right),\text{ \  when \ }%
n\rightarrow \infty.
\end{equation}

One of the most well-known summability method which is example of Norlund and $T$ means is Fej\'er means, which is given when $\{q_k=1:k\in \mathbb{N}\}$
\begin{equation*}
\sigma_nf:=\frac{1}{n}\sum_{k=1}^nS_kf.
\end{equation*}

The $\left(C,\alpha\right)$-means (Ces\`aro means) \index{S}{Ces\`aro means} \index{N}{$\left(C,\alpha\right)$} of the Vilenkin-Fourier
series are defined by \index{N}{$\sigma_n^{\alpha}f$}
\begin{equation*}
\sigma_n^{\alpha}f:=\frac{1}{A_n^{\alpha}}\overset{n}{\underset{k=1}{
		\sum}}A_{n-k}^{\alpha-1}S_kf,
\end{equation*}
where \index{N}{$A_n^{\alpha}$}
\begin{equation*}
A_0^{\alpha}:=0,\qquad A_n^{\alpha}:=\frac{\left(\alpha+1\right)...\left(\alpha+n\right)}{n!},\qquad \alpha \neq -1,-2,...
\end{equation*}

It is well-known that (see e.g. Zygmund \cite{13z}) \index{N}{$A_n^{\alpha}$}
\begin{equation} \label{node0}
A_n^{\alpha}=\overset{n}{\underset{k=0}{\sum}}A_{n-k}^{\alpha-1},
\end{equation}%
\begin{equation} \label{node01}
A_n^{\alpha}-A_{n-1}^{\alpha}=A_n^{\alpha-1},\ \
A_{n}^{\alpha }\backsim n^{\alpha }.
\end{equation}

We also consider the "inverse" $\left(C,\alpha\right)$-means, which is an example of a $T$-means:
\begin{equation*}
U_n^{\alpha}f:=\frac{1}{A_n^{\alpha}}\overset{n-1}{\underset{k=0}{\sum}}A_{k}^{\alpha-1}S_kf, \qquad 0<\alpha<1.
\end{equation*}

Let $V_n^{\alpha}$ denote
the $T$ mean, where $	\left\{q_0=1, \  q_k=k^{\alpha-1}:k\in \mathbb{N}_+\right\} ,$
that is 
\begin{equation*}
V_n^{\alpha}f:=\frac{1}{Q_n}\overset{n}{\underset{k=1}{\sum }}k^{\alpha-1}S_kf,\qquad 0<\alpha<1.
\end{equation*}

The $n$-th N\"orlund logarithmic mean \index{S}{N\"orlund logarithmic mean} $L_n$ and the Riesz logarithmic mean \index{S}{Riesz logarithmic mean} $%
R_{n}$ are defined by  \index{N}{$L_nf$} \index{N}{$R_nf$}
\begin{eqnarray*}
	L_nf&:=&\frac{1}{l_n}\sum_{k=1}^{n-1}\frac{S_kf}{n-k},\\
	R_nf&:=&\frac{1}{l_n}\sum_{k=1}^{n-1}\frac{S_kf}{k},
\end{eqnarray*}
respectively, where \index{N}{$l_n$}
\begin{equation*}
l_n:=\sum_{k=1}^{n-1}\frac{1}{k}.
\end{equation*}

Kernels of N\"orlund logarithmic mean and Riesz logarithmic mean are  respectively defined by
\index{N}{$P_nf$} \index{N}{$Y_nf$}
\begin{eqnarray*}
	P_nf&:=&\frac{1}{l_n}\sum_{k=1}^{n-1}\frac{D_kf}{n-k},\\
	Y_nf&:=&\frac{1}{l_n}\sum_{k=1}^{n-1}\frac{D_kf}{k}.
\end{eqnarray*}

It is well-known that 
\begin{eqnarray} \label{reisz}
\Vert	Y_nf \Vert_1<c<\infty.
\end{eqnarray}

Up to now we have considered N\"orlund and $T$ means in the case when the sequence $\{q_k:k\in\mathbb{N}\}$ is bounded but now we consider N\"orlund and $T$ summabilities with unbounded sequence $\{q_k:k\in\mathbb{N}\}$.

Let $\alpha\in
\mathbb{R}_+,\ \ \beta\in\mathbb{N}_+$ and
\begin{equation*}
\log^{(\beta)}x:=\overset{\beta\text{times}}{\overbrace{\log ...\log}}x.
\end{equation*}

If we define the sequence $\{q_k:k\in \mathbb{N}\}$ by
\begin{equation*}
\left\{q_0=0\text{ \  and \ }q_k=\log^{\left(\beta \right)}k^{\alpha
}:k\in\mathbb{N}_+\right\},
\end{equation*}%
then we get the class of N\"orlund means with non-decreasing coefficients: \index{N}{$\kappa_n^{\alpha,\beta}f$}
\begin{equation*}
\kappa_n^{\alpha,\beta}f:=\frac{1}{Q_n}\sum_{k=1}^{n}\log^{\left(
	\beta\right)}\left( n-k\right)^{\alpha}S_kf.
\end{equation*}%
First we note that $\kappa_n^{\alpha,\beta}$ are
well-defined for every
$n\in\mathbb{N}_+$, if we rewrite them as:
\begin{equation*}
\kappa_n^{\alpha,\beta}f=\sum_{k=1}^n\frac{\log^{\left(\beta
		\right)}\left(n-k\right)^{\alpha }}{Q_n}S_kf.
\end{equation*}

It is obvious that
\begin{equation*}
\frac{n}{2}\log^{\left(\beta \right)}\frac{n^{\alpha }}{2^{\alpha }}\leq Q_n\leq n\log^{\left(\beta\right)}n^{\alpha}.
\end{equation*}

It follows that
\begin{eqnarray} \label{node00}
\frac{q_{n-1}}{Q_n}&\leq&\frac{c\log^{\left(\beta\right)}\left(
	n-1\right)^{\alpha}}{n\log^{\left(\beta\right) }n^{\alpha}} \\ \notag
&=& O\left(\frac{1}{n}\right)\rightarrow 0,\text{ \ as \ }n\rightarrow \infty.
\end{eqnarray}

If we define the sequence $\{q_k:k\in \mathbb{N}\}$ by
$	\left\{q_0=0, \ q_k=\log^{\left(\beta \right)}k^{\alpha
}:k\in\mathbb{N}_+\right\},$ 
then we get the class of $T$ means with non-decreasing coefficients:
\begin{equation*}
B_n^{\alpha,\beta}f:=\frac{1}{Q_n}
\sum_{k=1}^{n}\log^{\left(\beta\right)}k^{\alpha}S_kf.
\end{equation*}%

We note that $B_n^{\alpha,\beta}$ are
well-defined for every $n \in \mathbb{N}$

\begin{equation*}
B_n^{\alpha,\beta}f=\sum_{k=1}^{n}\frac{\log^{\left(\beta\right)}k^{\alpha }}{Q_n}S_kf.
\end{equation*}

It is obvious that $\frac{n}{2}\log^{\left(\beta \right)}\frac{n^{\alpha }}{2^{\alpha }}\leq Q_n\leq n\log^{\left(\beta\right)}n^{\alpha}\rightarrow 0,\text{ \ as \ }n\rightarrow \infty.$

Let us define maximal operator of N\"orlund and $T$ means by \index{N}{$t^{\ast}f$} \index{N}{$T^{\ast}f$}
\index{S}{maximal operator of N\"orlund  means}
\index{S}{maximal operator of$T$ means}
\begin{eqnarray*}
	t^{\ast}f:=\sup_{n\in\mathbb{N}}\left\vert t_nf\right\vert, \\
	T^{\ast}f:=\sup_{n\in\mathbb{N}}\left\vert T_nf\right\vert
\end{eqnarray*}

The well-known examples of maximal operator of N\"orlund and $T$ means are maximal operator of  Ces\'aro means, N\"orlund anD Reisz logarithmic means which are defined by:  
\index{S}{maximal operator of Ces\'aro means} 
\index{S}{maximal operator of N\"orlund logarithmic means}
\index{S}{maximal operator of Reisz logarithmic means} 
\index{N}{$\sigma^{\alpha,\ast}f$}  
\index{N}{$L^{\ast}f$}
\index{N}{$R^{\ast}f$}
\begin{eqnarray*}
	\sigma^{\alpha ,\ast}f&:=&\sup_{n\in\mathbb{N}}\left\vert \sigma _n^{\alpha }f\right\vert, \\
	L^{\ast}f&:=&\sup_{n\in\mathbb{N}}\left\vert L_nf\right\vert, \\
	R^{\ast}f&:=&\sup_{n\in\mathbb{N}}\left\vert R_{n}f\right\vert.
\end{eqnarray*}

We also define some new maximal operators: \index{S}{maximal operators} \index{N}{$\kappa^{\alpha,\beta,\ast}f$} \index{N}{$\beta
	^{\alpha,\ast }f$}
\begin{eqnarray*}
	\kappa ^{\alpha,\beta,\ast}f&:=&\sup_{n\in
		\mathbb{N}}\left\vert \kappa_n^{\alpha ,\beta }f\right\vert, \\ \beta^{\alpha,\ast }f&:=&\sup_{n\in\mathbb{N}}\left\vert \beta_{n}^{\alpha}f\right\vert.
\end{eqnarray*}

\subsection{Weak-type and strong-type inequalities and  a.e convergence}

The convolution of two functions
$f,g\in L_{1}(G_m)$ is defined by
\begin{equation*}
\left( f\ast g\right) \left( x\right) :=\int_{G_m}f\left( x-t\right)
g\left( t\right) dt\text{ \ \ }\left( x\in G_m\right).
\end{equation*}%
It is easy to see that
\begin{equation*}
\left( f\ast g\right) \left( x\right) =\int_{G_m}f\left( t\right)
g\left( x-t\right) dt\text{ \ \ }\left( x\in G_m\right).
\end{equation*}
It is well-known that if $f\in L_{p}\left( G_m\right) ,$ $g\in L_{1}\left(  G_m\right) $ and $1\leq p<\infty ,$ then $f\ast g\in L_{p}\left(  G_m\right) $ and
\begin{equation} \label{covstrong}
\left\Vert f\ast g\right\Vert_{p}\leq \left\Vert f\right\Vert_{p}\left\Vert g\right\Vert_{1}.
\end{equation}

In classical Fourier analysis  a point $x\in (-\infty,\infty)$ is called a Lebesgue point of integrable function $f$ is defined by
\begin{eqnarray*}
	\lim_{h\to 0}\frac{1}{h}\int_{x}^{x+h}\left\vert f(t)-f\left( x\right)\right\vert d\mu(t)=0
\end{eqnarray*}

On $G_m$ we have the following definition of Lebesgue point:
A point $x$ on the Vilenkin group is called Lebesgue point of   $f\in L_{1}\left( G_m\right),$ if
\[\lim_{n\rightarrow \infty }M_n\int_{I_{n}(x)}f\left(t\right) dt=f\left( x\right) \ \ \ \ \ a.e.\ x\in G_m.\]

It is well-known that if  $f\in L_{1}\left( G_m\right)$ then a.e point is Lebesgue point and
\begin{equation} \label{ae}
\underset{n\rightarrow \infty }{\lim }S _{M_n}f(x)=f(x), \ \ \text{a.e. on} \ \ G_m.
\end{equation}
where  $S _{M_n}$ is $M_n$-th partial sums with respect to Vilenkin systems.

Let introduced the operator
\begin{eqnarray*}
	W_Af(x):= \sum_{s=0}^{A-1}M_s\sum_{r_s=1}^{m_s-1}\int_{I_A(x-r_se_s)}\left\vert f(t)-f\left( x\right)\right\vert d\mu(t)
\end{eqnarray*}

A point $x \in G_m $ is a Vilenkin-Lebesgue point of $f \in L_1 (G_m),$ if
\begin{eqnarray*}
	\lim_{A \rightarrow \infty} W_Af (x)=0.
\end{eqnarray*}

In most applications the a.e. convergence of $\left\{ T_{n}:n\in \mathbb{N}\right\} $ can be established for $f$ in some dense class of $L_{1}\left(G_{m}\right).$ In particular, the following result play an important role for studying this type of questions (see e.g. the books \cite{gol}, \cite{sws} and
\cite{13z}).

\begin{lemma}\label{lemmaae} 
	Let $f\in L_{1}$ and $T_{n}:L_{1}\rightarrow L_{1}$ be some
	sub-linear operators and
	\begin{equation*}
	T^{\ast }:=\sup_{n\in\mathbb{N}
	}\left\vert T_{n}\right\vert.
	\end{equation*}
	If
	\begin{equation*}
	T_{n}f\rightarrow f\text{ a.e. \ \ for every }f\in S
	\end{equation*}%
	where the set $S$ is dense in the space $\ L_{1}$ and the maximal operator $%
	T^{\ast }$ is bounded from the space $L_{1}$ to the space $weak-L_{1},$ that
	is%
	\begin{equation*}
	\sup_{\lambda >0}\lambda \mu \left\{ x\in G_{m}:\text{ }\left\vert T^{\ast
	}f\left( x\right) \right\vert >\lambda \right\} \leq \left\Vert f\right\Vert
	_{1},
	\end{equation*}%
	then
	\begin{equation*}
	T_{n}f\rightarrow f,\text{ a.e. \ \ for every }f\in L_{1}\left( G_{m}\right)
	.
	\end{equation*}
\end{lemma}

\begin{remark}
	\label{remarkae}Since the Vilenkin function $\psi _{m}$ is constant on $%
	I_{n}(x)$ for every $x\in G_{m}$ and $0\leq m<M_{n},$ it is clear that each
	Vilenkin function is a complex-valued step function, that is, it is a finite
	linear combination of the characteristic functions%
	\begin{equation*}
	\chi \left( E\right) =\left\{
	\begin{array}{ll}
	1, & x\in E, \\
	0, & x\notin E.%
	\end{array}%
	\right.
	\end{equation*}%
	On the other hand, notice that, by \eqref{9dn}, it
	yields that%
	\begin{equation*}
	\chi \left( I_{n}(t)\right) \left( x\right) =\frac{1}{M_{n}}%
	\sum_{j=0}^{M_{n}-1}\psi _{j}\left( x-t\right) ,\text{ \ \ }x\in I_{n}(t),
	\end{equation*}%
	for each $x,t\in $ $G_{m}$ and $n\in \mathbb{N}$. Thus each step function is
	a Vilenkin polynomial. Consequently, we obtain that the collection of step functions
	coincides with a collection of Vilenkin polynomials $\mathcal{P}$. Since the
	Lebesgue measure is regular it follows from Lusin theorem that given $%
	f\in L_{1}$ there exist Vilenkin polynomials $P_{1},P_{2}...,$ such that $%
	P_{n}\rightarrow f$ \ a.e. when $n\rightarrow \infty .$ This means that the
	Vilenkin polynomials are dense in the space $L_{1}.$
\end{remark}

\subsection{Basic notations concerning Walsh groups and Functions}

\ \ \ Let us define by  $Q_2$	the set of rational numbers of the form $p 2^{-n}$, where $0\leq p\leq 2^{n}-1$ for some $p\in\mathbb{N}$ and $n\in\mathbb{N}$.

Any $x \in [0,1]$ can be written in the form

\begin{equation}
\label{18} x=\sum_{k=0}^{\infty} x_k 2^{-(k+1)},
\end{equation}
where each $x_k = 0$ or $1$. For each $x \in [0,1]\setminus Q_2$ there is only one expression of this form. We shall call it the dyadic expansion of $x.$ When $x \in Q_2$ there are two expressions of this form, one which terminates in $0 \text{'s}$ and one which terminates in $1 \text{'s}.$ By the dyadic expansion of an $x \in Q_2$ we shall mean the one which terminates in $0 \text{'s}.$ Notice that $1\nsim Q_2$ so the dyadic expansion of $x = 1$ terminates in $1 \text{'s}.$

If $m_k= 2, \ \text{for all} \ k\in \mathbb{N},$ we have dyadic group $$G_2=\prod_{j=0}^{\infty}Z_{2},$$ which is called the Walsh group 

Let $G_0^* := \{x \in G_2: x = y^* \ \text{for some} \ y \in G_0\}.$ Define Fine's map  $\varrho : [0,1] \rightarrow G_2 $ by

\begin{equation}
\label{19} \varrho (x):=(x_0, x_1, \dots)
\end{equation}
where $x$ has dyadic expansion \eqref{18}. Then $\varrho$ is a strictly increasing, $1-1$ map from $[0,1)$ onto $G \setminus G_0^*$. Moreover, it is easy to prove that $\varrho$ satisfies
\begin{equation}
\label{20} \left\{ \begin{array}{l}\varrho (x+)=\varrho (x-)=\varrho (x) \ \ \ \ \ \ \ \ \ \ \ \ x\in (0,1)\setminus Q \\
\varrho (x+) \varrho (x) , \varrho (x-)=\varrho^* (x) \ \ \ \ \ \ \ \ x\in Q \\
\varrho  (0+)=0, \varrho (1-)=0^* . \end{array}\right.
\end{equation}
Here $\varrho (x+)$ represents the limit of $\varrho$ at $x$ from the right in the usual topology on $[0,1]$ and $\varrho (x-)$ that from the left.

Let $C(G_2)$ represent the collection of functions $f : G_2 \rightarrow \mathbb{R}$ which are continuous from the dyadic topology on $G_2$ to the usual topology on $\mathbb{R}.$ Let $C_w$ represent the collection of functions $g : [0,1) \rightarrow \mathbb{R}$ which are continuous at every dyadic irrational, continuous from the right on $[0,1)$ and have a finite limit from the left on $(0,1]$ all this in the usual topology.

We shall call the map $f \rightarrow f \circ \varrho $ the canonical isomorphism. It is easy to see that this map is a vector space isomorphism from $C(G_2)$ onto $C_w.$ First, it is clear by \eqref{20} that if $f \in C(G_2)$ and $g := f\circ \varrho$ then
\begin{equation}
\label{21} \left\{ \begin{array}{l}g (x+)=g (x-)=g (x) \ \ \ \ \ \ \ \ \ \ \ \ \ \ \ \ \ x\in (0,1)\setminus Q_2 \\
g (x+) g (x) , g (x-)=f(\varrho^* (x)) \ \ \ \ \ \ \ \ \ x\in Q_2 \\
g  (0+)=f(0), g (1-)=f(0^*) . \end{array}\right.
\end{equation}
Thus the canonical isomorphism takes $C(G_2)$ into $C_w.$ Next, notice by construction that the canonical isomorphism is a vector space homomorphism, i.e., preserves function addition and scalar multiplication. Finally, if $g\in C_w$ then the map $f$ defined on $G_2$ by
\begin{equation}
\label{22} \left\{ \begin{array}{l}f(y):=g (x) \ \ \ \ \ \ \ \ \ y=\varrho(x),  \ \ x\in (0,1)\setminus Q_2 \\
f(y^*):=g (x-) \ \ \ \ \ y=\varrho(x), \ \ \ \ x\in Q_2 \\
f  (0^*)=g(1-) \end{array}\right.
\end{equation}
is continuous on $G_2.$

It is also easy to see that the canonical isomorphism takes the character system $\hat{G_2}$ onto the Walsh system $w.$ Indeed if $x\in [0,1)$ has dyadic expansion \eqref{18} then the definition of the Rademacher functions given in 1.1 can be rewritten as
$$\upsilon_n(x)=(-1)^{x_n} \ \ \ \ \ \ (n\in \mathbb{N})$$

\begin{equation}
\label{11} \rho_n(x):=(-1)^{x_n}
\end{equation}

Comparing this with \eqref{11}, we see that $\upsilon_n = \rho_n \circ \varrho $ on $[0,1)$ and $\rho_n(x)=v_n(\left|x\right|)$ for $x \in G_2 \setminus G_0^*$ and every $n \in \mathbb{N}.$ 






\begin{equation}\label{4mag} 
\phi_n :=\prod_{k=0}^{\infty}\upsilon_k ^{n_k}={\upsilon_{\left\vert n\right\vert }}\left( x\right) \left( -1\right) ^{\underset{k=0}{\overset{\left\vert n\right\vert -1}{\sum }}n_{k}x_{k}}\text{\qquad }\left( n\in \mathbb{N}\right),
\end{equation}

\begin{equation}
\label{12mag} w_n:=\prod_{k=0}^{\infty}\rho_k ^{n_k}
\end{equation}

It follows from  \eqref{4mag} and \eqref{12mag} that
$$\phi_n = w_n \circ \varrho $$
and
$$ w_n(x)=\phi_n(|x|) \ \ \ \ \ \ (x\in G_2\setminus G_0^*) $$
for every $n\in \mathbb{N}.$






Fine's map can be used to define a new addition and a new topology on $[0,1)$ which are closely related to those on $G.$ Indeed, define the dyadic sum of two numbers $x, y \in [0,1)$ and the dyadic distance between these numbers by
$$x\overset{\cdot}+ y=|\varrho(x)+\varrho(y)|$$

In terms of the dyadic expansions of $x$ and $y$ we have

$$x\overset{\cdot }+y=\sum\limits_{k=0}^{\infty }{{}}\left| {{x}_{k}}-{{y}_{k}} \right|{{2}^{-k-1}}$$

Hence $\overset{\cdot}+$ is evidently a metric and a commutative binary operation on $[0,1)$ which satisfies $x \overset{\cdot}+ x = 0.$ We shall call the topology generated by $\overset{\cdot}+$ on $[0, 1)$ the dyadic topology. Note, $[0,1)$ is not a group under $\overset{\cdot}+.$

The Walsh functions almost behave like characters with respect to dyadic addition, namely,
\begin{equation}
\label{23} \phi_n(x\overset{\cdot }+y) = \phi_n(x)\phi_n(y) \ \ \ \ \ (n\in \mathbb{N},\ \ x,y\in [0,1), x\overset{\cdot }+y \notin Q_2).
\end{equation}
To prove \eqref{23} fix $m \in \mathbb{N}$ and $x, y \in [0,1).$ Notice that
$$\upsilon_m(x)\upsilon_m(y) = \rho_m \circ \varrho(x)\rho_m\circ\varrho(y) = \rho_m(\varrho(x)+\varrho(y)) ,$$
and that
$$\upsilon_m(x\overset{\cdot }+y)=\rho_m \left(\varrho \left(\left|\varrho(x)+\varrho(y)\right|\right) \right).$$
Since $\left|\varrho(x)+\varrho(y)\right|$  is a dyadic irrational when $x\overset{\cdot }+y$ is, it is clear that
$$\varrho \left(\left|\varrho(x)+\varrho(y)\right|\right)$$
for $x\overset{\cdot }+y \notin Q_2.$ Consequently, \eqref{23} holds for the Rademacher case, i.e., for $n = 2^m.$ But the general case follows immediately since Walsh functions are finite products of Rademacher functions. Since for each fixed $y \in [0,1)$ the set of points $x$ which satisfy $x\overset{\cdot }+y \in Q_2$ is a countable set, we observe that \eqref{23} holds for a.e. $x, y \in [0,1).$

By a dyadic interval in $[0,1)$ we shall always mean an interval of the form

\begin{equation}
\label{24} I(p,n):=\left[p 2^{-n}, (p+1)  2^{-n}, \right) \ \ \ \ \ \ \ (0\leq p<2^n, n,p, \in \mathbb{N}).
\end{equation}
Clearly, the dyadic topology is generated by the collection of dyadic intervals. Moreover, each dyadic interval is both open and closed in the dyadic topology. It follows that each Walsh function is continuous in the dyadic topology. Thus the dyadic topology differs from the usual topology in an essential way. 

For each $x \in [0,1)$ and $n \in \mathbb{N}$ we shell denote the dyadic interval of length $2^{-n}$ which contains $x$ by $I_n (x).$ Thus
$$I_n(x):=I(p,n)$$
where $0\leq p < 2^n$ is uniquely determined by the relationship $x \in I(p, n).$ This is the same notation used for dyadic intervals on the group but will not cause problems because context will make it clear whether the dyadic interval is in the group or inside the unit interval.

A function $I : [0,1) \rightarrow R$ which is continuous from the dyadic topology to the usual topology will be called $W-continuou$. Since
$$\left|x-y\right|\leq x\overset{\cdot }+y \ \ \ (x,y \in [0,1)),$$
it is clear that every classically continuous function on $[0,1)$ is $W- continuous.$ In fact, every function in $C_w$ is unifoflIuy $W-continuous$ on the unit interval. On the other hand, not every $W-continuous$ function belongs to $C_w.$

Let $L^0$ represent the collection of a.e. finite, Lebesgue measurable functions from $[0,1)$ into $[-\infty,\infty]$. For $0 < p < \infty$ let $L^p$ represent the collection of $f\in L^0$ for which
$$\|f\|_p := \left(\int_{0}^{1}\left|f\right|^p \right)^{1/p}$$
is finite. Let $L^{\infty}$ represent the collection of $f \in L^0$ for which
$$ \|f\|_\infty := \inf \{ y\in \mathbb{R} :\left|f(x)\right|\leq y \ \text{for a.e. } \ \ x\in [0,1) \}$$
is finite. It is well known that $L^p$ is a Banach space for each $1 \leq p \leq \infty.$

If $f\in L_{1}\left(G_2\right),$  then we can establish the Fourier
coefficients, the partial sums of the Fourier series, the Fejér means, the
Dirichlet and Fejér kernels with respect to the Vilenkin system $\psi $ (Walsh system $w$) in
the usual manner:%
\begin{eqnarray*}
	\widehat{f}^{w }\left( k\right) &:&=\int_{E}f\overline{w }_kd\mu,\text{\ \ \ \ \ \ }  \qquad \left(k\in \mathbb{N}\right) , \\
	S^w_n f&:&=\sum_{k=0}^{n-1}\widehat{f}\left(k\right)w_k, \qquad  \text{ \ \
		\ }\left(n\in \mathbb{N}_{+}, \ S^{\alpha}_0 f:=0\right) , \\
	D^{w}_n &:&=\sum_{k=0}^{n-1}w_{k},\text{ \qquad\ \ \ \ \ }    \left(n\in \mathbb{N}_{+}\right), \\
\end{eqnarray*}

We state well-known equalities for Dirichlet kernels (for details see e.g. \cite{gol} and \cite{sws}):

\begin{equation} \label{1dn}
D^w_{2^{n}}\left( x\right)=\left\{ 
\begin{array}{ll}
2^{n}, & \,\text{if\thinspace \thinspace \thinspace }x\in I_{n} \\ 
0, & \text{if}\ \ \ x\notin I_{n}
\end{array}
\right.  
\end{equation}
and
\begin{equation}\label{2dn}
D^w_n=w_n\overset{\infty }{\underset{k=0}{\sum }}n_kr_kD^w_{2^k}=w_n
\overset{\infty }{\underset{k=0}{\sum }}n_{k}\left(
D^w_{2^{k+1}}-D^w_{2^{k}}\right),\ \ \text{ for \ }\ \ n=\overset{\infty }{\underset{i=0}{\sum }}n_{i}2^{i},
\end{equation}







The most properties of
Lebesgue constants with respect to the Walsh-Paley system were obtained
by Fine in \cite{fi}. In (\cite{sws}, p. 34), the two-sided estimate
$$\frac{V(n)}{8}\leq L_n \leq V(n)$$
is proved, where
$n=\sum_{j=1}^{\infty}n_j2^j$  and $V\left(n\right)$ is defined by
$$V\left(n\right):=\sum_{j=1}^{\infty}\left\vert n_{j+1}-n_j\right\vert+n_0.$$

If $f\in L_{1}\left(G_2\right),$  then  Fejér means and
Fejér kernels with respect to the Walsh system $\psi $ (Walsh system $w$) is defined by
\begin{eqnarray*}
	\sigma^{w}_n f &:&=\frac{1}{n}\sum_{k=0}^{n-1}S^{w}_k f,\text{ \ \ \ \ \ \ }  \left(n\in \mathbb{N}_{+}\right) , \\
	K^{w}_n &:&=\frac{1}{n}\overset{n-1}{\underset{k=0}{\sum }}D^{w}_{k}, \ \ \ \ \ \ \ \ \ \left(n\in \mathbb{N}_{+}\right) .
\end{eqnarray*}	

The $n$-th N\"orlund logarithmic mean \index{S}{N\"orlund logarithmic mean} $L^{\alpha}_n$ and the Riesz logarithmic mean \index{S}{Riesz logarithmic mean} $%
R^{\alpha}_{n}$ with respect to the Walsh system $\psi $ (Walsh system $w$) are defined by  \index{N}{$L_nf$} \index{N}{$R_nf$}
\begin{eqnarray*}
	L^{w}_nf&:=&\frac{1}{l_n}\sum_{k=1}^{n-1}\frac{S^w_kf}{n-k}, \text{ \ \ \ \ \ \ } \left(n\in \mathbb{N}_{+}\right) ,\\
	R^{w}_nf&:=&\frac{1}{l_n}\sum_{k=1}^{n-1}\frac{S^w_kf}{k}, \text{ \ \ \ \ \ \ }  \left(n\in \mathbb{N}_{+}\right) ,
\end{eqnarray*}
respectively, where \index{N}{$l_n$}
\begin{equation*}
l_n:=\sum_{k=1}^{n-1}\frac{1}{k}.
\end{equation*}

Kernels of N\"orlund logarithmic mean and Riesz logarithmic mean are  respectively defined by
\index{N}{$P_nf$} \index{N}{$Y_nf$}
\begin{eqnarray*}
	P^{w}_nf&:=&\frac{1}{l_n}\sum_{k=1}^{n-1}\frac{D^{w}_kf}{n-k}, \text{ \ \ \ \ \ \ }  \left(n\in \mathbb{N}_{+}\right) ,\\
	Y^{w}_nf&:=&\frac{1}{l_n}\sum_{k=1}^{n-1}\frac{D^{w}_kf}{k}, \text{ \ \ \ \ \ \ }  \left(n\in \mathbb{N}_{+}\right).
\end{eqnarray*}

\subsection{Theory of Martingale Hardy spaces for $0<p\leq 1$}

The $\sigma $-algebra \index{S}{$\sigma $-algebra} generated by the intervals
\begin{equation*}
\left\{I_n\left(x\right):x\in G_m\right\}
\end{equation*}
will be denoted by
$\digamma_n\left(n\in\mathbb{N}\right).$  \index{N}{$\digamma_n\left(n\in\mathbb{N}\right)$}

A sequence $f=\left(f^{\left(n\right)}:n\in\mathbb{N}\right)$
\index{N}{$f=\left(f^{\left(n\right)}:n\in\mathbb{N}\right)$} of integrable functions $f^{\left(n\right)}$ is said to be a martingale \index{S}{martingale} with
respect to the $\sigma $-algebras
$\digamma_n\left(n\in\mathbb{N}\right)$ if (for details see e.g. Weisz \cite{We1})

$1) \qquad f_n$ is $\digamma _{n}$\ measurable for all $n\in \mathbb{N},$

$2) \qquad S_{M_n}f^{(m)}=f^{(n)}$ for all $n\leq m.$

The martingale $f=\left(f^{\left(n\right)},n\in \mathbb{N}\right)$ is said to be $L_p$-bounded ($0<p\leq\infty$) if $f^{\left(n\right)}\in L_p$ and

\begin{equation*}
\left\Vert f\right\Vert_p:=\underset{n\in\mathbb{N}}{\sup}\left\Vert
f_n\right\Vert_p<\infty.
\end{equation*}

If $f\in L_1\left(G_m\right),$ then it is easy to show that the
sequence $F=\left(E_nf:n\in \mathbb{N}\right)$ is a martingale. This
type of martingales is called regular. If $1\leq p\leq \infty $ and $f\in L_p\left(G_m\right)$ then $f=\left(f^{\left(n\right)},n\in \mathbb{N}\right) $ is $L_p$-bounded and
\begin{equation*}
\underset{n\rightarrow\infty}{\lim}\left\Vert E_nf-f\right\Vert_p=0
\end{equation*}
and consequently
$\left\Vert F\right\Vert_p=\left\Vert f\right\Vert_p$
(see \cite{W}). The converse of the latest statement holds also if $%
1<p\leq\infty $ (see \cite{W}): for an arbitrary $L_p$-bounded martingale $f=\left(f^{\left(n\right)},n\in \mathbb{N}\right)$ there exists a function $f\in L_p\left( G_m\right)$ for which $f^{\left( n\right)}=E_nf.$ If $p=1,$ then there exists a function $f\in L_1\left(
G_m\right)$ of the preceding type if and only if $f$ is uniformly
integrable (see \cite{W}), namely, if
\begin{equation*}
\underset{y\rightarrow\infty}{\lim}\underset{n\in \mathbb{N}}{\sup}
\int_{\left\{\left\vert f_n\right\vert>y\right\}}\left\vert f_n\left(
x\right)\right\vert d\mu\left(x\right)=0.
\end{equation*}

Thus the map $f\rightarrow f:=$ $\left(E_nf:n\in \mathbb{N}\right) $ is
isometric from $L_p$ onto the space of $L_p$-bounded martingales when $
1<p\leq\infty.$ Consequently, these two spaces can be identified with each other. Similarly, the space $L_{1}\left( G_{m}\right) $ can be identified with the space of uniformly integrable martingales.

Analogously, the martingale
$f=\left(f^{\left(n\right)},n\in\mathbb{N}\right) $ is said to be $weak-L_p$-bounded ($0<p\leq\infty$) if $f^{\left(n\right)}\in L_p$ and
\begin{equation*}
\left\Vert f\right\Vert_{weak-L_p}:=\underset{n\in \mathbb{N}}{\sup }
\left\Vert f_n\right\Vert_{weak-L_p}<\infty.
\end{equation*}

The maximal function of a martingale  \index{S}{maximal function of a martingale} $f$  is defined by
\index{N}{$f^{\ast}$}
\begin{equation*}
f^{\ast}:=\sup_{n\in\mathbb{N}}\left\vert f^{(n)}\right\vert.
\end{equation*}

In the case $f\in L_1(G_m),$ the maximal functions are also given by
\begin{equation*}
f^{\ast}\left(x\right):=\sup_{n\in\mathbb{N}}\frac{1}{\left\vert
	I_n\left(x\right)\right\vert}\left\vert\int_{I_n\left(x\right)} f\left(u\right)d\mu\left(u\right)\right\vert.
\end{equation*}

For $0<p<\infty $ the Hardy martingale spaces \index{S}{Hardy martingale spaces} $H_p$ \index{N}{$H_p$} consist of all martingales for which
\begin{equation*}
\left\Vert f\right\Vert_{H_p}:=\left\Vert f^{\ast}\right\Vert_p<\infty .
\end{equation*}

Vilenkin-Fourier coefficients of martingale
\index{S}{Vilenkin-Fourier coefficients of martingale} $f=\left(f^{\left(n\right)}:n\in\mathbb{N}\right)$ must be defined in a slightly different manner:
\begin{equation*}
\widehat{f}\left(i\right):=\lim_{k\rightarrow\infty}\int_{G_m}f^{\left(k\right)}\overline{\psi}_id\mu.
\end{equation*}

Investigation of the classical Fourier analysis, definition of several variable Hardy spaces and real Hardy spaces and related theorems of atomic decompositions of these spaces can be found in Fefferman and Stein \cite{FS} (see also Later \cite{La}, Torchinsky \cite{Tor1}, Wilson \cite{Wil}).

A bounded measurable function $a$ \index{N}{$a$} is a p-atom \index{S}{p-atom} if there exist an interval $I$
such that \qquad
\begin{equation*}
\int_{I}ad\mu =0,\text{ \ \ }\left\Vert a\right\Vert _{\infty }\leq \mu
\left( I\right) ^{-1/p},\text{ \ \ supp}\left( a\right) \subset I.
\end{equation*}

Explicit constructions of $p$-atoms can be found in the papers \cite{BGG} and \cite{BGG2} by Blahota, G\'{a}t and Goginava.

Next, we note that the Hardy martingale spaces $H_{p}\left( G_{m}\right) $
for $0<p\leq 1$ have atomic characterizations \index{S}{atomic characterizations}:

The proof of Lemma \ref{lemma2.1} was proved by Weisz \cite{We1,We3}.
\begin{lemma}\label{lemma2.1}
	A martingale $f=\left( f^{\left( n\right) }:n\in \mathbb{N}%
	\right) $ is in $H_{p}\left( 0<p\leq 1\right) $ if and only if there exist a
	sequence $\left( a_{k},k\in \mathbb{N}\right) $ of p-atoms and a sequence $%
	\left( \mu _{k}:k\in \mathbb{N}\right) $ of real numbers such that, for
	every $n\in \mathbb{N},$%
	\begin{equation}
	\qquad \sum_{k=0}^{\infty }\mu _{k}S_{M_{n}}a_{k}=f^{\left( n\right) },\text{
		\ \ a.e.,}  \label{condmart}
	\end{equation}%
	where
	\begin{equation*}
	\qquad \sum_{k=0}^{\infty }\left\vert \mu _{k}\right\vert ^{p}<\infty .
	\end{equation*}%
	Moreover,
	\begin{equation*}
	\left\Vert f\right\Vert _{H_{p}}\backsim \inf \left( \sum_{k=0}^{\infty
	}\left\vert \mu _{k}\right\vert ^{p}\right) ^{1/p},
	\end{equation*}%
	where the infimum is taken over all decomposition of $f=\left( f^{\left(
		n\right) }:n\in \mathbb{N}\right) $ of the form (\ref{condmart}).
\end{lemma}

Explicit constructions of $H_p$ martingales  can be found in the papers \cite{BPTW}, \cite{ptw20}, \cite{tep91}, \cite{tep20},  \cite{tep18}, \cite{tep19},  \cite{tep13},  \cite{tepPhDGeo}, \cite{tep16}, \cite{tep21}, \cite{ptc}.

By using atomic characterization it can be easily proved that the following Lemmas hold:

\begin{lemma}\label{lemma2.2} 
	Suppose that an operator $T$ is sub-linear and for some $0<p_0\leq 1$
	\begin{equation*}
	\int\limits_{\overset{-}{I}}\left\vert Ta\right\vert ^{p_0}d\mu \leq
	c_{p}<\infty
	\end{equation*}%
	for every $p_0$-atom  $a$, where $I$
	denotes the support of the atom. If $T$ is bounded from $
	L_{p_1}$  to $L_{p_1},$ $\left(1<p_{1}\leq \infty
	\right) $
	then
	\begin{equation} \label{p_0}
	\left\Vert Tf\right\Vert _{p_0}\leq c_{p_0}\left\Vert f\right\Vert _{H_{p_0}}.
	\end{equation}%
	Moreover, if $p_0<1,$ then we have the weak (1,1) type estimate
	\begin{equation*}
	\lambda \mu \left\{ x\in G_{m}:\text{ }\left\vert Tf\left( x\right)
	\right\vert >\lambda \right\} \leq \left\Vert f\right\Vert _{1}
	\end{equation*}%
	for all $f\in L_{1}.$
	
	The proof of Lemma \ref{lemma2.2} can be found in Weisz \cite{We3}.
\end{lemma}

\begin{lemma}
	\label{lemma2.3} Suppose that an operator $T$ is sub-linear and for some $%
	0<p_0\leq 1$%
	\begin{equation*}
	\underset{\lambda >0}{\sup }\lambda ^{p_0}\mu \left\{ x\in \overset{-}{I}:\left\vert Tf\right\vert >\lambda \right\} \leq c_{p_0}<+\infty
	\end{equation*}%
	for every $p_0$-atom $a$, where $I$  denote
	the support of the atom. If $T$ is bounded from
	$L_{p_1}$  to $L_{p_1},$ $\left(1<p_{1}\leq \infty\right) $
	then
	\begin{equation*}
	\left\Vert Tf\right\Vert _{weak-L_{p_0}}\leq c_{p_0}\left\Vert f\right\Vert
	_{H_{p_0}}.
	\end{equation*}%
	Moreover, if $p_0<1,$ then
	\begin{equation*}
	\lambda \mu \left\{ x\in G_{m}:\text{ }\left\vert Tf\left( x\right)
	\right\vert >\lambda \right\} \leq \left\Vert f\right\Vert _{1},
	\end{equation*}%
	for all $f\in L_{1}.$
\end{lemma}

The concept of modulus of continuity in martingale Hardy space \index{S}{ martingale Hardy space} $H_{p}$ $\left( p>0\right) $
is defined by \index{N}{$\omega_{H_p}\left(\frac{1}{M_n},f\right)$}
\begin{equation*}
\omega_{H_p}\left(\frac{1}{M_n},f\right):=\left\Vert
f-S_{M_n}f\right\Vert _{H_p}.
\end{equation*}

We need to understand the meaning of the expression $f-S_{M_n}f$ where $f$ is a martingale and $S_{M_n}f$ is function. So, we give an explanation in the following remark:

\begin{remark}
	\label{lemma2.3.6} Let $0<p\leq 1.$ Since
	\begin{equation*}
	S_{M_{n}}f=f^{\left( n\right) },\text{ \ for }f=\left( f^{\left( n\right)
	}:n\in
	\mathbb{N}
	\right) \in H_{p}
	\end{equation*}%
	and
	\begin{equation*}
	\left( S_{M_{k}}f^{\left( n\right) }:k\in
	\mathbb{N}
	\right) =\left( S_{M_{k}}S_{M_{n}}f,k\in
	\mathbb{N}
	\right)
	\end{equation*}%
	\begin{equation*}
	=\left( S_{M_{0}}f,\ldots ,S_{M_{n-1}}f,S_{M_{n}}f,S_{M_{n}}f,\ldots \right)
	\end{equation*}%
	\begin{equation*}
	=\left( f^{\left( 0\right) },\ldots ,f^{\left( n-1\right) },f^{\left(
		n\right) },f^{\left( n\right) },\ldots \right)
	\end{equation*}%
	we obtain that%
	\begin{equation*}
	f-S_{M_{n}}f=\left( f^{\left( k\right) }-S_{M_{k}}f:k\in
	\mathbb{N}
	\right)
	\end{equation*}%
	is a martingale, for which%
	\begin{equation}
	\left( f-S_{M_{n}}f\right) ^{\left( k\right) }=\left\{
	\begin{array}{ll}
	0, & k=0,.\ldots ,n, \\
	f^{\left( k\right) }-f^{\left( n\right) }, & k\geq n+1,%
	\end{array}%
	\right.  \label{g100}
	\end{equation}
\end{remark}

Watari \cite{wat} showed that there are strong connections between
\begin{equation*}
\omega _{p}\left( \frac{1}{M_{n}},f\right) ,\text{ \ }E_{M_{n}}\left(
L_{p},f\right) \text{ \ \ and\ \ \ }\left\Vert f-S_{M_{n}}f\right\Vert
_{p},\ \ p\geq 1,\text{ \ }n\in \mathbb{N}.
\end{equation*}

In particular,%
\begin{equation}
\frac{1}{2}\omega _{p}\left( \frac{1}{M_{n}},f\right) \leq \left\Vert
f-S_{M_{n}}f\right\Vert _{p}\leq \omega _{p}\left( \frac{1}{M_{n}},f\right)
\label{eqvi}
\end{equation}%
and%
\begin{equation*}
\frac{1}{2}\left\Vert f-S_{M_{n}}f\right\Vert _{p}\leq E_{M_{n}}\left(
L_{p},f\right) \leq \left\Vert f-S_{M_{n}}f\right\Vert _{p}.
\end{equation*}

\begin{remark}
	\label{lemma2.3.7}Since
	\begin{equation*}
	\left\Vert f\right\Vert _{H_{p}}\sim \left\Vert f\right\Vert _{p},
	\end{equation*}%
	when $p>1$, by applying (\ref{eqvi}), we obtain that%
	\begin{equation*}
	\omega _{H_{p}}\left( \frac{1}{M_{n}},f\right) \sim \omega _{p}\left( \frac{1%
	}{M_{n}},f\right) .
	\end{equation*}
\end{remark}

Next lemma gives answer what happens when $p>1$. The proof  can be found in Neveu \cite{W} (see also Weisz \cite{We5}).

\begin{lemma} \label{hpLp_eqvi}
	Let $p>1.$ Then
	\begin{equation*}
	H_{p}\sim L_{p}.
	\end{equation*}
\end{lemma}

The proof of Lemma \ref{lemma2.3.4} can be found in \cite{We3} (see also book \cite{sws}).
\begin{lemma}
	\label{lemma2.3.4}If $f\in L_{1},$ then the sequence $F:=\left(
	S_{M_{n}}f:n\in \mathbb{N}\right) $ is a martingale and
	\begin{equation*}
	\left\Vert F\right\Vert _{H_{p}}\sim \left\Vert \sup_{n\in \mathbb{N}%
	}\left\vert S_{M_{n}}f\right\vert \right\Vert _{p}.
	\end{equation*}%
	Moreover, if $F:=$ $\left( S_{M_{n}}f:n\in \mathbb{N}\right) $ is a regular
	martingale \index{S}{regular martingale} generated by $f\in L_{1},$ then \qquad \qquad \qquad \qquad
	\begin{equation*}
	\widehat{F}\left( k\right) =\int_{G_{m}}f\left( x\right) \psi_{k}\left(
	x\right)d\mu\left(x\right) =\widehat{f}\left(k\right),\text{ \qquad }k\in \mathbb{N}.
	\end{equation*}
\end{lemma}	

	\newpage

\section[Partial sums and Fej\'er means on $H_p$ spaces]{Partial sums and Fej\'er means of Vilenkin-Fourier series on Martingale Hardy spaces}

\vspace{0.5cm}

\subsection{Some classical results on the partial sums and and Fej\'er means of Vilenkin-Fourier series}

According to the  Riemann-Lebesgue lemma (for details see e.g. the book \cite{sws}) we obtain that $
\widehat{f}\left( k\right) \rightarrow 0,\text{ \ \ when \ \ \ }k\rightarrow
\infty ,
$
for each $f\in L_{1}.$

It is well-known (see e.g. the books \ \cite{AVD}\textit{ }and \cite{sws}
) that if $%
f\in L_{1}$ and the Vilenkin series
$
T\left( x\right) =\sum_{j=0}^{\infty }c_{j}\psi _{j}\left( x\right)
$
convergences to $f$ \ in $L_{1}$-norm, then
$
c_{j}=\int_{G_{m}}f\overline{\psi }_{j}d\mu :=\widehat{f}\left( j\right) ,
$
i.e. in this case the approximation series must be a Vilenkin-Fourier
series. An analogous result is true also if the Vilenkin series convergences
uniformly on $G_{m}$ to an integrable function $f$ .

\bigskip By using the Lebesgue constants we easily obtain that $S_{n_{k}}f$
\ convergence to $f$ \ in $L_{1}$-norm, for every integrable function $f$,
if and only if $\sup_{k}L_{n_{k}}\leq c<\infty .$ There are various results
when $p>1.$

It is also well-known that (see e.g. the book
\ \cite{sws})
\begin{equation*} \label{snp}
\left\Vert S_{n}f\right\Vert _{p}\leq c_{p}\left\Vert f\right\Vert _{p},%
\text{ \ when \ }p>1,
\end{equation*}%
but it can be proved also a more stronger result (see e.g. the book \ \cite{sws}):
\begin{equation*}
\left\Vert S^{\ast }f\right\Vert _{p}\leq c_{p}\left\Vert f\right\Vert _{p},%
\text{ \ when \ }f\in L_{p},\text{ \ }p>1.
\end{equation*}

Moreover, in the case $ p=1 $ Watari \cite{Wat1} (see also Gosselin \cite{goles} and Young \cite{Yo}) proved that there exists an absolute constant $ c $ such that, for $ n =1,2,..., $
\begin{eqnarray*}
	\lambda \mu
	\left( \vert S_n f\vert>\lambda \right)&\leq& c\left\Vert f\right\Vert_{1}, \ \ \ f\in L_1(G_m), \ \ \lambda>0.
\end{eqnarray*}

Uniform and point-wise convergence and some approximation properties of
partial sums with respect to Vilenkin (Walsh) and trigonometric systems in $L_{1}$ norms was investigated by Antonov \cite{Ant}, Avdispahi\'c and Memi\'c \cite{am}, \index{A}{Avdispahi\'c and Memi\'c}  \index{A}{Baramidze} Baramidze \cite{Bar1}, Goginava \index{A}{Goginava} \cite{gog1,gog2}, Shneider \index{A}{Shneider} \cite{Sne},   Sjölin \index{A}{\cite{sj}} \cite{sj}, Onneweer and Waterman \cite{OW1,OW2}. \index{A}{Onneweer and Waterman}
Fine \cite{fi} derived sufficient conditions for the uniform convergence,
which are in complete analogy with the Dini-Lipschitz \index{A}{Dini-Lipschitz} conditions. Guli\'cev
\cite{9} estimated the rate of uniform convergence of a Walsh-Fourier series by using Lebesgue \index{A}{Lebesgue} constants and modulus of continuity. Uniform convergence of
subsequences of partial sums was investigated also in Goginava and Tkebuchava \index{A}{Tkebuchava } \cite{gt}, Fridli \index{A}{Fridli} \cite{4} and G\'at \index{A}{G\'at}  \cite{5}. Approximation properties of the two-dimensional partial sums with respect to Vilenkin and trigonometric systems can be found \cite{sws} and \cite{13z}. 

A.e. convergence of Vilenkin-fourier series can be found in \cite{PSTW}. Divergence of Vilenkin Fourier series on the sets of measure zero and a.e can be found in Bitsadze \index{A}{Bitsadze} \cite{Bit1,Bit2},  Bugadze \index{A}{Bugadze} \cite{Bug1,Bug2}, Fej\'{e}r \index{A}{Fej\'{e}r} \cite{Fejer} Gosselin, \cite{goles} \index{A}{Gosselin}  Kahane \cite{Kahane2}, Katznelson \index{A}{Katznelson} \cite{katz}, Karagulian \cite{Kar1,Kar2}, Kheladze \index{A}{Kheladze} \cite{Khe1,Khe2},  Lebesgue \index{A}{Lebesgue} \cite{Leb1}   Stechkin \index{A}{Stechkin} \cite{Ste} Young \cite{Young5,Young6,Young7} \index{A}{Young} and Zhizhiashvili \cite{Zh4}.

Some estimates of Fourier coefficients and an absolute, point-wise and a.e. convergence and divergence of Fourier Series with respect to complete orthonormal systems were studied in Bochkarev \cite{boch2}, Gogoladze and Tsagareishvili \cite{GTs1,GTs2,GTs3}, Kashin and Saakyan \cite{KS}, Oniani \cite{On2}, Tsagareishvili and Tutberidze \cite{tsatu1}.   Approximation of functions on locally compact Abelian groups was investigate by Ugulava \cite{Dug1,Dug2} (see also \cite{CKT}).

Since $H_{1}\subset L_{1}$, according to Riemann-Lebesgue theorem, it yields that $\widehat{f}\left( k\right)\rightarrow 0$ when $k\rightarrow \infty ,$ for every $f\in H_{1}.$ 
The classical inequality of Hardy type is well known in the trigonometric as well as in
the Vilenkin-Fourier analysis and was proved in the trigonometric case by Hardy and
Littlewood \cite{hl} (see also the book \cite{cw}) and for the
Walsh system it was proved in the book \cite{sws}. Some inequalities relative to Vilenkin-Fourier coefficients were considered in \cite{11}, \cite{Si7},  \cite{SW0,SW5}, \cite{tep4}, \cite{We6,We4,We1}.

It is known (for details see e.g. the books \ \cite{sws}
and \cite{We1}) that the subsequence $S_{M_{n}}$ of the partial
sums is bounded from the martingale Hardy space $H_{p}$ to the Lebesgue
space $L_{p},$ for all $p>0.$ However, (see Tephnadze \cite{tep7}) there
exists a martingale $f\in H_{p}$ $\left( 0<p<1\right) ,$ such that
$$
\underset{n\in \mathbb{N}}{\sup }\left\Vert S_{M_{n}+1}f\right\Vert
_{weak-L_{p}}=\infty .
$$
The reason of the divergence of $S_{M_{n}+1}f$ \ is that when $0<p<1$ the
Fourier coefficients of $f\in H_{p}$ are not uniformly bounded (see
Tephnadze \cite{tep6}). On the other hand, there exists an
absolute constant $c_{p},$ depending only on $p,$ such that%
\begin{equation}\label{smn}
\left\Vert S_{M_n}f\right\Vert _{p}\leq c_{p}\left\Vert f\right\Vert _{H_{p}},\text{ \ } p>0,\text{ \ }n\in \mathbb{N}_{+}.
\end{equation}

Tephnadze \cite{tep7} (see also \cite{tep9} and \cite{tep12}) proved that for every $0<p<1,$  the  maximal operator
\begin{equation*}
\widetilde{S}_{p}^{\ast }f:=\sup_{n\in \mathbb{N}}\frac{\left\vert
	S_{n}f\right\vert }{\left( n+1\right) ^{1/p-1}}
\end{equation*}%
is bounded from the Hardy space $H_{p}$\ to the Lebesgue space $L_{p}.$
Moreover, the rate of the sequence $\left( n+1\right) ^{1/p-1}$ is in the sense sharp.

It follow that for any $0<p<1$ and $f\in H_{p},$  there exists an
absolute constant $c_{p},$ depending only on $p,$ such that%
\begin{equation*}
\left\Vert S_{n}f\right\Vert _{p}\leq c_{p}\left( n+1\right)
^{1/p-1}\left\Vert f\right\Vert _{H_{p}},\text{ \ }n\in \mathbb{N}_{+}.
\end{equation*}

Blahota, Persson, Nagy and Tephnadze (\cite{bnpt}) proved that for any $0<p\leq 1$ and a sub-sequence of	positive numbers $\left\{ \alpha _{k}:k\in \mathbb{N}\right\} $, satisfying the condition
\begin{equation}
\sup_{k\in \mathbb{N}}\rho \left( \alpha _{k}\right) =\varkappa <\infty ,  \label{ak0}
\end{equation}%
the maximal operator
$
\widetilde{S}^{\ast ,\vartriangle }f:=\sup_{k\in \mathbb{N}}\left\vert S_{\alpha
	_{k}}f\right\vert
$
is bounded from the Hardy space $H_{p}$ to the space $L_{p}.$
Moreover, for every $0<p<1\ \ $and any sub-sequence of positive numbers $\left\{ \alpha _{k}:k\in \mathbb{N}\right\} $ satisfying the condition
\begin{equation}
\sup_{k\in \mathbb{N}}\rho \left( \alpha _{k}\right) =\infty,  \label{ak}
\end{equation}%
there exists a martingale $f\in H_{p}, \ \left( 0<p<1\right) $ such that
$
\sup_{k\in \mathbb{N}}\left\Vert S_{\alpha _{k}}f\right\Vert
_{weak-L_{p}}=\infty .
$

It follow that for any $p>0$ and $f\in H_{p}$, the maximal operators
\begin{equation}
\widetilde{S}_{\#}^{\ast }f:=\sup_{n\in \mathbb{N}}\left\vert S_{M_{n}}f\right\vert   \ \
\text{ and } \ \ 	\sup_{n\in \mathbb{N_{+}}}\left\vert S_{M_{n}+M_{n-1}}f\right\vert
\end{equation}
are bounded from the Hardy space $H_{p}$ to the space $L_{p}.$
We also obtain that if $p>0$ and $f\in H_{p},$ the maximal operator
$$
\sup_{n\in \mathbb{N_{+}}}\left\vert S_{M_{n}+1}f\right\vert
$$
is not bounded from the Hardy space $H_{p}$ to the space $L_{p}.$

Tephnadze \cite{tep7} (see also \cite{tep9} and \cite{tep12}) proved that for any  $0<p<1$ and $f\in H_{p}.$ Then there exists an
absolute constant $c_{p}$ depending only on $p$ such that
\begin{equation*}
\left\Vert S_{n}f\right\Vert _{H_{p}}\leq c_{p}n^{1/p-1}\left\Vert
f\right\Vert _{H_{p}}.
\end{equation*}

Tephnadze (\cite{tep7}) proved that for any $0<p<1,$ $f\in H_{p}$ and $M_{k}<n\leq M_{k+1}$. Then
there is an absolute constant $c_{p}$ depending only on $p$ such that
\begin{equation*}
\left\Vert S_{n}f -f\right\Vert _{H_{p}}\leq c_{p}n^{1/p-1}\omega _{H_{p}}\left(
\frac{1}{M_{k}},f\right) .
\end{equation*}

This estimate immediately follows that if $0<p<1,$ $f\in H_{p}$ and%
\begin{equation*}
\omega _{H_{p}}\left( \frac{1}{M_{n}},f\right) =o\left( \frac{1}{M_{n}^{1/p-1}}\right) ,\text{ when }n\rightarrow \infty,
\end{equation*}
then
$$\left\Vert S_{k}f-f\right\Vert _{H_{p}}\rightarrow 0,\,\,\,\text{as \ }\,\,\,k\rightarrow \infty .$$
Moreover, For every $0<p<1$ there exists a martingale $f\in H_{p}$, for which
\begin{equation*}
\omega _{H_{p}}\left( \frac{1}{M_{n}},f\right) =O\left( \frac{1}{%
	M_{n}^{1/p-1}}\right) ,\text{ \ when \ \ \ }n\rightarrow \infty
\end{equation*}%
\textit{and}
$$
\left\Vert S_{k}f-f\right\Vert _{weak-L_{p}}\nrightarrow 0,\,\,\,\text{%
	when\thinspace \thinspace \thinspace }k\rightarrow \infty .
$$

Tephadze \cite{tep17} proved that for any $0<p<1$ and $f\in H_{p},$ there exists an absolute constant $c_{p},$ depending only on $p,$ such that
\begin{equation*}
\text{ }\left\Vert S_{n}f\right\Vert _{H_{p}}\leq \frac{c_{p}M_{\left\vert n\right\vert }^{1/p-1}}{M_{\left\langle n\right\rangle }^{1/p-1}}\left\Vert f\right\Vert _{H_{p}}.
\end{equation*}

Moreover, for every $0<p<1\ \ $ and $\left\{ n_{k}:\text{ }k\in \mathbb{N}\right\} $ any an increasing sequence of nonnegative integers such that condition (\ref{ak}) is satisfied and $\left\{ \Phi _{n}:n\in \mathbb{N}\right\} $ for any non-decreasing sequence, satisfying the condition
\begin{equation} \label{12e}
\overline{\underset{k\rightarrow \infty }{\lim }}\frac{M_{\left\vert
		n_{k}\right\vert }^{1/p-1}}{M_{\left\langle n_{k}\right\rangle }^{1/p-1}\Phi_{n_{k}}}=\infty,  
\end{equation}
there exists a martingale $f\in H_{p},$ such that
\begin{equation*}
\underset{k\in \mathbb{N}}{\sup }\left\Vert \frac{S_{n_{k}}f}{\Phi _{n_{k}}}\right\Vert _{L_{p,\infty }}=\infty .
\end{equation*}
Moreover, if $0<p<1$, $f\in H_{p}$ and $\left\{ n_{k}:\text{ }k\in\mathbb{N}\right\} $ be an increasing sequence of  nonnegative integers, then
$
\left\Vert S_{n_{k}}f\right\Vert _{H_{p}}\leq c_{p}\left\Vert f\right\Vert_{H_{p}}
$
holds true	if and only if condition (\ref{ak0}) is satisfied.

In \cite{tep14} (see also \cite{tep17}) was proved that if $0<p<1,$ $f\in H_{p}\ \ $and $M_{k}<n\leq M_{k+1}.$ Then there exists an absolute constant $c_{p}$, depending only on $p,$ such that
\begin{equation*} 
\left\Vert S_{n}f-f\right\Vert _{H_{p}}\leq \frac{c_{p}M_{\left\vert
		n\right\vert }^{1/p-1}}{M_{\left\langle n\right\rangle }^{1/p-1}}\omega_{H_p} \left( \frac{1}{M_{k}},f\right),\ \ \ \left( 0<p<1\right) .
\end{equation*}
It follows that  if $\{n_{k}:k\in\mathbb{N}\}$  be increasing sequence of nonnegative integers such that
\begin{equation*}
\omega_{H_{p}} \left( \frac{1}{M_{\left\vert n_{k}\right\vert }},f\right)=o\left( \frac{M_{\left\langle n_{k}\right\rangle }^{1/p-1}}{M_{\left\vert n_{k}\right\vert }^{1/p-1}}\right) ,\text{ as \ }k\rightarrow \infty,
\end{equation*}
then
$\left\Vert S_{n_{k}}f-f\right\Vert _{H_{p}}\rightarrow 0,\,\,\,\text{as}\,\,\,k\rightarrow \infty.  $
Moreover, if $\{n_{k}:k\in\mathbb{N}\}$ be an increasing sequence of nonnegative
integers such that condition (\ref{ak}) is satisfied. Then
there exists a martingale $f\in H_{p}$ and a subsequence
$\{\alpha _{k}:k\in%
\mathbb{N}\}\subset \{n_{k}:k\in\mathbb{N}\},$ for which
\begin{equation*}  
\omega_{H_{p}}\left( \frac{1}{M_{\left\vert \alpha _{k}\right\vert }},f\right)=O\left( \frac{M_{\left\langle \alpha _{k}\right\rangle }^{1/p-1}}{M_{\left\vert \alpha _{k}\right\vert }^{1/p-1}}\right),\text{ as \ }k\rightarrow \infty
\end{equation*}
and
$
\limsup\limits_{k\rightarrow \infty }\left\Vert S_{\alpha_k}
f-f\right\Vert _{weak-L_{p}}>c>0,\,\,\,\text{as\thinspace \thinspace\thinspace }k\rightarrow \infty .  
$

In Tephnadze \cite{tep7} (see also \cite{tep12}) was proved that for every $f\in H_{1},$ the maximal operator
\begin{equation*}
\widetilde{S}^{\ast }f:=\sup_{n\in \mathbb{N}_{+}}\frac{\left\vert
	S_{n}f\right\vert }{\log \left( n+1\right) }\text{\ }
\end{equation*}%
is bounded from the Hardy space $H_{1}$ to the space $L_{1}.$ Moreover, the rate of the sequence $\log \left( n+1\right) $ is in the sense sharp.
Hence, for any $f\in H_{1},$ there exists an absolute
constant $c,$ such that%
\begin{equation} \label{dnl}
\left\Vert S_{n}f\right\Vert _{1}\leq c\log \left( n+1\right) \left\Vert
f\right\Vert _{H_{1}},\text{ \ }n\in \mathbb{N}_{+}.
\end{equation}

This estimate immediately follow that if $f\in H_{1}$ and $M_{k}<n\leq M_{k+1}$. Then there is
an absolute constant $c$ such that
\begin{equation*}
\left\Vert S_{n}f-f\right\Vert _{H_{1}}\leq c\lg n\omega _{H_{1}}\left( \frac{1}{%
	M_{k}},f\right) .
\end{equation*}

By using this estimate we obtain that if $f\in H_{1}\ $and
\begin{equation*}
\omega _{H_{1}}\left( \frac{1}{M_{n}},f\right) =o\left( \frac{1}{n}\right),%
\text{ when }n\rightarrow \infty ,
\end{equation*}%
Then
$
\left\Vert S_{k}f-f\right\Vert _{H_{1}}\rightarrow 0,\,\,\,\text{when\thinspace
	\thinspace \thinspace }k\rightarrow \infty .
$
Moreover (for details see \cite{tep7}) there exists a martingale $f\in H_{1},$ \ for which%
\begin{equation*}
\omega _{H_{1}}\left( \frac{1}{M_{2M_{n}}},f\right) =O\left( \frac{1}{M_{n}}%
\right) ,\ \ \ \text{when}\ \ n\rightarrow \infty
\end{equation*}%
and
$
\left\Vert S_{k}f-f\right\Vert _{1}\nrightarrow 0,\,\,\,\text{when\thinspace
	\thinspace \thinspace }k\rightarrow \infty .
$

In \cite{tep14} (see also \cite{tep17}) was proved that if $f\in H_{1}$ and $M_{k}<n\leq M_{k+1}$, then there
exists an absolute constant $c$ such that 
\begin{equation*}
\left\Vert S_{n}f\right\Vert _{H_{1}}\leq c\left( v\left( n\right) +v^{\ast
}\left( n\right) \right) \left\Vert f\right\Vert _{H_{1}}.
\end{equation*}

Moreover, if $\left\{ \Phi _{n}:n\in \mathbb{N}\right\} $ be any non-decreasing
and non-negative sequence satisfying condition 
$
\underset{n\rightarrow \infty }{\lim }\Phi _{n}=\infty
$
and $\ \left\{ n_{k}\geq 2:k\in \mathbb{N}\right\} $ be a subsequence such
that 
\begin{equation*}
\lim_{k\rightarrow \infty }\frac{v\left( n_{k}\right) +v^{\ast }\left(
	n_{k}\right) }{\Phi _{n_{k}}}=\infty,
\end{equation*}%
then there exists a martingale $f\in H_{1}$ such that 
\begin{equation*}
\sup_{k\in \mathbb{N}}\left\Vert \frac{S_{n_{k}}f}{\Phi _{n_{k}}}\right\Vert
_{1}\rightarrow \infty ,\text{ as }k\rightarrow \infty .
\end{equation*}

In \cite{tep14} (see also \cite{tep17}) it was also proved that if $f\in H_{1}$ and $M_{k}<n\leq M_{k+1},$ then there
exists an absolute constant $c$ such that 
\begin{equation*}
\left\Vert S_{n}f-f\right\Vert _{H_{1}}\leq c\left( v\left( n\right)
+v^{\ast }\left( n\right) \right) \omega _{H_{1}}\left( \frac{1}{M_{k}}%
,f\right) . 
\end{equation*}

It follows that if $f\in H_{1}$ and $\{n_{k}:k\in \mathbb{N}\}$ be a sequence of non-negative integers such that 
\begin{equation*}
\omega _{H_{1}}\left( \frac{1}{M_{\left\vert n_{k}\right\vert }},f\right)=o\left( \frac{1}{v\left( n_{k}\right) +v^{\ast }\left( n_{k}\right)}\right),\text{ as }k\rightarrow \infty,
\end{equation*}
then 
$
\left\Vert S_{n_{k}}f-f\right\Vert _{H_{1}}\rightarrow 0,\text{ when }k\rightarrow \infty.
$
Moreover, if $\{n_{k}:k\geq 1\}$ be sequence of non-negative integers such that 
$
\sup_{k\in \mathbb{N}}\left( v\left( n_{k}\right) +v^{\ast }\left(
n_{k}\right) \right) =\infty ,
$
then there exists a martingale $f\in H_{1}$ and a sequence $\{\alpha
_{k}:k\in \mathbb{N}\}\subset \{n_{k}:\in \mathbb{N}\}$ for which 
\begin{equation*}
\omega _{H_{1}}\left( \frac{1}{M_{\left\vert \alpha _{k}\right\vert }},f\right) =O\left( \frac{1}{v\left( \alpha _{k}\right) +v^{\ast }\left(\alpha _{k}\right) }\right)
\end{equation*}%
and 
$
\limsup\limits_{k\rightarrow \infty }\left\Vert S_{\alpha_{k}}f-f\right\Vert _{1}>c>0\text{ when }k\rightarrow \infty .
$

Simon \cite{Si11} proved that for any $f\in H_{p},$ there exists
an absolute constant $c_{p},$ depending only on $p,$ such that
\begin{equation*}
\overset{\infty }{\underset{k=1}{\sum }}\frac{\left\Vert S_{k}f\right\Vert
	_{p}^{p}}{k^{2-p}}\leq c_{p}\left\Vert f\right\Vert _{H_{p}}^{p},\text{ \ \
	\ }\left( 0<p<1\right).
\end{equation*}%

In Tephnadze \cite{tep4}) was proved  sharpness of this result in the special sense. In particular, if  $0<p<1$ and $\left\{ \Phi _{n}:n\in \mathbb{N}\right\} $ be any
non-decreasing sequence satisfying the condition
$
\underset{n\rightarrow \infty }{\overline{\lim }}\Phi _{n}=+\infty ,
$
there exists a martingale $f\in H_{p}$ such that
\begin{equation*}
\text{ }\underset{k=1}{\overset{\infty }{\sum }}\frac{\left\Vert
	S_{k}f\right\Vert _{weak-L_{p}}^{p}\Phi _{k}}{k^{2-p}}=\infty.
\end{equation*}

In G\'{a}t \cite{Ga1} the following strong convergence result was
obtained for all $f\in H_{1}$:
\begin{equation*}
\underset{n\rightarrow \infty }{\lim }\frac{1}{\log n}\overset{n}{\underset{%
		k=1}{\sum }}\frac{\left\Vert S_{k}f-f\right\Vert _{1}}{k}=0
\end{equation*}%
For the trigonometric analogue see Smith \cite{sm} and for the Walsh-Paley
system see Simon \cite{Si3}, for Vilenkin-like systems see Blahota \cite{bl1}. Moreover, for all $f\in H_{1}$, there exists an absolute constant $c,$ such that%
\begin{equation*}
\frac{1}{\log n}\overset{n}{\underset{k=1}{\sum }}\frac{\left\Vert
	S_{k}f\right\Vert _{1}}{k}\leq c\left\Vert f\right\Vert _{H_{1}} 
\text{ \ \ \ and \ \ \ }
\underset{n\rightarrow \infty }{\lim }\frac{1}{\log n}\overset{n}{\underset{k=1}{\sum }}\frac{\left\Vert S_{k}f\right\Vert _{1}}{k}=\left\Vert
f\right\Vert _{H_{1}}
\text{ \ }%
\left( n=2,3...\right).
\end{equation*}%

In the one-dimensional case Yano \cite{Yano} proved that
\begin{equation*}
\left\Vert K_{n}\right\Vert \leq 2\ \ \ (n\in \mathbb{N}).
\end{equation*}%
Consequently,
\begin{equation*}
\left\Vert \sigma _{n}f-f\right\Vert _{p}\rightarrow 0,\text{ \ \ \ when \ \ \
	\ }n\rightarrow \infty ,\text{ \ }(f\in L_{p},\text{ \ }1\leq p\leq \infty ).
\end{equation*}%
However (see \cite{JOO, sws}) the rate of convergence can not be better then
$O\left( n^{-1}\right) $ $\left( n\rightarrow \infty \right) $ for
non-constant functions. a.e, if $f\in L_{p},$ $1\leq p\leq \infty $ and
\begin{equation*}
\left\Vert \sigma _{M_{n}}f-f\right\Vert _{p}=o\left( \frac{1}{M_{n}}\right)
,\text{ \ when \ \ }n\rightarrow \infty ,
\end{equation*}%
then \textit{\ }$f$ \ is a constant function.

Fridli \cite{FR} used dyadic modulus of continuity to characterize the set
of functions in the space $L_{p}$, whose Vilenkin-Fej\'{e}r means converge
at a given rate. It is also known that (see e.g books \cite{AVD} and
\cite{sws})%
\begin{equation*}
\left\Vert \sigma _{n}f-f\right\Vert _{p}
\end{equation*}%
\begin{equation*}
\leq c_{p}\omega _{p}\left( \frac{1}{M_{N}},f\right) +c_{p}\sum_{s=0}^{N-1}%
\frac{M_{s}}{M_{N}}\omega _{p}\left( \frac{1}{M_{s}},f\right) ,\text{ \ }%
\left( 1\leq p\leq \infty ,\text{ \ }n\in \mathbb{N}\right) .
\end{equation*}

By applying this estimate, we immediately obtain that if $f\in lip\left(
\alpha ,p\right) ,$ i.e.,%
\begin{equation*}
\omega _{p}\left( \frac{1}{M_{n}},f\right) =O\left( \frac{1}{M_{n}^{\alpha }}%
\right) ,\text{ \ \ }n\rightarrow \infty ,
\end{equation*}%
then%
\begin{equation*}
\left\Vert \sigma _{n}f-f\right\Vert _{p}=\left\{
\begin{array}{ll}
O\left( \frac{1}{M_{N}}\right) , & \text{if }\alpha >1, \\
O\left( \frac{N}{M_{N}}\right) , & \text{if }\alpha =1, \\
O\left( \frac{1}{M_{N}^{\alpha }}\right) , & \text{if }\alpha <1\text{. }%
\end{array}%
\right.
\end{equation*}
On the other hand, if $1\leq p\leq \infty ,$ $f\in L_{p}$ and
\begin{equation*}
\left\Vert \sigma _{M_{n}}f-f\right\Vert_{p}=o\left( 1/M_{n}\right), \ \text{as} \ n\rightarrow \infty,
\end{equation*}
then $f$  is constant function $f=const.$

Weisz \cite{We2} considered the norm convergence of Fej\'er means of
Vilenkin-Fourier series and proved that
\begin{equation}
\left\Vert \sigma _{k}f\right\Vert _{p}\leq c_{p}\left\Vert f\right\Vert
_{H_{p}},\text{ \ \ }p>1/2\text{ \ \ and \ \ \ }f\in H_{p}.  \label{f100}
\end{equation}

This result implies that%
\begin{equation*}
\frac{1}{n^{2p-1}}\overset{n}{\underset{k=1}{\sum }}\frac{\left\Vert \sigma
	_{k}f\right\Vert _{p}^{p}}{k^{2-2p}}\leq c_{p}\left\Vert f\right\Vert
_{H_{p}}^{p},\text{ \ \ \ }\left( 1/2<p<\infty \right) .
\end{equation*}

If (\ref{f100}) hold for $0<p\leq 1/2,$ then we would have that
\begin{equation}
\frac{1}{\log ^{\left[ 1/2+p\right] }n}\overset{n}{\underset{k=1}{\sum }}%
\frac{\left\Vert \sigma _{k}f\right\Vert _{p}^{p}}{k^{2-2p}}\leq
c_{p}\left\Vert f\right\Vert _{H_{p}}^{p},\text{ \ \ \ }\left( 0<p\leq
1/2\right) .  \label{2cc}
\end{equation}

However, in Tephnadze \cite{tep1} it was shown that the assumption $p>1/2$ in (\ref%
{f100}) is essential. In particular, is was proved that there exists a
martingale $f\in H_{1/2}$ such that
\begin{equation*}
\sup_{n\in \mathbb{N}}\left\Vert \sigma _{n}f\right\Vert _{1/2}=+\infty .
\end{equation*}

For Vilenkin systems in \cite{tep5} it was proved that (\ref{2cc}) holds,
though inequality (\ref{f100}) is not true for $0<p\leq 1/2.$

Some new strong convergence result for Fejer means was considered in Persson, Tephnadze, Tutberidze and Wall \cite{pttw}.

In the one-dimensional case the weak type inequality
\begin{equation*}
\mu \left( \sigma ^{\ast }f>\lambda \right) \leq \frac{c}{\lambda }%
\left\Vert f\right\Vert _{1},\text{ \qquad }\left( f\in L_{1},\text{ \ \ \ }%
\lambda >0\right)
\end{equation*}%
can be found in Zygmund \cite{13z} for the trigonometric series, in Schipp
\cite{Sc} for Walsh series and in P\'al, Simon \cite{PS} for bounded Vilenkin
series. Fujji \cite{Fu} and Simon \cite{Si2} verified that $\sigma ^{\ast }$
is bounded from $H_{1}$ to $L_{1}$. Weisz \cite{We2} generalized this result
and proved the boundedness of $\sigma ^{\ast }$ from the martingale space $%
H_{p}$ to the Lebesgue space $L_{p}$ for $p>1/2$. Simon \cite{Si11} gave a
counterexample, which shows that boundedness does not hold for $0<p<1/2.$
The counterexample for $p=1/2$ due to Goginava \cite{GoAMH}, (see also \cite{BGG} and \cite{BGG2}). In \cite{tep1} Tephnadze proved that there exist a martingale $f\in H_{1/2}$ such that
\begin{equation*}
\sup\limits_{n\in \mathbb{N}}\left\Vert \sigma _{n}f\right\Vert
_{1/2}=+\infty.
\end{equation*}
Moreover, there exist a martingale $f\in H_{p},$  for $0<p<1/2,$ such that
\begin{equation*}
\sup\limits_{n\in \mathbb{N}}\left\Vert \sigma _{n}f\right\Vert
_{weak-L_{p}}=+\infty .
\end{equation*}

It follows that there exist a martingale $f\in H_{1/2}$ such that
\begin{equation*}
\left\Vert \sigma ^{\ast }f\right\Vert _{1/2}=+\infty .
\end{equation*}
Moreover, there exist a martingale $f\in H_{p}$ for $0<p<1/2,$ such that
\begin{equation*}
\left\Vert \sigma ^{\ast }f\right\Vert _{weak-L_{p}}=+\infty .
\end{equation*}

Weisz \cite{We4} proved that $\sigma ^{\ast }$ is
bounded from the Hardy space $H_{1/2}$ to the space $weak-L_{1/2}$. In  \cite{tep3}  it was proved
that the maximal operator $\widetilde{\sigma }_{p}^{\ast }$ with respect to
Vilenkin systems defined by
\begin{equation*}
\widetilde{\sigma }_{p}^{\ast }:=\sup_{n\in \mathbb{N}}\frac{\left\vert
	\sigma _{n}\right\vert }{\left( n+1\right) ^{1/p-2}},
\end{equation*}%
where $0<p< 1/2$, is bounded from the Hardy space $H_{p}$ to the Lebesgue space $L_{p}.$
Moreover, the order of deviant behavior of the $n$-th Fej\'er mean was given exactly. That is, for any  be any non-decreasing
sequence $\left\{ \Phi _{n}:n\in \mathbb{N}\right\} $ satisfying the condition%
\begin{equation*}
\overline{\lim_{n\rightarrow \infty }}\frac{\left( n+1\right) ^{1/p-2}}{\Phi
	_{n}}=+\infty,
\end{equation*}
we have that
\begin{equation*}
\sup_{k\in \mathbb{N}}\frac{\left\Vert \frac{\sigma _{M_{_{2n_{k}}}+1}f_{k}}{%
		\Phi _{M_{2n_k}+1}}\right\Vert _{weak-L_{p}}}{\left\Vert
	f_k\right\Vert_{H_p}}=\infty.
\end{equation*}

As a corollary we immediately get that
\begin{equation*}
\left\Vert \sigma _{n}f\right\Vert _{p}\leq c_{p}\left( n+1\right)
^{1/p-2}\left( n+1\right) \left\Vert
f\right\Vert _{H_{p}},
\end{equation*}
but it was proved more stronger result (for details see e.g. \cite{tep12}). In particular, 
if $0<p<1/2$ and $f\in H_{p},$  there exists an
absolute constant $c_{p}$, depending only on $p,$ such that
\begin{equation*}
\left\Vert \sigma _{n}f\right\Vert _{H_{p}}\leq c_{p}n^{1/p-2}\left\Vert
f\right\Vert _{H_{p}}.
\end{equation*}

In \cite{tep2}  (for Walsh system see \cite{GoSzeged}) it was proved
that the maximal operator $\widetilde{\sigma }^{\ast }$ with respect to
Vilenkin systems defined by
\begin{equation*}
\widetilde{\sigma }^{\ast }f:=\sup_{n\in \mathbb{N}}\frac{\left\vert
	\sigma _{n}f\right\vert }{\log ^{2 }\left( n+1\right) },
\end{equation*}%
is bounded from the Hardy space $H_{1/2}$ to the Lebesgue space $L_{1/2}.$

Moreover, for  any
non-decreasing sequence $\left\{ \Phi _{n}:n\in \mathbb{N}\right\} $ satisfying the condition%
\begin{equation*}
\overline{\lim_{n\rightarrow \infty }}\frac{\log ^{2}\left( n+1\right) }{\Phi _{n}}=+\infty,
\end{equation*}%
we have that
\begin{equation*}
\sup_{k\in \mathbb{N}}\frac{\left\Vert \frac{\sigma _{q_{n_{k}}}f_{k}}{\Phi_{q_{n_{k}}}}\right\Vert _{1/2}}{\left\Vert f_{k}\right\Vert _{H_{1/2}}}=\infty.
\end{equation*}
As a corollary we get that
\begin{equation*}
\left\Vert \sigma _{n}f\right\Vert _{1/2}\leq c\log^2 \left( n+1\right) \left\Vert
f\right\Vert _{H_{1/2}}.
\end{equation*}
but it was proved more stronger result (for details see e.g. \cite{tep12}). In particular, 
if $f\in H_{1/2},$  there exists an
absolute constant $c, $ such that
\begin{equation*} \label{strongsgn}
\left\Vert \sigma _{n}f\right\Vert _{H_{1/2}}\leq c\log^2 \left( n+1\right) \left\Vert
f\right\Vert _{H_{1/2}}.
\end{equation*}

For Walsh-Kaczmarz system analogical theorems were proved in \cite{GNCz} and
\cite{tep4}.

For the one-dimensional Vilenkin-Fourier series Weisz \cite{We2} proved that
the maximal operator
\begin{equation*}
\sigma ^{\#}f=\sup_{n\in \mathbb{N}}\left\vert \sigma _{M_{n}}f\right\vert
\end{equation*}%
is bounded from the martingale Hardy space $H_{p}$ to the Lebesgue space $L_{p}$ for $p>0.$ Moreover,  the
operator $\left\vert \sigma _{M_{n}}f\right\vert $ is not bounded from the
space $H_{p}$ to the space $H_{p},$ for $0<p\leq 1.$ This result for the Walsh system can be found in Goginava \cite{gog1} and for bounded Vilenkin systems in the paper of Persson and Tephnadze \cite{pt1}.

Approximation properties of subsequences of Fej\'er means with respect to the two-dimensional Walsh-Fourier series was considered in Persson, Tephnadze, Tutberidze \cite{ptt1} and Tutberidze \cite{tut2}.

Tephnadze \cite{pt2} proved that if  $0<p\leq 1/2$ and $\left\{  \alpha _{k}:k\in \mathbb{N}\right\} $ be a subsequence of positive numbers such that
\begin{equation*}
\sup_{k\in \mathbb{N}}\rho \left(  \alpha _{k}\right) =\varkappa <c<\infty, 
\end{equation*}
Then the maximal operator
\begin{equation*}
\widetilde{\sigma }^{\ast ,\vartriangle }f:=\sup_{k\in \mathbb{N}}\left\vert \sigma
_{ \alpha _{k}}f\right\vert
\end{equation*}%
is bounded from the Hardy space $H_{p}$ to the Lebesgue space $L_{p}.$

Moreover, if $0<p\leq 1/2\ \ $and $\left\{  \alpha _{k}:k\in \mathbb{N}\right\} $ be a
subsequence of positive numbers satisfying the condition
\begin{equation*}
\sup_{k\in \mathbb{N}}\rho \left(  \alpha _{k}\right) =\infty .  \label{fenk1}
\end{equation*}%
then there exists an martingale $f\in H_{p}$ such that%
\begin{equation*}
\sup_{k\in \mathbb{N}}\left\Vert \sigma _{ \alpha _{k}}f\right\Vert
_{weak-L_{p}}=\infty ,\text{ \ \ }\left( 0<p<1/2\right) .
\end{equation*}

It immediately follows that for $0<p\leq 1/2,$ and  $f\in H_{p},$  there exists an absolute constant $c_{p}$, depending only on $p$, such that 
\begin{equation*}
\text{ }\left\Vert \sigma _{n_k}f\right\Vert _{p}\leq c_{p}\left\Vert
f\right\Vert _{H_{p}},\text{ \ \ }k\in \mathbb{N}
\end{equation*}
if and only if 
\begin{equation*}
\sup_{k\in \mathbb{N}}\rho\left( n_{k}\right)<c<\infty.
\end{equation*}

As a consequence, for  $p>0$ and $f\in H_{p}$, then there exists an absolute constant $c_{p}$, depending only on $p,$ such that
\begin{equation} \label{sigmamn}
\left\Vert \sigma _{M_{n}}f\right\Vert _{p}\leq c_{p}\left\Vert
f\right\Vert _{H_{p}},\ \ \ \left( p>0\right).
\end{equation}

In \cite{tep6} it was proved that if  $0<p<1/2,$ $f\in H_{p}$ and
\begin{equation*}
\omega _{p}\left( \frac{1}{M_{n}},f\right) =o\left( \frac{1}{M_{n}^{1/p-2}}\right) \text{when \ }n\rightarrow \infty,  
\end{equation*}
then
\begin{equation*}
\left\Vert \sigma _{n}f-f\right\Vert_{H_{p}}\rightarrow 0,\text{ when }n\rightarrow \infty .
\end{equation*}

Moreover, there exists a martingale $f\in
H_{p}$ $(0<p<1/2)$ for which
\begin{equation*}
\omega \left( \frac{1}{M_{n}},f\right) _{H_{p}}=O\left( \frac{1}{%
	M_{n}^{1/p-2}}\right) \text{ \ when \ }n\rightarrow \infty
\end{equation*}%
\textit{and}
\begin{equation*}
\left\Vert \sigma _{n}f-f\right\Vert _{weak-L_{p}}\nrightarrow 0,\,\,\,\text{%
	when\thinspace \thinspace \thinspace }n\rightarrow \infty .
\end{equation*}

When $p=1/2$ we have the following results:
If  $f\in H_{1/2}$ and
\begin{equation}
\omega _{H_{1/2}}\left( \frac{1}{M_{n}},f\right) =o\left( \frac{1}{n^{2}}%
\right) ,\text{ when \ }n\rightarrow \infty,  \label{fecon2}
\end{equation}%
then
\begin{equation*}
\left\Vert \sigma _{n}f-f\right\Vert_{H_{1/2}}\rightarrow 0,\text{ when }%
n\rightarrow \infty .
\end{equation*}%
Moreover, there exists a martingale $f\in H_{1/2}$ for which
\begin{equation*}
\omega _{H_{1/2}}\left( \frac{1}{M_{n}},f\right) =O\left( \frac{1}{n^{2}}%
\right) ,\text{ \ when \ }n\rightarrow \infty
\end{equation*}%
and
\begin{equation*}
\left\Vert \sigma _{n}f-f\right\Vert _{1/2}\nrightarrow 0,\,\,\,\text{%
	when\thinspace \thinspace \thinspace }n\rightarrow \infty .
\end{equation*}

We state consequences of this result for Walsh system to clearly see difference of divergence rates for the various subsequences:
Let $0<p<1/2,$ \  $f\in H_{p}.$ Then there exists an absolute constant $c_{p}$, depending only on $p$, such that 
\begin{equation} \label{aaaa}
\left\Vert \sigma _{M_{n}+1}f\right\Vert _{H_{p}}\leq c_{p}M_n ^{1/p-2}\left\Vert f\right\Vert _{H_{p}},\text{ \ \ }n\in \mathbb{N}
\end{equation}
and
\begin{equation} \label{aaa11}
\left\Vert \sigma _{M_{n}+M_{\left[n/2\right]}}f\right\Vert _{H_{p}}\leq c_{p}M_n ^{1/2p -1}\left\Vert f\right\Vert _{H_{p}},\text{ \ \ }n\in \mathbb{N}.
\end{equation}
Moreover, the rates $M_n ^{1/p-2}$ and  $ M_n ^{1/2p -1}$ in inequalities (\ref{aaaa}) and (\ref{aaa11}) are sharp in the same sense.

Blahota and Tephnadze \cite{bt1} proved that if  $0<p<1/2$ and $f\in H_{p},$ then there exists an absolute constant $c_{p},$ depending only on $p$, such that
\begin{equation*}
\overset{\infty }{\underset{k=1}{\sum }}\frac{\left\Vert \sigma
	_{k}f\right\Vert _{p}^{p}}{k^{2-2p}}\leq c_{p}\left\Vert f\right\Vert
_{H_{p}}^{p},
\end{equation*}
Moreover, if $0<p<1/2$\textit{\ and }$\left\{ \Phi _{k}\ :k\in \mathbb{N}%
\right\} $ be  any non-decreasing sequence satisfying the conditions $\Phi _{n}\uparrow \infty $  and
\begin{equation}
\overline{\underset{k\rightarrow \infty }{\lim }}\frac{k^{2-2p}}{\Phi _{k}}%
=\infty,  \label{cond000}
\end{equation}%
there exists a martingale $f\in H_{p}$  such that 
\begin{equation*}
\underset{k=1}{\overset{\infty }{\sum }}\frac{\left\Vert \sigma_{k}f\right\Vert _{weak-L_{p}}^{p}}{\Phi _{k}}=\infty.
\end{equation*}

As a corollary we also get that if $0<p<1/2$ and $f\in H_{p},$ then there
exists an absolute constant $c_{p},$ depending only on $p$, such that
\begin{equation*}
\overset{\infty }{\underset{k=1}{\sum }}\frac{\left\Vert \sigma
	_{k}f\right\Vert _{H_{p}}^{p}}{k^{2-2p}}\leq c_{p}\left\Vert f\right\Vert
_{H_{p}}^{p},
\end{equation*}
\begin{equation*}
\frac{1}{n}\overset{n}{\underset{k=1}{\sum }}\frac{\left\Vert \sigma
	_{k}f\right\Vert _{H_{p}}^{p}}{k^{1-2p}}\leq c_{p}\left\Vert f\right\Vert
_{H_{p}}^{p},
\end{equation*}%
\begin{equation*}
\frac{1}{n}\overset{n}{\underset{k=1}{\sum }}\frac{\left\Vert \sigma
	_{k}f-f\right\Vert _{H_{p}}^{p}}{k^{1-2p}}=0,
\end{equation*}%
and	
\begin{equation*}
\frac{1}{n}\overset{n}{\underset{k=1}{\sum }}\frac{\left\Vert \sigma
	_{k}f\right\Vert _{H_{p}}^{p}}{k^{1-2p}}=\left\Vert f\right\Vert
_{H_{p}}^{p}.
\end{equation*}

In Blahota and Tephnadze \cite{bt1} also considered the endpoint case $p=1/2$ and they proved that if $f\in H_{1/2}$ then there exists an absolute constant $c$ such that
\begin{equation*}
\frac{1}{\log n}\overset{n}{\underset{k=1}{\sum }}\frac{\left\Vert \sigma
	_{k}f\right\Vert _{1/2}^{1/2}}{k}\leq c\left\Vert f\right\Vert
_{H_{1/2}}^{1/2}.
\end{equation*}
As a corollary we also get that that if  $f\in H_{1/2},$ then
\begin{equation*}
\frac{1}{\log n}\overset{n}{\underset{k=1}{\sum }}\frac{\left\Vert \sigma
	_{k}f\right\Vert _{H_{1/2}}^{1/2}}{k}\leq c\left\Vert f\right\Vert
_{H_{1/2}}^{1/2},\text{ }
\end{equation*}%
\begin{equation*}
\lim_{n\rightarrow \infty }\frac{1}{\log n}\overset{n}{\underset{k=1}{\sum }}
\frac{\left\Vert \sigma _{k}f-f\right\Vert _{H_{1/2}}^{1/2}}{k}=0
\end{equation*}%
and
\begin{equation*}
\lim_{n\rightarrow \infty }\frac{1}{\log n}\overset{n}{\underset{k=1}{\sum }}
\frac{\left\Vert \sigma _{k}f\right\Vert _{H_{1/2}}^{1/2}}{k}=\left\Vert
f\right\Vert _{H_{1/2}}^{1/2}.
\end{equation*}

\subsection{Estimates of Dirichlet and Fej\'er Kernels  with respect to	Vilenkin systems}

The proof of Lemma can be found in Tephnadze \cite{tep6}.
\begin{lemma} \label{dn2.6}
	Let $x\in I_s\backslash I_{s+1}, \ \  s=0,...,N-1.$  Then
	\begin{equation*}
	\int_{I_N}\left\vert D_n\left(x-t\right)\right\vert d\mu\left(
	t\right)\leq\frac{cM_s}{M_N},
	\end{equation*}%
	where $c$ is an absolute constant.
\end{lemma}

\textbf{Proof:}
Let $x\in I_s\backslash I_{s+1}, \ \  s=0,...,N-1.$
By combining \eqref{3aa} and \eqref{2dna}we
have that
\begin{equation*}
\left\vert D_n\left(x\right)\right\vert \leq\underset{j=0}{\overset{s} {\sum}}n_jD_{M_j}\left( x\right)=\underset{j=0}{\overset{s}{\sum}}
n_jM_j\leq cM_s.
\end{equation*}

Since $t\in I_N$ and $x\in I_s\backslash I_{s+1}, \ \ s=0,...,N-1,$ we
obtain that $x-t\in I_s\backslash I_{s+1}$. By using the estimate above we get
that
\begin{equation*}
\left\vert D_n\left(x-t\right)\right\vert \leq cM_s
\end{equation*}%
and
\begin{equation*}
\int_{I_N}\left\vert D_n\left(x-t\right)\right\vert d\mu\left(t\right) \leq\frac{cM_s}{M_N}.
\end{equation*}

The proof is complete.

\QED

The proof of the next lemma \ref{lemma222} can be found in Tephnadze \cite{tep2, tep3}.
\begin{lemma} \label{lemma222}
	Let $n\in \mathbb{N}$ and $x\in I_N^{k,l},$ where $k<l.$
	Then
	\begin{equation} \label{star1}
	K_{M_n}\left(x\right)=0,\text{ \ if \ } n>l.
	\end{equation}
	and
	\begin{equation} \label{star2}
	\left\vert K_{M_n}\left(x\right)\right\vert\leq cM_k.
	\end{equation}%
	Moreover,%
	\begin{equation} \label{star3}
	\int_{G_m}\left\vert K_{M_n}\right\vert d\mu\leq c<\infty,
	\end{equation}%
	where $c$ is an absolute constant.
\end{lemma}

Next Lemma is proved in the book of G. N. Agaev, N. Ya. Vilenkin, G. M. Dzhafarly and A. I.Rubinshtein \cite{AVD}.

\begin{lemma}\label{lemma7kn}
	Let $n\in \mathbb{N}.$ Then
	\begin{equation} \label{fn5}
	n\left\vert K_n\right\vert\leq c\sum_{l=\left\langle n\right\rangle}^ {\left\vert n\right\vert}M_l\left\vert K_{M_l}\right\vert\leq
	c\sum_{l=0}^{\left\vert n\right\vert }M_l\left\vert K_{M_l}\right\vert
	\end{equation}%
	where $c$ is an absolute constant.
\end{lemma}

The proof of Lemmas \ref{lemma5} and \ref{lemma5aa} is due to Tephnadze \cite{tep2, tep3}  (see also Blahota, Tephnadze \cite{bt1}).
\begin{lemma} \label{lemma5}
	Let $x\in I_N^{k,l}, \ k=0,\dots ,N-2, \ l=k+1,\dots ,N-1.$ \ Then
	\begin{equation*}
	\int_{I_N}\left\vert K_n\left(x-t\right)\right\vert d\mu\left(t\right) \leq\frac{cM_lM_k}{nM_N}.
	\end{equation*}
	
	Let \ $x\in I_N^{k,N}, \ k=0,\dots ,N-1.$  Then
	\begin{equation*}
	\int_{I_{N}}\left\vert K_{n}\left( x-t\right) \right\vert d\mu \left(
	t\right) \leq \frac{cM_{k}}{M_{N}},
	\end{equation*}%
	where $c$ is an absolute constant.
\end{lemma}

The next lemma is a simple consequence of Lemma \ref{lemma5}.

\begin{lemma} \label{lemma5aa}
	Let $x\in I_N^{k,l}, \qquad k=0,\dots ,N-1, \qquad l=k+1,\dots ,N.$ \qquad Then
	\begin{equation*}
	\int_{I_N}\left\vert K_n\left(x-t\right)\right\vert d\mu\left(t\right) \leq \frac{cM_lM_k}{M_N^2},\text{ \ for \ }n\geq M_N,
	\end{equation*}%
	where $c$ is an absolute constant.
\end{lemma}

Now we prove some bellow estimate of Fej\'er kernel, wich will be used to prove some negative results. The proof of Lemma\ref{lemma6kn} is proved by Blahota and Tephnadze \cite{bt1}.

\begin{lemma}\label{lemma6kn}
	Let $t,s_n,$ $1\leq s_n\leq m_n-1$ $n\in \mathbb{N}.$ Then 
	\begin{equation} \label{100kn1}
	\left\vert K_{s_nM_n}\left(x\right)\right\vert\geq \frac{M_{n}}{2\pi s_n}, \text{ for }x\in I_{n+1}\left( e_{n-1}+e_n\right).
	\end{equation}%
	Moreover, if $x\in I_t\backslash I_{t+1},$  $x-x_te_t\notin I_n$
	and $n>t,$ then
	\begin{equation} \label{kn1}
	K_{s_nM_n}(x)=0.
	\end{equation}
\end{lemma}

\QED

The proof of Lemma \ref{lemma8ccc} can be found in Persson and Tephnadze \cite{pt2}.
\begin{lemma}\label{lemma8ccc} 
	Let $n\in \mathbb{N}$, $\left\langle n\right\rangle \neq\left\vert n\right\vert $ and $x\in I_{\left\langle n\right\rangle
		+1}\left(e_{\left\langle n\right\rangle-1}+e_{\left\langle n\right\rangle}\right).$ Then
	\begin{equation*}
	\left\vert nK_n(x)\right\vert =\left\vert \left(n-M_{\left\vert n\right\vert}\right)K_{n-M_{\left\vert n\right\vert}}(x)\right\vert \geq \frac{M_{\left\langle n\right\rangle }^2}{2\pi\lambda},
	\end{equation*}%
	where $\lambda :=\sup m_n.$
\end{lemma}

For the Walsh system analogical of Lemma \ref{lemma3} was proved in Tephnadze \cite{tep5} and Corollary \ref{cor3a} is simple consequence of it. Similar bellow estimate is proved in Blahota, G\'{a}t and Goginava \cite{BGG2} and \cite{BGG}.
\begin{lemma}\label{lemma3} Let 
	\begin{equation*}
	n=\sum_{i=1}^{s}\sum_{k=l_{i}}^{m_{i}}n_{k}M_{k},
	\end{equation*}
	\ where 
	\begin{equation*}
	0\leq l_{1}\leq m_{1}\leq l_{2}-2<l_{2}\leq m_{2}\leq ...\leq
	l_{s}-2<l_{s}\leq m_{s}.
	\end{equation*}
	
	Then 
	\begin{equation*}
	n\left\vert K_{n}\left( x\right) \right\vert \geq cM_{l_{i}}^{2},\text{ \ \
		for \ \ }x\in I_{l_{i}+1}\left( e_{l_{i}-1}+e_{l_{i}}\right) ,
	\end{equation*}
	where $\lambda =\sup_{n\in \mathbb{N}}m_{n}$ and $c$ is an absolute constant.
\end{lemma}
\textbf{Proof:}
Let $x\in I_{l_i+1}\left( e_{l_i-1}+e_{l_i}\right).$ By combining \eqref{3aa} and  \eqref{9dn}  with equality (\ref{star1}) in Lemma \ref{lemma3} we obtain that 
\begin{equation*}
D_{M_{l_i}}=0
\end{equation*}
and 
\begin{equation*}
D_{s_{n_k}M_{s_{n_k}}}=K_{s_{n_k}M_{{s_{n_k}}}}=0, \ \ s_{n_k}> l_{i}.
\end{equation*}

Since $s_{n_1}>s_{n_2}>\dots >s_{n_r}\geq 0$ we find that 
\begin{eqnarray*}
	n^{(k)}&=& n-\sum_{i=1}^{k}s_{n_i}M_{n_{i}}\\
	&=&\sum_{i=k+1}^{s}s_{n_i}M_{n_{i}}
	\leq  \sum_{i=0}^{n_{k+1}}(m_{i}-1)M_{i} \\
	&=& m_{n_{k+1}}M_{n_{k+1}}-1\leq M_{n_{k}}.
\end{eqnarray*}

According to \eqref{kn10} we have that 
\begin{eqnarray*}
	n\left\vert K_{n}\right\vert &\geq &\left\vert
	s_{l_{i}}M_{l_{i}}K_{s_{l_{i}}M_{l_{i}}}\right\vert \\
	&-&\sum_{r=1}^{i-1}\sum_{k=l_{r}}^{m_{r}}\left\vert
	s_{k}M_{k}K_{s_{k}M_{k}}\right\vert  \notag \\
	&-&\sum_{r=1}^{i-1}\sum_{k=l_{r}}^{m_{r}}\left\vert
	M_{k}D_{s_{k}M_{k}}\right\vert \\
	&=& I_1-I_2-I_3.
\end{eqnarray*}

Let $x\in I_{l_{i}+1}\left( e_{l_{i}-1}+e_{l_{i}}\right) $ and $1\leq
s_{l_{i}}\leq m_{l_{i}}-1$. By using Lemma \ref{lemma3} we get that 
\begin{eqnarray}  \label{10.0}
I_1=\left\vert s_{l_{i}}M_{l_{i}}K_{s_{l_{i}}M_{l_{i}}}\right\vert
\geq \frac{M_{l_{i}}^{2}}{2\pi }\geq \frac{2M_{l_{i}}^{2}}{9}.  \notag
\end{eqnarray}

It is easy to see that%
\begin{eqnarray*}
	\sum_{s=0}^{k}n_{s}^{2}M_{s}^{2}&\leq& \sum_{s=0}^{k}\left( m_{s}-1\right) ^{2}M_{s}^{2} \\
	&\leq
	&\sum_{s=0}^{k}m_{s}^{2}M_{s}^{2}-2\sum_{s=0}^{k}m_{s}M_{s}^{2}+%
	\sum_{s=0}^{k}M_{s}^{2} \\
	&=&\sum_{s=0}^{k}M_{s+1}^{2}-2\sum_{s=0}^{k}M_{s+1}M_{s}+%
	\sum_{s=0}^{k}M_{s}^{2} \\
	&=&
	M_{k+1}^{2}+2\sum_{s=0}^{k}M_{s}^{2}-2\sum_{s=0}^{k}M_{s+1}M_{s}-M_{0}^{2} \\
	&\leq & M_{k+1}^{2}-1.
\end{eqnarray*}%
and%
\begin{eqnarray*}
	\sum_{s=0}^{k}n_{s}M_{s}&\leq& \sum_{s=0}^{k}\left( m_{s}-1\right) M_{s} \\
	&=&m_{k}M_{k}-m_{0}M_{0} \\
	&\leq& M_{k+1}-2.
\end{eqnarray*}

Since $m_{i-1}\leq l_{i}-2$ if we use the estimates above, then we obtain that

\begin{eqnarray}  \label{10.1}
I_2 &\leq &\sum_{s=0}^{l_{i}-2}\left\vert n_{s}M_{s}K_{n_{s}M_{s}}\left(
x\right) \right\vert  \\  \notag
&\leq&\sum_{s=0}^{l_{i}-2}n_{s}M_{s}\frac{\left( n_{s}M_{s}+1\right) }{2} \\ \notag
&\leq &\frac{\left( m_{l_{i}-2}-1\right) M_{l_{i}-2}}{2}\sum_{s=0}^{l_{i}-2}%
\left( n_{s}M_{s}+1\right) \\ \notag
&\leq &\frac{\left( m_{l_{i}-2}-1\right) M_{l_{i}-2}}{2}M_{l_{i}-1} \\
&+& \frac{\left( m_{l_{i}-2}-1\right) M_{l_{i}-2}}{2}l_{i}  \notag \\
&\leq &\frac{M_{l_{i}-1}^{2}}{2}-\frac{M_{l_{i}-2}M_{l_{i}-1}}{2}%
+M_{l_{i}-1}l_{i}.  \notag
\end{eqnarray}

For $I_3$ we have that%
\begin{eqnarray}  \label{10.3}
I_3 &\leq &\sum_{k=0}^{l_{i}-2}\left\vert M_{k}D_{n_{k}M_{k}}\left( x\right)
\right\vert \\
&\leq& \sum_{k=0}^{l_{i}-2}n_{k}M_{k}^{2}  \notag \\ \notag
&\leq& M_{l_{i}-2}\sum_{k=0}^{l_{i}-2}n_{k}M_{k}  \\
&\leq& M_{l_{i}-1}M_{l_{i}-2}-2M_{l_{i}-2}.  \notag
\end{eqnarray}

By combining (\ref{10.0})-(\ref{10.3}) we have that%
\begin{eqnarray*}
	n\left\vert K_{n}\left( x\right) \right\vert &\geq & I_1-I_2-I_3 \\
	&\geq &\frac{M_{l_{i}}^{2}}{2\pi }+\frac{3}{2}+2M_{l_{i}-2} \\
	&-&\frac{M_{l_{i}-1}M_{l_{i}-2}}{2}-\frac{M_{l_{i}-1}^{2}}{2}%
	-M_{l_{i}-1}l_{i} \\
	&\geq &\frac{M_{l_{i}}^{2}}{2\pi }-\frac{M_{l_{i}}^{2}}{16}-\frac{%
		M_{l_{i}}^{2}}{8}+\frac{7}{2}-M_{l_{i}-1}l_{i} \\
	&\geq &\frac{2M_{l_{i}}^{2}}{9}-\frac{3M_{l_{i}}^{2}}{16}+\frac{7}{2}%
	-M_{l_{i}-1}l_{i} \\
	&\geq &\frac{M_{l_{i}}^{2}}{144}-M_{l_{i}-1}l_{i}.
\end{eqnarray*}

Suppose that $l_{i}\geq 4$. Then 
\begin{eqnarray*}
	n\left\vert K_{n}\left( x\right) \right\vert &\geq& \frac{M_{l_{i}}^{2}}{36}-\frac{M_{l_{i}}}{4} \\
	&\geq& \frac{M_{l_{i}}^{2}}{36}-\frac{M_{l_{i}}^{2}}{64}\\
	&\geq& \frac{5M_{l_{i}}^{2}}{36\cdot 16} \\
	&\geq& \frac{M_{l_{i}}^{2}}{144}.
\end{eqnarray*}
The proof is complete.

\QED

Next corollary is simple consequence of Lemma \ref{lemma3}:
\begin{corollary} \label{cor3a}
	Let $2<n\in \mathbb{N}_+$ and $
	q_n=M_{2n}+M_{2n-2}+...+M_2+M_0.$ Then \index{N}{$K_{q_{n-1}}$}
	\begin{equation*}
	q_{n-1}\left\vert K_{q_{n-1}}(x)\right\vert\geq \frac{M^2_{2k}}{144} ,\text{ \ \ for \ \ } x\in I_{2k+1}\left(e_{2k-1}+e_{l_{2k}}\right),
	\end{equation*}
	where $k=0,1,...,n.$
\end{corollary}

\subsection{Strong convergence of partial sums of Vilenkin-Fourier series on martingale Hardy spaces}

In this section we investigate some new strong convergence result of of partial sums of Vilenkin-Fourier series:

\begin{theorem}
	\label{theorem1}a) Let $f\in H_{1}.$ Then there exists an absolute constant $%
	c,$ such that
	\begin{equation*}
	\sup_{n\in \mathbb{N}}\frac{1}{n\log n}\overset{n}{\underset{k=1}{\sum }}%
	\left\Vert S_{k}f\right\Vert _{1}\leq \left\Vert f\right\Vert_{H_1}.
	\end{equation*}
	
	b) Let $\varphi :\mathbb{N_+}\rightarrow \lbrack 1,$ $\infty )$ be a nondecreasing
	function satisfying the condition
	\begin{equation}
	\overline{\lim_{n\rightarrow \infty }}\frac{\log n}{\varphi _{n}}=+\infty .
	\label{cond1aaa}
	\end{equation}
	
	Then there exists a function $f\in H_{1},$ such that
	\begin{equation*}
	\sup_{n\in \mathbb{N}}\frac{1}{n\varphi_n}\overset{n}{\underset{k=1}{\sum }}%
	\left\Vert S_{k}f\right\Vert _{1}=\infty .
	\end{equation*}
\end{theorem}

\textbf{Proof:}
By using (\ref{dnl}) we can conclude that 
\begin{equation*}
\frac{1}{n\log n}\overset{n}{\underset{k=1}{\sum }}%
\left\Vert S_{k}f\right\Vert _{1}\leq  \frac{c\left\Vert f\right\Vert_{H_1}}{n\log n}\overset{n}{\underset{k=1}{\sum }}{\log k}\leq c\left\Vert f\right\Vert_{H_1}.
\end{equation*}
and the proof of part a) is complete.

Under the condition \eqref{cond1aaa} there exists an increasing sequence of the positive integers $\left\{\alpha_{k}:k\in \mathbb{N}\right\} $ such that
\begin{equation*}
\overline{\lim_{k\rightarrow \infty }}\frac{\log M_{{\alpha _{k}}}}{\varphi
	_{2M_{\alpha _{k}}}}=+\infty
\end{equation*}%
and
\begin{equation} \label{69a}
\sum_{k=0}^{\infty }\frac{\varphi _{2M_{\alpha _{k}}}^{1/2}}{\log ^{1/2}M_{{%
			\alpha _{k}}}}<c<\infty .  
\end{equation}

Let $f=( f^{(n)},\ \ n\in \mathbb{N}) $ be martingale, defined by 
\begin{equation*}
f^{(n)}:=\sum_{\left\{ k;\text{ }2\alpha_{k}<n\right\} }\lambda _{k}a_{k},
\end{equation*}
where
\begin{equation*}
a_{k}=r_{\alpha _{k}}D_{M_{_{\alpha _{k}}}}=D_{2M_{_{\alpha
			_{k}}}}-D_{M_{_{\alpha _{k}}}}
\end{equation*}%
and
\begin{equation*}
\lambda _{k}=\frac{\varphi _{2M_{\alpha _{k}}}^{1/2}}{\log ^{1/2}M_{{\alpha
			_{k}}}}.
\end{equation*}

By the definition of $H_1$ and Lemma \ref{lemma2.1}, if we apply (\ref{69a}) we can conclude that $ f\in H_{1}.$ Moreover,
\begin{equation} \label{600}
\widehat{f}(j)=\left\{
\begin{array}{l}
{\lambda _{k}},\,\,\text{\ \ \ \ \thinspace \thinspace }j\in \left\{
M_{\alpha _{k}},...,2M_{\alpha _{k}}-1\right\} ,\text{ }k\in \mathbb{N} \\
0\text{ },\text{ \thinspace \qquad \thinspace \thinspace \thinspace
	\thinspace \thinspace }j\notin \bigcup\limits_{k=1}^{\infty }\left\{
M_{\alpha _{k}},...,2M_{\alpha _{k}}-1\right\} .\text{ }%
\end{array}%
\right.   
\end{equation}

Since
\begin{equation*}
D_{j+M_{\alpha _{k}}}=D_{M_{\alpha _{k}}}+\psi _{_{M_{\alpha _{k}}}}D_{j},%
\text{ \qquad when \thinspace \thinspace }j\leq M_{\alpha _{k}},
\end{equation*}
if we apply (\ref{600}) we obtain that
\begin{eqnarray}
S_{j}f &=&S_{M_{\alpha _{k}}}f+\sum_{v=M_{\alpha _{k}}}^{j-1}\widehat{f}%
(v)\psi _{v}  \label{8} \\
&=&S_{M_{\alpha _{k}}}f+{\lambda _{k}}\sum_{v=M_{\alpha _{k}}}^{j-1}\psi _{v}
\notag \\
&=&S_{M_{\alpha _{k}}}f+{\lambda _{k}}\left( D_{j}-D_{M_{\alpha _{k}}}\right)
\notag \\
&=&S_{M_{\alpha _{k}}}f+{\lambda _{k}}\psi _{M_{\alpha _{k}}}D_{j-M_{\alpha
		_{k}}}  \notag \\
&=&I_{1}+I_{2}.  \notag
\end{eqnarray}

In view of (\ref{smn}) we can write that
\begin{eqnarray}  \label{8bbb2}
\left\Vert I_{1}\right\Vert_{1} \leq\left\Vert
S_{M_{\alpha_{k}}}f\right\Vert _{1} \leq c\left\Vert f\right\Vert_{H_{1}}.
\end{eqnarray}

By combining (\ref{8bbb2}) with  lower estimate in \ref{var1} we get that
\begin{eqnarray*}
	\left\Vert S_{n}f\right\Vert _{1} &\geq& \left\Vert I_{2}\right\Vert
	_{1}-\left\Vert I_{1}\right\Vert _{1} \\
	&\geq& {\lambda _{k}}{L\left( {%
			n-M_{\alpha _{k}}}\right)}-c\left\Vert f\right\Vert _{H_{1}} \\
	&\geq& c{\lambda _{k}}{v\left( {%
			n-M_{\alpha _{k}}}\right)}-c\left\Vert f\right\Vert _{H_{1}}.
\end{eqnarray*}

Hence, by applying \eqref{var2} we find that
\begin{eqnarray*}
	&&\underset{n\in \mathbf{\mathbb{N}}_{+}}{\sup }\frac{1}{n\varphi _{n}}%
	\underset{k=1}{\overset{n}{\sum }}\left\Vert S_{k}f\right\Vert _{1} \\
	&\geq &\frac{1}{2M_{\alpha _{k}}\varphi _{2M_{\alpha _{k}}}}\underset{%
		\left\{ M_{\alpha _{k}}\leq l\leq 2M_{\alpha _{k}}\right\} }{\sum }%
	\left\Vert S_{l}f\right\Vert _{1} \\
	&\geq &\frac{c}{2M_{\alpha _{k}}\varphi _{2M_{\alpha _{k}}}}\underset{%
		\left\{ M_{\alpha _{k}}\leq l\leq 2M_{\alpha _{k}}\right\} }{\sum }\left(
	\frac{v\left( l-M_{\alpha _{k}}\right) \varphi _{2M{\alpha _{k}}}^{1/2}}{%
		\log ^{1/2}M_{\alpha _{k}}}-c\left\Vert f\right\Vert _{H_{1}}\right) \\
	&\geq &\frac{c\varphi _{2M{\alpha _{k}}}^{1/2}}{2M_{\alpha _{k}}\log
		^{1/2}M_{\alpha _{k}}\varphi _{2M_{\alpha _{k}}}}\underset{l=1}{\overset{%
			M_{\alpha _{k}}-1}{\sum }}v\left( l\right) -c\left\Vert f\right\Vert
	_{H_{1}}^{1/2} \\
	&\geq &\frac{c\varphi _{2M{\alpha _{k}}}^{1/2}\log M_{\alpha _{k}}}{\log
		^{1/2}M_{\alpha _{k}}\varphi _{2M_{\alpha _{k}}}}\\
	&\geq& \frac{c\log
		^{1/2}M_{\alpha _{k}}}{\varphi _{2M_{\alpha _{k}}}^{1/2}}\rightarrow \infty ,%
	\text{ as \ }k\rightarrow \infty .
\end{eqnarray*}	
The proof is complete. 

\QED

Using theorem of Weisz  \cite{We2} in the case $p=1$ we get that there exists an absolute constant $c,$ such that 
\begin{equation*} 
\sup_{n\in\mathbb{N}}\left\Vert\sigma_nf\right\Vert_{1}<c\left\Vert f \right\Vert_{H_1}.
\end{equation*}
That is,
\begin{equation} \label{sinma}
\sup_{n\in\mathbb{N}}\left\Vert\frac{1}{n}\overset{n}{\underset{k=1}{\sum }}
S_kf\right\Vert_1<c\left\Vert f \right\Vert_{H_1}.
\end{equation}
This estimate arise an interesting question, if there exists an absolute absolute constant $c,$ such that the following strong convergence result
\begin{equation*}
\sup_{n\in\mathbb{N}}\frac{1}{n}\overset{n}{\underset{k=1}{\sum }}\left\Vert
S_{k}f\right\Vert_{1}<c\left\Vert f \right\Vert_{H_1}.
\end{equation*}
holds true, which is stronger inequality than \eqref{sinma}. In particular, we have negative answer on this question:
\begin{corollary}
	\label{theorem2sn} There exists a function $f\in H_{1},$ such that
	\begin{equation*}
	\sup_{n\in\mathbb{N}}\frac{1}{n}\overset{n}{\underset{k=1}{\sum }}\left\Vert
	S_{k}f\right\Vert_{1}=\infty .
	\end{equation*}
\end{corollary}

\subsection{Strong convergence of Vilenkin-Fej\'er means on  martingale Hardy spaces}
The main result of this section reads:
\begin{theorem}
	\label{theorem1sigma}a) Let $f\in H_{1/2}.$ Then there exists an absolute constant $%
	c,$ such that
	\begin{equation*}
	\sup_{n\in \mathbb{N}}\frac{1}{n\log n}\overset{n}{\underset{k=1}{\sum }}%
	\left\Vert\sigma_{k}f\right\Vert_{1/2}^{1/2}\leq c \left\Vert f\right\Vert _{H_{1/2}}^{1/2}. 
	\end{equation*}
	
	b) Let $\varphi :\mathbb{N}_{+}\rightarrow \lbrack 1,$ $\infty )$ be a nondecreasing
	function satisfying the condition
	\begin{equation}
	\overline{\lim_{n\rightarrow \infty }}\frac{\log n}{\varphi _{n}}=+\infty .
	\label{cond1a}
	\end{equation}
	
	Then there exists a function $f\in H_{1/2},$ such that
	\begin{equation*}
	\sup_{n\in \mathbb{N}}\frac{1}{n\varphi_n}\overset{n}{\underset{k=1}{\sum }}\left\Vert\sigma_{k}f\right\Vert _{H_{1/2}}^{1/2}=\infty .
	\end{equation*}
\end{theorem}

\begin{corollary}
	\label{theorem1a}There exists a martingale $f\in H_{1/2},$ such that 
	\begin{equation*}
	\sup_{n\in \mathbb{N}}\frac{1}{n}\overset{n}{\underset{k=1}{\sum }}%
	\left\Vert \sigma _{k}f\right\Vert _{1/2}^{1/2}=\infty .
	\end{equation*}
\end{corollary}

\qquad
\textbf{Proof:}
By using \ref{strongsgn} was proved that there exists an absolute constant $c$, such that
\begin{equation*}
\left\Vert\sigma_{k}f\right\Vert_{H_{1/2}}^{1/2}\leq c \log k \left\Vert f\right\Vert _{H_{1/2}}^{1/2}, \ \ k=1,2,... 
\end{equation*} 
Hence,
\begin{equation*}
\frac{1}{n\log n}\overset{n}{\underset{k=1}{\sum }}%
\left\Vert \sigma_{k}f\right\Vert _{H_{1/2}}^{1/2}\leq  \frac{c\left\Vert f\right\Vert_{H_{1/2}}^{1/2}} {n\log n}\overset{n}{\underset{k=1}{\sum }}{\log k}\leq c\left\Vert f\right\Vert_{H_{1/2}}^{1/2}.
\end{equation*}

The proof of part a) is complete.

Under the condition \eqref{cond1a} there exists an increasing sequence of the positive integers $\left\{\alpha_{k}:k\in \mathbb{N}\right\} $ such that
\begin{equation*}
\overline{\lim_{k\rightarrow \infty }}\frac{\log M_{{\alpha _{k}}}}{\varphi
	_{2M_{\alpha _{k}}}}=+\infty
\end{equation*}%
and
\begin{equation} \label{69}
\sum_{k=0}^{\infty }\frac{\varphi _{2M_{\alpha _{k}}}^{1/2}}{\log ^{1/2}M_{{%
			\alpha _{k}}}}<c<\infty .  
\end{equation}

Let $f=( f^{(n)},\ \ n\in \mathbb{N}) $ be martingale, defined by 
\begin{equation*}
f^{(n)} :=\sum_{\left\{ k;\text{ }2\alpha_{k}<n\right\} }\lambda _{k}a_{k},
\end{equation*}
where
\begin{equation*}
a_{k}=M_{\alpha_{k}}r_{\alpha _{k}}D_{M_{_{\alpha _{k}}}}=M_{\alpha_{k}}(D_{2M_{_{\alpha
			_{k}}}}-D_{M_{_{\alpha _{k}}}})
\end{equation*}%
and
\begin{equation*}
\lambda _{k}=\frac{\varphi _{2M_{\alpha _{k}}}}{\log M_{{\alpha
			_{k}}}}.
\end{equation*}

Since$\ $%
\begin{equation}
S_{2^{A}}a_{k}=\left\{ 
\begin{array}{ll}
a_{k}, & \alpha _{k} <A, \\ 
0, &  \alpha _{k} \geq A,%
\end{array}%
\right.  \label{4aa}
\end{equation}%
\begin{equation*}
\text{supp}(a_{k})=I_{\alpha _{k} },\text{\ }%
\int_{I_{\alpha _{k} }}a_{k}d\mu =0,\text{\ }\left\Vert
a_{k}\right\Vert _{\infty }\leq M^2_{\alpha _{k} }=\mu (%
\text{supp }a_{k})^{-2},
\end{equation*}%
if we apply Lemma \ref{lemma2.1} and (\ref{69}) we conclude that $f\in H_{1/2}.$

Moreover,
\begin{equation} \label{6}
\widehat{f}(j)=\left\{
\begin{array}{l}
M_{\alpha_{k}}{\lambda _{k}},\,\,\text{\ \ \ \ \thinspace \thinspace }j\in \left\{
M_{\alpha _{k}},...,2M_{\alpha _{k}}-1\right\} ,\text{ }k\in \mathbb{N} \\
0\text{ },\text{ \thinspace \qquad \thinspace \thinspace \thinspace
	\thinspace \thinspace }j\notin \bigcup\limits_{k=1}^{\infty }\left\{
M_{\alpha _{k}},...,2M_{\alpha _{k}}-1\right\} .\text{ }%
\end{array}%
\right.   
\end{equation}

We have that 
\begin{eqnarray}  \label{7}
\sigma _{n}f &=&\frac{1}{n}\sum_{j=0}^{M_{\alpha _{k}}-1}S_{j}f+\frac{1}{n}%
\sum_{j=M_{\alpha _{k}}}^{n-1}S_{j}f \\
&=&I+II.  \notag
\end{eqnarray}

Let $M_{\alpha _{k}}\leq j<2M_{\alpha _{k}}.$ Since
\begin{equation*}
D_{j+M_{\alpha _{k}}}=D_{M_{\alpha _{k}}}+\psi _{_{M_{\alpha _{k}}}}D_{j},%
\text{ \qquad when \thinspace \thinspace }j\leq M_{\alpha _{k}},
\end{equation*}
if we apply (\ref{6}) we obtain that
\begin{eqnarray}  \label{8sas}
S_{j}f &=&S_{M_{\alpha _{k}}}f+\sum_{v=M_{\alpha _{k}}}^{j-1}\widehat{f}%
(v)\psi _{v}  \\
&=&S_{M_{\alpha _{k}}}f+M_{\alpha_{k}}{\lambda _{k}}\sum_{v=M_{\alpha _{k}}}^{j-1}\psi _{v}
\notag \\
&=&S_{M_{\alpha _{k}}}f+M_{\alpha_{k}}{\lambda _{k}}\left( D_{j}-D_{M_{\alpha _{k}}}\right)
\notag \\
&=&S_{M_{\alpha _{k}}}f+{\lambda _{k}}\psi _{M_{\alpha _{k}}}D_{j-M_{\alpha
		_{k}}}  \notag
\end{eqnarray}

According to (\ref{8sas}) concerning $II$ we conclude can that 
\begin{eqnarray*}
	II &=&\frac{n-M_{\alpha _{k}}}{n}S_{M_{\alpha_k}}f \\
	&&+\frac{\lambda _{k}M_{\alpha _{k}}}{ n}\sum_{j=M_{2\alpha
			_{k}}}^{n-1}\psi_{M_{\alpha_{k}}}D_{j-M_{\alpha_{k}}} \\
	&&=II_{1}+II_{2}.
\end{eqnarray*}

We can estimate $II_{2}$ as fallows:
\begin{eqnarray*}
	\left\vert II_{2}\right\vert &=&\frac{\lambda _{k}M_{\alpha _{k}}}{n}\left\vert \psi _{M_{\alpha _{k}}}\sum_{j=0}^{n-M_{\alpha
			_{k}}-1}D_{j}\right\vert \\
	&=&\frac{\lambda _{k}M_{\alpha _{k}}}{n}{(n-M_{\alpha _{k}})}%
	\left\vert K_{n-M_{\alpha _{k}}}\right\vert \\
	&\geq &{\lambda_{k}}\left( n-M_{\alpha _{k}}\right)
	\left\vert K_{n-M_{\alpha _{k}}}\right\vert.
\end{eqnarray*}

Let $n=\sum_{i=1}^{s}\sum_{k=l_{i}}^{m_{i}}M_{k},$ \ where 
\begin{equation*}
0\leq l_{1}\leq m_{1}\leq l_{2}-2<l_{2}\leq m_{2}\leq ...\leq
l_{s}-2<l_{s}\leq m_{s}.
\end{equation*}

By applying Lemma \ref{lemma5} we get that 
\begin{eqnarray*}
	\left\vert II_{2}\right\vert &\geq &{c\lambda _{k}\left\vert \left( n-M_{\alpha
			_{k}}\right) K_{n-M_{\alpha _{k}}}\left( x\right) \right\vert } \\
	&\geq &{c\lambda _{k}M_{l_{i}}^{2}},\text{ \ \ for \ \ }x\in
	I_{l_{i}+1}\left( e_{l_{i}-1}+e_{l_{i}}\right) .
\end{eqnarray*}

Hence 
\begin{eqnarray}\label{8aaa0} 
&&\int_{G_{m}}\left\vert II_{2}\right\vert ^{1/2}d\mu  \\  \notag
&\geq &\sum_{i=1}^{s-1}\int_{I_{l_{i}+1}\left( e_{l_{i}-1}+e_{l_{i}}\right)
}\left\vert II_{2}\right\vert ^{1/2}d\mu   \\ \notag
&\geq &c\sum_{i=1}^{s-1}\int_{I_{l_{i}+1}\left( e_{l_{i}-1}+e_{l_{i}}\right)
}{\lambda _{k}^{1/2}M_{l_{i}}}d\mu   \\ \notag
&\geq &{c\lambda _{k}^{1/2}\left( s-1\right) } \\
&\geq&{c\lambda _{k}^{1/2}v\left( n-M_{\alpha _{k}}\right) }.  \notag
\end{eqnarray}

In view of  (\ref{smn}), (\ref{sigmamn}) and (\ref{7})  we find that 
\begin{eqnarray}  \label{8aaa}
\left\Vert I\right\Vert ^{1/2}
&=&\left\Vert \frac{M_{\alpha _{k}}}{n}%
\sigma_{M_{\alpha _{k}}}f\right\Vert _{1/2}^{1/2} \\ \notag
&\leq & \left\Vert \sigma_{M_{\alpha _{k}}}f\right\Vert _{1/2}^{1/2} \\ \notag
&\leq& c\left\Vert f\right\Vert_{H_{1/2}}^{1/2} 
\end{eqnarray}
and 
\begin{eqnarray}  \label{8bbb}
\left\Vert II_{1}\right\Vert ^{1/2} 
&=&\left\Vert \frac{n-M_{\alpha _{k}}}{n}S_{M_{\alpha _{k}}}f\right\Vert
_{1/2}^{1/2}  \\ \notag
&\leq&  \left\Vert S_{M_{\alpha _{k}}}f\right\Vert _{1/2}^{1/2}   \\ \notag
&\leq&  c\left\Vert f\right\Vert_{H_{1/2}}^{1/2}.  
\end{eqnarray}

By combining (\ref{8aaa0}), (\ref{8aaa}) and (\ref{8bbb}) we get that 
\begin{eqnarray*}
	&&\left\Vert \sigma _{n}f\right\Vert _{1/2}^{1/2} \\
	&\geq &\left\Vert II_{2}\right\Vert _{1/2}^{1/2}-\left\Vert
	II_{1}\right\Vert _{1/2}^{1/2}-\left\Vert I\right\Vert _{1/2}^{1/2} \\
	&\geq &{c\lambda _{k}^{1/2}v\left( {n-M_{\alpha _{k}}}\right) }%
	-c\left\Vert f\right\Vert _{H_{1/2}}^{1/2}.
\end{eqnarray*}

By using estimates with the above we can conclude that 
\begin{eqnarray}
&&\underset{n\in \mathbf{\mathbb{N}}_{+}}{\sup }\frac{1}{n\varphi_n}\underset{k=1}{%
	\overset{n}{\sum }}\left\Vert \sigma _{k}f\right\Vert _{1/2}^{1/2}
\label{8bbb100} \\
&\geq &\frac{1}{M_{\alpha _{k}+1}\varphi_{2M_{\alpha _{k}}}}\underset{\left\{ M_{\alpha _{k}}\leq
	l\leq 2M_{\alpha _{k}}\right\} }{\sum }\left\Vert \sigma _{l}f\right\Vert
_{1/2}^{1/2}  \notag \\
&\geq &\frac{c}{M_{\alpha _{k}+1}\varphi_{2M_{\alpha _{k}}}}\underset{\left\{ M_{\alpha _{k}}\leq
	l\leq 2M_{\alpha _{k}}\right\} }{\sum }\left( {\lambda _{k}^{1/2}v\left( l-M_{\alpha
		_{k}}\right) }-c\left\Vert f\right\Vert
_{H_{1/2}}^{1/2}\right)  \notag \\ \notag
&\geq &\frac{c\lambda _{k}^{1/2}}{M_{\alpha _{k}}\varphi_{2M_{\alpha _{k}}}}\underset{l=1}{\overset{
		M_{\alpha _{k}}}{\sum }}v\left( l\right) \\ \notag
&-&\frac{c\left\Vert f\right\Vert _{H_{1/2}}^{1/2}}{M_{\alpha _{k}}\varphi_{2M_{\alpha _{k}}}}%
\underset{\left\{ M_{\alpha _{k}}\leq l\leq 2M_{\alpha _{k}}\right\} }{\sum }%
1   \\ \notag
&\geq &\frac{c\lambda _{k}^{1/2}}{M_{\alpha _{k}}\varphi_{2M_{\alpha _{k}}}}\underset{l=1}{\overset{
		M_{\alpha _{k}}-1}{\sum }}v\left( l\right)- c \\ \notag
&	\geq &c\frac{\log^{1/2} M_{{\alpha _{k}}}}{\varphi_{2M_{\alpha _{k}}}^{1/2}}\rightarrow \infty ,\text{ as \ }k\rightarrow
\infty .  \notag
\end{eqnarray}

The proof is complete.

\QED

\subsection{Convergence of subsequences of Vilenkin-Fej\'er means on the martingale Hardy spaces}

Our main result of this section reads:
\begin{theorem}\label{theorem1sub} 
	a) Let $0<p<1/2,$ $f\in H_{p}.$ Then there exists an
	absolute constant $c_{p}$, depending only on $p$, such that 
	\begin{equation*}
	\text{ }\left\Vert \sigma _{n_{k}}f\right\Vert _{H_{p}}\leq \frac{%
		c_{p}M_{\left\vert n_{k}\right\vert }^{1/p-2}}{M_{\left\langle
			n_{k}\right\rangle }^{1/p-2}}\left\Vert f\right\Vert _{H_{p}}.
	\end{equation*}
	
	\textit{b) (sharpness) Let }$0<p<1/2$ \textit{and} $\Phi \left( n\right) $ \textit{be
		any nondecreasing function,\ such that } 
	\begin{equation} \label{31aaa}
	\sup_{k\in \mathbb{N}}\rho\left( n_{k}\right) =\infty ,\text{ \ \ }\overline{\underset{%
			k\rightarrow \infty }{\lim }}\frac{M_{\left\vert n_{k}\right\vert }^{1/p-2}}{%
		M_{\left\langle n_{k}\right\rangle}^{1/p-2}\Phi\left( n_{k}\right) }
	=\infty.  
	\end{equation}%
	\textit{Then there exists a martingale }$f\in H_{p},$\textit{\ such that} 
	\begin{equation*}
	\underset{k\in \mathbb{N}}{\sup }\left\Vert \frac{\sigma _{n_{k}}f}{\Phi \left(
		n_{k}\right) }\right\Vert _{weak-L_{p}}=\infty .
	\end{equation*}
\end{theorem}

\textbf{Proof:} by using  \eqref{knbounded} we obtain that 
$$\frac{M_{\left\langle n_{k}\right\rangle
	}^{1/p-2}\left\vert \sigma _{n_{k}}a\left( x\right) \right\vert }{%
	M_{\left\vert n_{k}\right\vert }^{1/p-2}}$$   
is bounded from $L_{\infty}$ to $L_{\infty}.$ According to Lemma \ref{lemma2.2} we find that the proof of
Theorem will be complete, if we show that%
\begin{equation*}
\int_{\overline{I_{N}}}\left\vert \frac{M_{\left\langle n_{k}\right\rangle
	}^{1/p-2}\sigma _{n_{k}}a\left( x\right) }{M_{\left\vert n_{k}\right\vert
	}^{1/p-2}}\right\vert ^{p}<c<\infty ,
\end{equation*}%
for every $p$-atom $a,$ with support$\ I$ and $\mu \left( I\right)
=M_{N}^{-1}.$ We may assume that $I=I_{N}.$ It is easy to see that $\sigma
_{n_{k}}\left( a\right) =0$ when $n_{k}\leq M_{N}.$ Therefore, we can
suppose that $n_{k}>M_{N}$.

Since $\left\Vert a\right\Vert _{\infty }\leq M_{N}^{1/p}$ we find that
\begin{eqnarray} \label{400}
&&\frac{M_{\left\langle n_{k}\right\rangle }^{1/p-2}\left\vert \sigma
	_{n_{k}}a\left( x\right) \right\vert }{M_{\left\vert n_{k}\right\vert
	}^{1/p-2}}\\ \notag
&\leq& \frac{M_{\left\langle n_{k}\right\rangle }^{1/p-2}}{%
	M_{\left\vert n_{k}\right\vert }^{1/p-2}}\int_{I_{N}}\left\vert a\left(
t\right) \right\vert \left\vert K_{n_{k}}\left( x-t\right) \right\vert d\mu
\left( t\right)   \label{400a} \\
&\leq &\frac{M_{\left\langle n_{k}\right\rangle }^{1/p-2}\left\Vert
	a\right\Vert _{\infty }}{M_{\left\vert n_{k}\right\vert }^{1/p-2}}%
\int_{I_{N}}\left\vert K_{n_{k}}\left( x-t\right) \right\vert d\mu \left(
t\right)   \notag \\
&\leq &\frac{M_{\left\langle n_{k}\right\rangle }^{1/p-2}M_{N}^{1/p}}{%
	M_{\left\vert n_{k}\right\vert }^{1/p-2}}\int_{I_{N}}\left\vert
K_{n_{k}}\left( x-t\right) \right\vert d\mu \left( t\right)   \notag \\
&\leq &M_{\left\langle n_{k}\right\rangle }^{1/p-2}M_{\left\vert
	n_{k}\right\vert }^{2}\int_{I_{N}}\left\vert K_{n_{k}}\left( x-t\right)
\right\vert d\mu \left( t\right).  \notag
\end{eqnarray}

Without loss the generality we may assume that $i<j$. Let $x\in I_{N}^{i,j}\ 
$and $j<\left\langle n_{k}\right\rangle .$ Then $x-t\in I_{N}^{i,j}$ for $%
t\in I_{N}$ and, according to \eqref{kn8}, we obtain that 

\begin{equation} \label{kn0}
\left\vert K_{M_{l}}\left( x-t\right) \right\vert =0,\text{ \ for all }
\left\langle n_{k}\right\rangle \leq \text{ }l\leq \left\vert
n_{k}\right\vert.
\end{equation}

By combining (\ref{fn5}) in Lemma \ref{lemma7kn} with (\ref{400}) and \eqref{kn0}, for $ x \in I_{N}^{i,j},\text{ \ }0\leq i<j<\left\langle
n_{k}\right\rangle $ we can conclude that 

\begin{eqnarray} \label{401} 
&&\frac{M_{\left\langle n_{k}\right\rangle }^{1/p-2}\left\vert \sigma
	_{n_{k}}a\left( x\right) \right\vert }{M_{\left\vert n_{k}\right\vert
	}^{1/p-2}} \\ \notag
&\leq& M_{\left\langle n_{k}\right\rangle }^{1/p-2}M_{\left\vert
	n_{k}\right\vert }^{2}\overset{\left\vert n_{k}\right\vert }{\underset{%
		l=\left\langle n_{k}\right\rangle }{\sum }}\int_{I_{N}}\left\vert
K_{M_{l}}\left( x-t\right) \right\vert d\mu \left( t\right) =0.\text{ }
\end{eqnarray}

Let $x\in I_{N}^{i,j},\,$where $\left\langle n_{k}\right\rangle \leq j\leq N.
$ Then, in the view of Lemma \ref{lemma5}, we have that 
\begin{equation*}
\int_{I_{N}}\left\vert K_{n_{k}}\left( x-t\right) \right\vert d\mu \left(
t\right) \leq \frac{cM_{i}M_{j}}{M_{N}^{2}}.
\end{equation*}

By using again (\ref{400}) we find that 
\begin{eqnarray}  \label{402}
&&\frac{M_{\left\langle n_{k}\right\rangle }^{1/p-2}\left\vert \sigma
	_{n_{k}}a\left( x\right) \right\vert }{M_{\left\vert n_{k}\right\vert
	}^{1/p-2}} \\ \notag
&\leq &\frac{M_{\left\langle n_{k}\right\rangle
	}^{1/p-2}M_{N}^{1/p}}{M_{\left\vert n_{k}\right\vert }^{1/p-2}}%
\int_{I_{N}}\left\vert K_{n_{k}}\left( x-t\right) \right\vert d\mu \left(
t\right)   \\ \notag
&\leq &\frac{M_{\left\langle n_{k}\right\rangle }^{1/p-2}M_{N}^{1/p}}{%
	M_{\left\vert n_{k}\right\vert }^{1/p-2}}\frac{M_{i}M_{j}}{M_{N}^{2}}\\ \notag
&\leq&
M_{\left\langle n_{k}\right\rangle }^{1/p-2}M_{i}M_{j}.  \notag
\end{eqnarray}

By using (\ref{1.1}) we get that 
\begin{eqnarray*}
	&&\int_{\overline{I_{N}}}\left\vert \frac{M_{\left\langle n_{k}\right\rangle
		}^{1/p-2}\left\vert \sigma _{n_{k}}a\left( x\right) \right\vert }{%
		M_{\left\vert n_{k}\right\vert }^{1/p-2}}\right\vert ^{p}d\mu  \\
	&=&\overset{N-2}{\underset{i=0}{\sum }}\overset{N-1}{\underset{j=i+1}{\sum }}%
	\int_{I_{N}^{i,j}}\left\vert \frac{M_{\left\langle n_{k}\right\rangle
		}^{1/p-2}\left\vert \sigma _{n_{k}}a\left( x\right) \right\vert }{%
		M_{\left\vert n_{k}\right\vert }^{1/p-2}}\right\vert ^{p}d\mu \\ &+&\overset{N-1}{%
		\underset{i=0}{\sum }}\int_{I_{N}^{k,N}}\left\vert \frac{M_{\left\langle
			n_{k}\right\rangle }^{1/p-2}\left\vert \sigma _{n_{k}}a\left( x\right)
		\right\vert }{M_{\left\vert n_{k}\right\vert }^{1/p-2}}\right\vert ^{p}d\mu 
	\\
	&\leq &\overset{\left\langle n_{k}\right\rangle -1}{\underset{i=0}{\sum }}%
	\overset{N-1}{\underset{j=\left\langle n_{k}\right\rangle }{\sum }}%
	\int_{I_{N}^{i,j}}\left\vert \frac{M_{\left\langle n_{k}\right\rangle
		}^{1/p-2}\left\vert \sigma _{n_{k}}a\left( x\right) \right\vert }{%
		M_{\left\vert n_{k}\right\vert }^{1/p-2}}\right\vert ^{p}d\mu  \\
	&&+\overset{N-2}{\underset{i=\left\langle n_{k}\right\rangle }{\sum }}%
	\overset{N-1}{\underset{j=i+1}{\sum }}\int_{I_{N}^{i,j}}\left\vert \frac{%
		M_{\left\langle n_{k}\right\rangle }^{1/p-2}\left\vert \sigma
		_{n_{k}}a\left( x\right) \right\vert }{M_{\left\vert n_{k}\right\vert
		}^{1/p-2}}\right\vert ^{p}d\mu  \\
	&&+\overset{N-1}{\underset{i=0}{\sum }}\int_{I_{N}^{i,N}}\left\vert \frac{%
		M_{\left\langle n_{k}\right\rangle }^{1/p-2}\left\vert \sigma
		_{n_{k}}a\left( x\right) \right\vert }{M_{\left\vert n_{k}\right\vert
		}^{1/p-2}}\right\vert ^{p}d\mu  \\
	&\leq &\overset{\left\langle n_{k}\right\rangle -1}{\underset{i=0}{\sum }}%
	\overset{N-1}{\underset{j=\left\langle n_{k}\right\rangle }{\sum }}%
	\int_{I_{N}^{i,j}}\left\vert M_{\left\langle n_{k}\right\rangle
	}^{1/p-2}M_{i}M_{j}\right\vert ^{p}d\mu \\
	&+&\overset{N-2}{\underset{%
			i=\left\langle n_{k}\right\rangle }{\sum }}\overset{N-1}{\underset{j=i+1}{%
			\sum }}\int_{I_{N}^{i,j}}\left\vert M_{\left\langle n_{k}\right\rangle
	}^{1/p-2}M_{i}M_{j}\right\vert ^{p}d\mu  \\
	&&+\overset{N-1}{\underset{i=0}{\sum }}\int_{I_{N}^{i,N}}\left\vert
	M_{\left\langle n_{k}\right\rangle }^{1/p-2}M_{i}M_{N}\right\vert ^{p}d\mu 
\end{eqnarray*}

Hence, by combining (\ref{400}-\ref{402}) we get that 

\begin{eqnarray*}
	&&\int_{\overline{I_{N}}}\left\vert \frac{M_{\left\langle n_{k}\right\rangle
		}^{1/p-2}\left\vert \sigma _{n_{k}}a\left( x\right) \right\vert }{%
		M_{\left\vert n_{k}\right\vert }^{1/p-2}}\right\vert ^{p}d\mu  \\
	&\leq &c_{p}M_{\left\langle n_{k}\right\rangle }^{1-2p}\overset{\left\langle
		n_{k}\right\rangle -1}{\underset{i=0}{\sum }}\overset{N-1}{\underset{%
			j=\left\langle n_{k}\right\rangle }{\sum }}\frac{\left( M_{i}M_{j}\right)
		^{p}}{M_{j}} \\
	&+&c_{p}M_{\left\langle n_{k}\right\rangle }^{1-2p}\overset{N-2}{%
		\underset{i=\left\langle n_{k}\right\rangle }{\sum }}\overset{N-1}{\underset{%
			j=i+1}{\sum }}\frac{\left( M_{i}M_{j}\right) ^{p}}{M_{j}} \\
	&+&c_{p}M_{\left\langle n_{k}\right\rangle }^{1-2p}\underset{i=0}{\sum }%
	\frac{\left( M_{i}M_{N}\right) ^{p}}{M_{N}} \\
	&\leq &c_{p}M_{\left\langle n_{k}\right\rangle }^{1-2p}\overset{\left\langle
		n_{k}\right\rangle }{\underset{i=0}{\sum }}M_{i}^{p}\overset{N-1}{\underset{%
			j=\left\langle n_{k}\right\rangle +1}{\sum }}\frac{1}{M_{j}^{1-p}} \\
	&+&M_{\left\langle n_{k}\right\rangle }^{1-2p}\overset{N-2}{\underset{%
			i=\left\langle n_{k}\right\rangle }{\sum }}M_{i}^{p}\overset{N-1}{\underset{%
			j=i+1}{\sum }}\frac{1}{M_{j}^{1-p}} \\
	&&+c_{p}\overset{N-1}{\underset{i=0}{\sum }}\frac{M_{i}^{p}}{M_{N}^{p}} \\
	&\leq &c_{p}M_{\left\langle n_{k}\right\rangle }^{1-2p}M_{\left\langle
		n_{k}\right\rangle }^{p}\frac{1}{M_{\left\langle n_{k}\right\rangle }^{1-p}}\\
	&+&c_{p}M_{\left\langle n_{k}\right\rangle }^{1-2p}\overset{N-2}{\underset{%
			i=\left\langle n_{k}\right\rangle }{\sum }}\frac{1}{M_{i}^{1-2p}}+c_{p}\leq
	c_{p}<\infty .
\end{eqnarray*}

The proof of the a) part is complete.

b) Let 
$\left\{ n_{k}:k\geq 0\right\} $ be a sequence of positive numbers,
satisfying condition (\ref{31aaa}). Then 
\begin{equation}
\sup_{k\in \mathbb{N}}\frac{M_{\left\vert n_{k}\right\vert }}{M_{\left\langle
		n_{k}\right\rangle }}=\infty .  \label{12h}
\end{equation}

Under condition (\ref{12h}) there exists a sequence $\left\{ \alpha _{k}:%
\text{ }k\geq 0\right\} \subset \left\{ n_{k}:\text{ }k\geq 0\right\} $ such
that $\alpha _{0}\geq 3$ and 
\begin{equation}
\sum_{k=0}^{\infty }\frac{M_{\left\langle \alpha _{k}\right\rangle }^{\left(
		1-2p\right) /2}\Phi^{p/2}\left(\alpha_{k}\right)}{M_{\left\vert \alpha _{k}\right\vert }^{\left( 1-2p\right)
		/2}}<c<\infty .  \label{12hh}
\end{equation}

Let \qquad 
\begin{equation*}
f^{\left( n\right) }=\sum_{\left\{ k;\text{ }\left\vert \alpha
	_{k}\right\vert <n\right\} }\lambda _{k}a_{k},
\end{equation*}%
where 
\begin{equation*}
\lambda _{k}=\frac{\lambda M_{\left\langle \alpha _{k}\right\rangle
	}^{\left( 1/p-2\right) /2}\Phi^{1/2}\left(\alpha_{k}\right)}{M_{\left\vert \alpha _{k}\right\vert }^{\left(
		1/p-2\right) /2}}
\end{equation*}%
and%
\begin{equation*}
a_{k}=\frac{M_{\left\vert \alpha _{k}\right\vert }^{1/p-1}}{\lambda }\left(
D_{M_{\left\vert \alpha _{k}\right\vert +1}}-D_{M_{\left\vert \alpha
		_{k}\right\vert }}\right) .
\end{equation*}

B applying Lemma \ref{lemma2.1} we can conclude that $f\in H_{p}.$

It is evident that
\begin{equation}
\widehat{f}(j)=\left\{ 
\begin{array}{l}
M_{\left\vert \alpha _{k}\right\vert }^{1/2p}M_{\left\langle \alpha
	_{k}\right\rangle }^{\left( 1/p-2\right) /2}\Phi^{1/2}\left(\alpha_{k}\right),\,\,\text{ } \\ 
\text{if \thinspace \thinspace }j\in \left\{ M_{\left\vert \alpha
	_{k}\right\vert },...,\text{ ~}M_{\left\vert \alpha _{k}\right\vert
	+1}-1\right\} ,\text{ }k=0,1,2..., \\ 
0\text{ },\text{ \thinspace \qquad \thinspace\ \ \ \ \thinspace\ \ \ \ \ }
\\ 
\text{\ if \thinspace \thinspace \thinspace }j\notin
\bigcup\limits_{k=0}^{\infty }\left\{ M_{\left\vert \alpha _{k}\right\vert
},...,\text{ ~}M_{\left\vert \alpha _{k}\right\vert +1}-1\right\} .\text{ }%
\end{array}%
\right.  \label{6aacharp}
\end{equation}

Moreover,%
\begin{eqnarray*}
	\frac{\sigma _{_{\alpha _{k}}}f}{\Phi\left(\alpha_{k}\right)}=\frac{1}{\alpha_{k}\Phi\left(\alpha_{k}\right)}\sum_{j=1}^{M_{\left\vert
			\alpha _{k}\right\vert }}S_{j}f+\frac{1}{\alpha _{k}\Phi\left(\alpha_{k}\right)}\sum_{j=M_{\left\vert
			\alpha _{k}\right\vert }+1}^{\alpha _{k}}S_{j}f 
	&:=&I+II.
\end{eqnarray*}%
Let $M_{\left\vert \alpha _{k}\right\vert }<j\leq \alpha _{k}.$ Then, by
applying (\ref{6aacharp}) we get that 
\begin{equation}
S_{j}f=S_{M_{\left\vert \alpha _{k}\right\vert }}f+M_{\left\vert \alpha
	_{k}\right\vert }^{1/2p}M_{\left\langle \alpha _{k}\right\rangle }^{\left(
	1/p-2\right)/2}\Phi^{1/2}\left(\alpha_{k}\right)\left( D_{j}-D_{M_{\left\vert \alpha _{k}\right\vert
}}\right).  \label{8aafn}
\end{equation}

By using (\ref{8aafn}) we can rewrite $II$ as
\begin{eqnarray*}
	II &=&\frac{\alpha _{k}-M_{\left\vert \alpha _{k}\right\vert }}{\alpha _{k}\Phi\left(\alpha_{k}\right)}%
	S_{M_{\left\vert \alpha _{k}\right\vert }}f\\
	&+&\frac{M_{\left\vert \alpha
			_{k}\right\vert }^{1/2p}M_{\left\langle \alpha _{k}\right\rangle }^{\left(
			1/p-2\right)/2}}{\alpha _{k}\Phi^{1/2}\left(\alpha_{k}\right)}\sum_{j=M_{\left\vert \alpha _{k}\right\vert
	}}^{\alpha _{k}}\left( D_j-D_{M_{\left\vert \alpha _{k}\right\vert
	}}\right) \\
	&:=&II_{1}+II_{2}.
\end{eqnarray*}

If we combine \eqref{smn} and \eqref{sigmamn} it is easy to show that%
\begin{eqnarray*}
	\left\Vert II_{1}\right\Vert _{weak-L_p}^p 
	&\leq&\left(\frac{\alpha_k-M_{\left\vert \alpha _k\right\vert }}{\alpha _k\Phi\left(\alpha_{k}\right)}\right)^{p}\left\Vert S_{M_{\left\vert \alpha _k\right\vert }}f\right\Vert_{weak-L_p}^{p} \\
	&\leq&\left(\frac{\alpha
		_{k}-M_{\left\vert \alpha _{k}\right\vert }}{\alpha _{k}\Phi\left(\alpha_{k}\right)}\right)
	^{p}\left\Vert S_{M_{\left\vert \alpha _{k}\right\vert }}f\right\Vert
	_{p}^{p} \\
	&\leq& c_p\left\Vert f\right\Vert _{H_p}^p<\infty .
\end{eqnarray*}
and
\begin{eqnarray*}
	\left\Vert I\right\Vert _{weak-L_{p}}^{p}&=&\left( \frac{M_{\left\vert \alpha
			_{k}\right\vert }}{\alpha _{k}\Phi\left(\alpha_{k}\right)}\right) ^{p}\left\Vert \sigma _{M_{\left\vert \alpha _{k}\right\vert }}f\right\Vert _{weak-L_p}^p\\
	&\leq&\left( \frac{M_{\left\vert \alpha_{k}\right\vert }}{\alpha _{k}\Phi\left(\alpha_{k}\right)}\right) ^{p}\left\Vert \sigma _{M_{\left\vert \alpha _{k}\right\vert}}f\right\Vert_p^p\\
	&\leq& c_p\left\Vert f\right\Vert _{H_p}^p<\infty .
\end{eqnarray*}

Let $x\in $ $I_{_{\left\langle \alpha _{k}\right\rangle +1}}^{\left\langle
	\alpha _{k}\right\rangle -1,\left\langle \alpha _{k}\right\rangle }.$ Under
condition (\ref{31aaa}) we can conclude that $\left\langle \alpha
_{k}\right\rangle \neq \left\vert \alpha _{k}\right\vert $ and $\left\langle
\alpha _{k}-M_{\left\vert \alpha _{k}\right\vert }\right\rangle
=\left\langle \alpha _{k}\right\rangle .$ Since 
\begin{equation}
D_{j+M_{n}}=D_{M_{n}}+\psi _{M_{n}}D_{j}=D_{M_{n}}+r_{n}D_{j},\text{ when }%
\,\,j<M_{n}  \label{8k}
\end{equation}
if we apply estimate Lemma \ref{lemma8ccc} for $II_{2}$ we obtain
that 
\begin{eqnarray*}
	\left\vert II_{2}\right\vert &=&\frac{M_{\left\vert \alpha _{k}\right\vert
		}^{1/2p}M_{\left\langle \alpha _{k}\right\rangle }^{\left( 1/p-2\right) /2}}{
		\alpha _{k}\Phi^{1/2}\left(\alpha_{k}\right)}\left\vert \sum_{j=1}^{\alpha _{k}-M_{\left\vert \alpha
			_{k}\right\vert }}\left( D_{j+M_{\left\vert \alpha _{k}\right\vert
	}}-D_{M_{\left\vert \alpha _{k}\right\vert }}\right) \right\vert \\
	&=&\frac{M_{\left\vert \alpha _{k}\right\vert }^{1/2p}M_{\left\langle \alpha
			_{k}\right\rangle }^{\left( 1/p-2\right) /2}}{\alpha _{k}\Phi^{1/2}\left(\alpha_{k}\right)}\left\vert \psi
	_{M_{\left\vert \alpha _{k}\right\vert }}\sum_{j=1}^{\alpha
		_{k}-M_{\left\vert \alpha _{k}\right\vert }}D_{j}\right\vert \\
	&\geq &\frac{c_{p}M_{\left\vert \alpha _{k}\right\vert }^{1/2p-1}M_{\left\langle
			\alpha _{k}\right\rangle }^{\left( 1/p-2\right) /2}}{\Phi^{1/2}\left(\alpha_{k}\right)}\left( \alpha
	_{k}-M_{\left\vert \alpha _{k}\right\vert }\right) \left\vert K_{\alpha
		_{k}-M_{\left\vert \alpha _{k}\right\vert }}\right\vert \\
	&\geq &\frac{c_{p}M_{\left\vert \alpha _{k}\right\vert }^{1/2p-1}M_{\left\langle
			\alpha _{k}\right\rangle }^{\left( 1/p+2\right) /2}}{\Phi^{1/2}\left(\alpha_{k}\right)}.
\end{eqnarray*}

It follows that%
\begin{eqnarray*}
	&&\left\Vert II_{2}\right\Vert _{weak-L_{p}}^{p} \\
	&\geq &c_{p}\left( \frac{M_{\left\vert \alpha _{k}\right\vert }^{\left(
			1/p-2\right) /2}M_{\left\langle \alpha _{k}\right\rangle }^{\left(
			1/p+2\right) /2}}{\Phi^{1/2}\left(\alpha_{k}\right)}\right) ^{p}\mu \left\{ x\in G_{m}:\left\vert
	IV_{2}\right\vert \geq c_{p}M_{\left\vert \alpha _{k}\right\vert }^{\left(
		1/p-2\right) /2}M_{\left\langle \alpha _{k}\right\rangle }^{\left(
		1/p+2\right) /2}\right\} \\
	&\geq &c_{p}\frac{M_{\left\vert \alpha _{k}\right\vert }^{1/2-p}M_{\left\langle
			\alpha _{k}\right\rangle }^{1/2+p}\mu \left\{I_{_{\left\langle \alpha
				_{k}\right\rangle +1}}^{\left\langle \alpha _{k}\right\rangle
			-1,\left\langle \alpha _{k}\right\rangle }\right\}}{\Phi^{p/2}\left(\alpha_{k}\right)} \\
	&\geq& \frac{c_{p}M_{\left\vert \alpha _{k}\right\vert}^{1/2-p}}{M_{\left\langle \alpha
			_{k}\right\rangle }^{1/2-p}\Phi^{p/2}\left(\alpha_{k}\right)}.
\end{eqnarray*}
Hence, for large $k$, 
\begin{eqnarray*}
	&&\left\Vert \sigma _{\alpha _{k}}f\right\Vert _{weak-L_{p}}^{p} \\
	&\geq &\left\Vert II_{2}\right\Vert _{weak-L_{p}}^{p}-\left\Vert
	II_{1}\right\Vert _{weak-L_{p}}^{p}-\left\Vert I\right\Vert _{weak-L_{p}}^{p}
	\\
	&\geq &\frac{1}{2}\left\Vert II_{2}\right\Vert _{weak-L_{p}}^{p}\\
	&\geq& \frac{c_{p}M_{\left\vert \alpha _{k}\right\vert }^{1/2-p}}{2M_{\left\langle \alpha_{k}\right\rangle }^{1/2-p}\Phi^{p/2}\left(\alpha_{k}\right)}\rightarrow \infty ,\text{ as }k\rightarrow
	\infty.
\end{eqnarray*}

The proof is complete.

\QED

\begin{corollary}
	\label{corollary1sub} Let $0<p<1/2,$ and  $f\in H_{p}.$ Then there exists
	an absolute constant $c_{p}$, depending only on $p$, such that 
	\begin{equation*}
	\text{ }\left\Vert \sigma _{n_k}f\right\Vert _{H_{p}}\leq c_{p}\left\Vert
	f\right\Vert _{H_{p}},\text{ \ \ }k\in \mathbb{N}
	\end{equation*}
	if and only if when
	\begin{equation*}
	\sup_{k\in \mathbb{N}}\rho\left( n_{k}\right)<c<\infty.
	\end{equation*}
\end{corollary}

As a application we also obtain the previous mentioned result by Weisz \cite%
{We1}, \cite{We3} (Theorem W).

\begin{corollary}\label{corollary2sub} 
	Let $0<p<1/2,$ \  $f\in H_{p}.$ Then there exists
	an absolute constant $c_{p}$, depending only on $p$, such that 
	\begin{equation*}
	\text{ }\left\Vert \sigma _{M_{n}}f\right\Vert _{H_{p}}\leq c_{p}\left\Vert
	f\right\Vert _{H_{p}},\text{ \ \ }n\in \mathbb{N}.
	\end{equation*}
\end{corollary}
On the other hand, the following unexpected result is true:
\begin{corollary}
	\label{corollary3sub} \textit{a)}Let $0<p<1/2,$ \  $f\in H_{p}.$ Then there exists an absolute constant $c_{p}$, depending only on $p$, such that 
	\begin{equation*}
	\text{ }\left\Vert \sigma _{M_{n}+1}f\right\Vert _{H_{p}}\leq c_{p}M^{1/p-2}_{n}\left\Vert
	f\right\Vert _{H_{p}},\text{ \ \ }n\in \mathbb{N}.
	\end{equation*}
	\textit{b)} Let $0<p<1/2$ and $\Phi \left( n\right) $ be any nondecreasing function, such that
	\begin{equation*} 
	\overline{\underset{	k\rightarrow \infty }{\lim }}\frac{M_{k }^{1/p-2}}{\Phi\left( k\right)}
	=\infty.  
	\end{equation*}%
	Then there exists a martingale  $f\in H_{p},$\textit{\ such that} 
	\begin{equation*}
	\underset{k\in \mathbb{N}}{\sup }\left\Vert \frac{\sigma _{M_{k}+1}f}{\Phi \left(k\right)}\right\Vert_{weak-L_p}=\infty .
	\end{equation*}
\end{corollary}
\begin{remark}
	From Corollary \ref{corollary2sub} we obtain that $ \sigma _{M_{n}} f$ are bounded from $ H_p $ to $ H_p $, but from Corollary \ref{corollary3sub} we conclude that $ \sigma _{M_{n}+1} f$ are not bounded from $ H_p $ to $ H_p $. The main reason is that Fourier coefficients of martingale $ f\in H_p $ are not uniformly bounded (for details see e.g. \cite{tep8}).
\end{remark}

We also state the following corollary which shows difference of the rate of divergence across to  different subsequences.

\begin{corollary}
	\label{corollary30sub} \textit{a)}Let $0<p<1/2,$ \  $f\in H_{p}.$ Then there exists an absolute constant $c_{p}$, depending only on $p$, such that 
	\begin{equation*}
	\left\Vert \sigma _{M_{n}+M_{[n/2]}}f\right\Vert _{H_{p}}\leq c_{p}\left(M_{n}/M_{[n/2]}\right)^{1/p-2}\left\Vert
	f\right\Vert _{H_{p}},\text{ \ \ }n\in \mathbb{N},
	\end{equation*}
	where $[n/2]$ denotes integer part of $n/2$.
	
	\textit{b)} Let $0<p<1/2$ and $\Phi \left( n\right) $ be any non-decreasing function, such that
	\begin{equation*} 
	\overline{\underset{k\rightarrow \infty }{\lim }}\frac{{\left(M_k/M_{[k/2]} \right)}^{1/p-2}}{\Phi\left(k\right)}=\infty.  
	\end{equation*}%
	Then there exists a martingale  $f\in H_{p},$ such that
	\begin{equation*}
	\underset{k\in \mathbb{N}}{\sup }\left\Vert \frac{\sigma _{M_{k}+/M_{[k/2]}}f}{\Phi \left(k\right)}\right\Vert_{weak-L_p}=\infty .
	\end{equation*}
\end{corollary}

In the next corollary we state Corollaries \ref{corollary3sub} and  \ref{corollary30sub} for Walsh system only to clearly see difference of divergence rates for the various subsequences:
\begin{corollary}
	\label{corollary4sub} a) Let $0<p<1/2,$ \  $f\in H_{p}.$ Then there exists an absolute constant $c_{p}$, depending only on $p$, such that 
	\begin{equation} \label{aaaaa}
	\text{ }\left\Vert \sigma^w _{2^{n}+1}f\right\Vert _{H_{p}}\leq c_{p}2^{(1/p-2)n}\left\Vert
	f\right\Vert _{H_{p}},\text{ \ \ }n\in \mathbb{N}
	\end{equation}
	and
	\begin{equation} \label{aaa1aa}
	\text{ }\left\Vert \sigma^w _{2^{n}+2^{[n/2]}}f\right\Vert _{H_{p}}\leq c_{p}2^{\frac{(1/p-2)n}{2}}\left\Vert
	f\right\Vert _{H_{p}},\text{ \ \ }n\in \mathbb{N},
	\end{equation}
	where $[n/2]$ denotes integer part of $n/2$.
	
	\textit{b)} The rates $2^{(1/p-2)n} $ and  $ 2^{\frac{(1/p-2)n}{2}}$ in inequalities (\ref{aaaaa}) and (\ref{aaa1aa}) are sharp in the same sense as in Theorem \ref{theorem1sub}.
\end{corollary}

\subsection{necessary and sufficient condition for the norm convergence of subsequences of partial sums in terms of modulus of continuity}

\begin{theorem}
	\label{theorem1}a) Let $0<p<1/2,$ $f\in H_{p}$, $ \sup_{k\in \mathbb{N}} d\left(n_{k}\right) =\infty $ and 
	\begin{equation}
	\omega _{H_{p}}\left( 1/M_{\left\vert n_{k}\right\vert},f\right) =o\left(\frac{M_{\left\langle n_{k}\right\rangle}^{1/p-2}}{M_{\left\vert n_{k}\right\vert }^{1/p-2}}\right) ,\text{ as \ }	k\rightarrow \infty .  \label{cond3}
	\end{equation}
	Then 
	\begin{equation}
	\left\Vert \sigma _{n_{k}}f-f\right\Vert_{H_p}\rightarrow 0,\text{ \ \ as \ }k\rightarrow \infty .  \label{fe2}
	\end{equation}
	
	b) Let $\sup_{k\in \mathbb{N}}\rho\left( n_{k}\right) =\infty .$ Then there exists a
	martingale $f\in H_{p}(G)$ $\left( 0<p<1/2\right) ,$\ \ for which 
	\begin{equation}
	\omega _{H_{p}}\left( 1/M_{\left\vert n_{k}\right\vert },f\right) =O\left(\frac{M_{\left\langle n_{k}\right\rangle}^{1/p-2}}{M_{\left\vert n_{k}\right\vert }^{1/p-2}}\right) ,\text{ \ as \ }
	k\rightarrow \infty  \label{cond4a}
	\end{equation}%
	and 
	\begin{equation}
	\left\Vert \sigma _{n_{k}}f-f\right\Vert _{weak-L_p}\nrightarrow
	0,\,\,\,\text{as\thinspace \thinspace \thinspace }k\rightarrow \infty .
	\label{kn4}
	\end{equation}
\end{theorem}

\textbf{Proof:}
Let $0<p<1/2, \ f\in H_{p}$ and $M_{k}<n\leq M_{k+1}.$ By applying part a) of Theorem \ref{theorem1sub} we can conclude that
\begin{equation*}
\left\Vert \sigma _{n}f-f\right\Vert _{H_p}^{p}
\end{equation*}
\begin{equation*}
\leq \left\Vert \sigma_{n}f-\sigma_{n}S_{M_{k}}f\right\Vert
_{H_p}^{p}+\left\Vert\sigma_{n}S_{M_{k}}f-S_{M_{k}}f\right\Vert
_{H_p}^{p}+\left\Vert S_{M_{k}}f-f\right\Vert_{H_p}^{p}
\end{equation*}%
\begin{equation*}
=\left\Vert \sigma _{n}\left( S_{M_{k}}f-f\right) \right\Vert
_{H_p}^{p}+\left\Vert S_{M_{k}}f-f\right\Vert _{H_p}^{p}+\left\Vert
\sigma _{n}S_{M_{k}}f-S_{M_{k}}f\right\Vert _{H_p}^{p}
\end{equation*}%
\begin{equation*}
\leq c_{p}\left( \frac{M_{\left\vert n\right\vert }^{1-2p}}{M_{\left\langle n\right\rangle}^{1-2p}} +1\right) \omega_{H_{p}}^{p}\left( 1/M_{n},f\right) +\left\Vert \sigma_{n}S_{M_{k}}f-S_{M_{k}}f\right\Vert_{H_p}^{p}.
\end{equation*}

By simple calculation we get that
\begin{eqnarray} \label{knsn1*}
&&\sigma _{n}S_{M_{N}}f-S_{M_{N}}f
\\ \notag
&=&\frac{1}{n}\sum_{k=0}^{M_{N}}S_{k}S_{M_{N}}f+\frac{1}{n}
\sum_{k=M_{N}+1}^{n}S_{k}S_{M_{N}}f-S_{M_{N}}f
\\ \notag
&=&\frac{1}{n}\sum_{k=0}^{M_{N}}S_{k}f+\frac{1}{n}%
\sum_{k=M_{N}+1}^{n}S_{M_{N}}f-S_{M_{N}}f
\\ \notag
&=&\frac{1}{n}\sum_{k=0}^{M_{N}}S_{k}f+\frac{n-M_{N}}{n}S_{M_{N}}f-S_{M_{N}}f
\\ \notag
&=&\frac{M_{N}}{n}\sigma _{M_{N}}f-\frac{M_{N}}{n}S_{M_{N}}f
\\ \notag
&=&\frac{M_{N}}{n}\left( S_{M_{N}}\sigma _{M_{N}}f-S_{M_{N}}f\right)
\\ \notag
&=&\frac{M_{N}}{n}S_{M_{N}}\left(\sigma _{M_{N}}f-f\right).
\end{eqnarray}

Let $p>0.$ By combining \eqref{smn} and \eqref{sigmamn} we can conclude that
\begin{equation} \label{qwer}
\left\Vert \sigma _{n}S_{M_{k}}f-S_{M_{k}}f\right\Vert_{H_p}^{p} 
\end{equation}
\begin{equation*}
\leq \frac{2^{M_{k}}}{n^{p}}\left\Vert S_{M_{k}}\left( \sigma
_{M_{k}}f-f\right) \right\Vert_{H_p}^{p}\leq c_p\left\Vert \sigma
_{M_{k}}f-f\right\Vert _{H_{p}}^{p}\rightarrow 0,\text{ as }k\rightarrow
\infty \text{.}
\end{equation*}
On the other hand, under the condition (\ref{cond3}) we also get that
\begin{equation}
c_{p}\left( \frac{M_{\left\vert n\right\vert }^{1-2p}}{M_{\left\langle n\right\rangle}^{1-2p}} +1\right) \omega_{H_{p}}^{p}\left( 1/M_{n},f\right) \rightarrow 0 \label{asd}
\end{equation}
by combining (\ref{qwer}) and (\ref{asd}) we complete the proof of theorem.

Now, prove second part of theorem. Since $\sup_{k\in\mathbb{N}}\rho\left(
n_{k}\right) =\infty ,$ we obtain that, for $0<p<1/2,$ 
\begin{equation*}
\frac{M_{\left<n_k\right>}^{1/p -2}}{M_{\left|n_k\right|}^{1/p -2}}=\frac{1} {\left(m_{\left<n_k\right>} \ldots m_{\left|n_k\right|-1} \right)^{1/p -2}}\leq  \frac{1}{2^{\rho\left(n_k\right)\left(1/p -2\right)}}\rightarrow 0.
\end{equation*}
It follows that there exists $\{\alpha _{k}:k\geq 1\}\subset
\{n_{k}:k\geq 1\}$ such that $\sup_{k\in\mathbb{N}}\rho\left( \alpha _{k}\right) =\infty $
and
\begin{equation} \label{17aa}
\frac{M_{\left<\alpha_{k}\right>}^{1/p -2}}{M_{\left|\alpha_{k}\right|}^{1/p -2}}\leq \left(\frac{M_{\left<\alpha_{k-1}\right>}^{1/p -2}}{M_{\left|\alpha_{k-1}\right|}^{1/p -2}}\right)^2 \quad \text{for all} \quad k\in\mathbb{N}.
\end{equation}

By using \eqref{17aa} we get that
\begin{equation*}
\frac{{M_{\left<\alpha_k\right>}^{1/p-2}}}{M_{\left|\alpha_k\right|}^{1/p-2}}\leq \left(\frac{{M_{\left<\alpha_{k-1}\right>}^{1/p-2}}}{M_{\left|\alpha_{k-1}\right|}^{1/p-2}}\right)^2\leq\ldots \leq\left(\frac{{M_{\left<\alpha_{0}\right>}^{1/p-2}}}{M_{\left|\alpha_{0}\right|}^{1/p-2}}\right)^{k+1}\leq\frac{1}{2^{(k+1)\left(\left|\alpha_0\right|-\left<\alpha_0\right>\right)(1/p-2)}}
\end{equation*}
and
\begin{equation} \label{17a}
\underset{k=0}{\overset{\infty}{\sum }}\left(\frac{{M_{\left<\alpha_k\right>}^{1/p-2}}}{M_{\left|\alpha_k\right|}^{1/p-2}}\right)^p
\leq\underset{k=0}{\overset{\infty}{\sum }}\frac{1}{2^{(k+1)\left(\left|\alpha_0\right|-\left<\alpha_0\right>\right)(1-2p)}}<c<\infty.
\end{equation}
We set $f=\left(f^{\left(n\right)}, n\in\mathbb{N}\right)$ where
\begin{equation*}
f^{\left(n\right)}=\sum_{\left\{ i:\text{ }\left\vert \alpha_{i}\right\vert <n\right\}
}\frac{\lambda{M_{\left<\alpha_{i}\right> }^{\left(1/p-2\right)}}}{{M_{\left|\alpha_{i}\right|}^{\left(1/p-2\right)}}}a_{i}^{\left(p\right)},
\end{equation*}

$$a_{k}^{(p)}:=\frac {M_{\left\vert \alpha_{k}\right\vert }^{1/p-1}}{\lambda}\left(
D_{M_{\left\vert \alpha _{k}\right\vert +1}}-D_{M_{\left\vert \alpha_{k}\right\vert }}\right)$$
and 
$$\lambda=\sup_{k\in \mathbb{N}} m_k.$$
Since $a_{i}^{(p)}(x)$ is $p$-atom if we apply Lemma \ref{lemma2.1}  and (\ref{17a}) we conclude that $f\in H_{p}.$ On the other hand, if we apply Remark \ref{lemma2.3.6} we find that
\begin{equation*}
f-S_{M_{\left\vert \alpha_{_{n}}\right\vert }}f=\left( 0,...,0,\underset{i=n}{	\overset{n+s}{\sum }}\frac{{M_{\left<\alpha_i\right>}^{1/p-2}}}{M_{\left|\alpha_i\right|}^{1/p -2}}a_{i}^{\left(p\right)} ,...\right) ,\text{ \ }s\in \mathbb{N}_{+} 
\end{equation*}
is martingale. By combining (\ref{17aa}) and Lemma \ref{lemma2.1} we get that%
\begin{equation*}
\omega _{H_{p}}(1/M_{\left\vert \alpha_{_{n}}\right\vert },f) 
\end{equation*}%
\begin{equation*}
\leq \sum\limits_{i=n}^{\infty }\frac{{M_{\left<\alpha_i\right>}^{1/p -2}}}{{M_{\left|\alpha_i\right|}^{1/p -2}}}
\leq \sum\limits_{i=1}^{\infty }\left(\frac{{M_{\left<\alpha_n\right>}^{1/p -2}}}{{M_{\left|\alpha_n\right|}^{1/p -2}}}\right)^i
=O\left( \frac{{M_{\left<\alpha_n\right>}^{1/p -2}}}{{M_{\left|\alpha_{_n}\right|}^{1/p -2}}}\right), \ \ \text{as} \ \ n\rightarrow \infty.
\end{equation*}

It is easy to show that%
\begin{equation}
\widehat{f}(j)=\left\{ 
\begin{array}{ll}
M_{\left|\alpha_k\right|} M_{\left<\alpha_k\right>}^{1/p -2}, & \text{\thinspace
	\thinspace }j\in \left\{ M_{\left\vert \alpha _{k}\right\vert
},...,M_{_{\left\vert \alpha _{k}\right\vert +1}}-1\right\} ,\text{ }
k=0,1,... \\ 
0\,, & \text{\thinspace }j\notin \bigcup\limits_{i=0}^{\infty }\left\{
M_{_{\left\vert \alpha _{k}\right\vert }},...,M_{_{\left\vert \alpha
		_{k}\right\vert +1}}-1\right\} .\text{ }
\end{array}%
\right.  \label{4.22}
\end{equation}

Let $M_{\left\vert \alpha _{k}\right\vert }<j<\alpha _{k}.$ By using (\ref{4.22}) we get that%

\begin{eqnarray*}
	S_{j}f&=&S_{M_{\left\vert \alpha _{k}\right\vert }}f+\sum_{v=M{_{\left\vert
				\alpha _{k}\right\vert }}}^{j-1}\widehat{f}(v)w_{v} \\ \notag
	&=&S_{M_{\left\vert \alpha
			_{k}\right\vert }}f+M_{\left|\alpha_k\right|} M_{\left<\alpha_k\right>}^{1/p -2}\left( D_j-D_{M_{\left\vert \alpha_{k}\right\vert }}\right) 
\end{eqnarray*}

Hence,
\begin{equation}
\sigma _{_{\alpha_{k}}}f=\frac{1}{\alpha _{k}}\sum_{j=1}^{M_{\left\vert \alpha_{k}\right\vert }}S_{j}f+\frac{1}{\alpha _{k}}\sum_{j=M_{\left\vert \alpha_{k}\right\vert }+1}^{\alpha _{k}}S_{j}f
\label{7aaa2}
\end{equation}%
\begin{eqnarray*}
	&=&\frac{M_{\left|\alpha_k\right|}}{\alpha_k}\sigma _{_{M_{\left\vert \alpha _{k}\right\vert }}}f+\frac{\left( \alpha_{k}-M_{\left\vert	\alpha_{k}\right\vert }\right) S_{M_{\left\vert \alpha _{k}\right\vert }}f}{ \alpha _{k}} \\
	&+&\frac{M_{\left|\alpha_k\right|} M_{\left<\alpha_k\right>}^{1/p -2)}}{\alpha _{k}}\sum_{j=M_{_{\left\vert \alpha _{k}\right\vert }}+1}^{\alpha _{k}}\left( D_{_{j}}-D_{M_{\left\vert \alpha_{k}\right\vert }}\right)\\
	&=& I+II+III.
\end{eqnarray*}%

Since%
\begin{equation*}
D_{j+M_{n}}=D_{M_{n}}+\psi_{M_{n}}D_{j},\text{ \thinspace when \thinspace }j<M_{n+1}  \label{f1}
\end{equation*}%
we obtain that%
\begin{eqnarray} \label{9aaa}
\left\vert III\right\vert &=&\frac{M_{\left|\alpha_k\right|} M_{\left<\alpha_k\right>}^{1/p -2}}{\alpha _{k}}\left\vert \sum_{j=1}^{\alpha_{k}-M_{{\left\vert \alpha _{k}\right\vert }}}\left(D_{j+M_{{\left\vert \alpha_{k}\right\vert }}}-D_{M_{\left\vert \alpha_{k}\right\vert }}\right)	\right\vert \\ \notag
&=&\frac{M_{\left|\alpha_k\right|} M_{\left<n_k\right>}^{1/p -2}}{\alpha _{k}}\left\vert \sum_{j=1}^{\alpha _{k}-M_{\left\vert \alpha_{k}\right\vert }}D_{j}\right\vert \\ \notag
&=&\frac{M_{\left|\alpha_k\right|} M_{\left<\alpha_k\right>}^{1/p-2}}{\alpha _{k}}\left( \alpha_{k}-M{_{\left\vert \alpha_{k}\right\vert }}\right) \left\vert K_{\alpha_{k}-M_{_{\left\vert \alpha _{k}\right\vert }}}\right\vert \\ \notag
&\geq & c M_{\left<\alpha_k\right>}^{1/p -2}\left( \alpha _{k}-M_{_{\left\vert \alpha _{k}\right\vert }}\right)
\left\vert K_{\alpha _{k}-M_{_{\left\vert \alpha _{k}\right\vert
}}}\right\vert.
\end{eqnarray}%

By combining \eqref{7aaa2} and \eqref{9aaa} we can conclude that
\begin{eqnarray*}
	&&\left\Vert \sigma _{\alpha _{k}}f-f\right\Vert _{weak-L_p}^{p}=\Vert I+II+III-f\Vert _{weak-L_p}^{p} \\
	&=&\Vert III+\frac{M_{\left|\alpha_k\right|}}{\alpha_k}\sigma _{_{M_{\left\vert \alpha _{k}\right\vert }}}f+\frac{\left( \alpha_{k}-M_{\left\vert	\alpha_{k}\right\vert }\right) S_{M_{\left\vert \alpha _{k}\right\vert }}f}{ \alpha _{k}}-f\Vert _{weak-L_p}^{p}\\
	&=&\Vert III+\frac{M_{\left|\alpha_k\right|}}{\alpha_k}\left(\sigma _{_{M_{\left\vert \alpha _{k}\right\vert }}}f-f\right)+\frac{ \alpha_{k}-M_{\left\vert\alpha_{k}\right\vert }}{ \alpha _{k}}\left(S_{M_{\left\vert \alpha _{k}\right\vert }}f-f\right)\Vert _{weak-L_p}^{p} \\
	&\geq &\Vert III \Vert _{weak-L_p}^{p}-\left( \frac{M_{\left\vert \alpha _{k}\right\vert }}{\alpha _{k}}\right)^{p}\Vert \sigma _{M_{_{\left\vert \alpha _{k}\right\vert }}}f-f\Vert_{weak-L_p}^{p} \\
	&-&\left( \frac{\alpha _{k}-M_{_{\left\vert\alpha_{k}\right\vert }}}{\alpha _{k}}\right) ^{p}\Vert S_{M_{_{\left\vert \alpha_{k}\right\vert }}}f-f\Vert _{weak-L_p}^{p} \\
	&\geq& \Vert III\Vert _{weak-L_p}^{p}-\Vert \sigma _{M_{_{\left\vert \alpha _{k}\right\vert }}}f-f\Vert
	_{weak-L_p}^{p}-\Vert S_{M_{_{\left\vert \alpha_{k}\right\vert }}}f-f\Vert _{weak-L_p}^{p}.
\end{eqnarray*}

If we combine \eqref{smn} and \eqref{sigmamn} it is easy to show that%
\begin{eqnarray*}
	\Vert \sigma _{M_{_{\left\vert \alpha _{k}\right\vert }}}f-f\Vert
	_{weak-L_p}^{p}\rightarrow 0,\text{ \ \ as \ }k\rightarrow \infty ,   \\
	\Vert S_{M_{\left\vert \alpha_{k}\right\vert }}f-f\Vert _{weak-L_p}^{p}\rightarrow 0,\text{ \ \ as \ }k\rightarrow \infty , 
\end{eqnarray*}

Hence,	for sufficiently large $k,$ we can write that
\begin{eqnarray*}
	&&\left\Vert \sigma _{\alpha _{k}}f-f\right\Vert _{weak-L_p}^{p}\\
	&\geq& \frac{1}{2}\Vert III\Vert _{weak-L_p}^{p}
	\geq \frac{M_{\left<\alpha_k\right>}^{1-2p}}{2}\Vert \left( \alpha_{k}-M_{\left\vert \alpha _{k}\right\vert }\right) K_{\alpha_{k}-M_{\left\vert \alpha _{k}\right\vert }}\Vert _{weak-L_p}^{p}
\end{eqnarray*}

Let $x\in E_{\left< \alpha _{k}\right> }.$ By using Lemma 
\ref{lemma8ccc} we have that
\begin{eqnarray*}
	&&\mu \left\{ x\in G_m:\left( \alpha_k-M_{\left\vert \alpha_k\right\vert
	}\right) \left\vert K_{\alpha_k-M_{\left\vert \alpha_k\right\vert
	}}\right\vert \geq c M_{\left< \alpha_k\right>}^{2}\right\} \\
	&\geq& \mu \left(E_{\left< \alpha_k\right> }\right) \geq c/M_{\left< \alpha _k\right>},
\end{eqnarray*}
and
\begin{equation*}
cM_{\left|\alpha_k\right|}^{2p}\mu \left\{ x\in G:\left( \alpha
_k-M_{_{\left\vert \alpha_k\right\vert}}\right)\left\vert K_{\alpha
	_k-M_{\left\vert \alpha_k\right\vert}}\right\vert \geq cM_{\left| \alpha_k\right|}^{2p}\right\}\geq cM_{\left| \alpha_k\right|}^{ 2p-1}.
\end{equation*}

Hence, 
\begin{equation*}
\left\Vert \sigma _{\alpha_{k}}f-f\right\Vert _{weak-L_p}\nrightarrow
0,\,\,\,\text{as}\,\,\,k\rightarrow \infty
\end{equation*}%
and Theorem is proved.

\QED

On the other hand, the following unexpected new result is also obtained:
\begin{corollary}
	\label{corollary3} \textit{a)} Let $0<p<1/2,$ \  $f\in H_{p}$ and \begin{equation*}
	\omega _{H_{p}}\left( 1/M_{ n_k},f\right) =o\left(\frac{1}{M_{ n_{k} }^{1/p-2}}\right) ,\text{ as \ }	k\rightarrow \infty . 
	\end{equation*}
	Then 
	\begin{equation*}
	\left\Vert \sigma _{M_{n_k}+1}f-f\right\Vert_{H_p}\rightarrow 0,\text{ \ \ as \ }k\rightarrow \infty . 
	\end{equation*}
	
	b) Let $\sup_{k\in \mathbb{N}}\rho\left( n_{k}\right) =\infty .$ Then there exists a
	martingale $f\in H_{p}(G)$ $\left( 0<p<1/2\right) ,$\ \ for which 
	\begin{equation*}
	\omega _{H_{p}}\left( 1/M_{\left\vert n_{k}\right\vert },f\right) =O\left(\frac{1}{M_{ n_{k} }^{1/p-2}}\right) ,\text{ \ as \ }
	k\rightarrow \infty 
	\end{equation*}%
	and 
	\begin{equation*}
	\left\Vert \sigma _{{M_{n_k}+1}}f-f\right\Vert _{weak-L_p}\nrightarrow
	0,\,\,\,\text{as\thinspace \thinspace \thinspace }k\rightarrow \infty .
	\end{equation*}
\end{corollary}

\begin{corollary}
	\label{corollary3a} \textit{a)} Let $0<p<1/2,$ \  $f\in H_{p}$ and \begin{equation*}
	\omega _{H_{p}}\left( 1/M_{ n_k},f\right) =o\left(\frac{{M_{ [n_{k}/2] }^{1/p-2}}}{M_{ n_{k} }^{1/p-2}}\right) ,\text{ as \ }	k\rightarrow \infty .  
	\end{equation*}
	Then 
	\begin{equation*}
	\left\Vert \sigma _{M_{n_k}+M_{ [n_{k}/2] }}f-f\right\Vert_{H_p}\rightarrow 0,\text{ \ \ as \ }k\rightarrow \infty .  
	\end{equation*}
	
	b) Let $\sup_{k\in \mathbb{N}}\rho\left( n_{k}\right) =\infty .$ Then there exists a
	martingale $f\in H_{p}(G_m)$ $\left( 0<p<1/2\right) ,$\ \ for which 
	\begin{equation*}
	\omega _{H_{p}}\left( 1/M_{\left\vert n_{k}\right\vert },f\right) =O\left(\frac{{M_{ [n_{k}/2] }^{1/p-2}}}{M_{ n_{k} }^{1/p-2}}\right),
	\end{equation*}%
	and 
	\begin{equation*}
	\left\Vert \sigma _{M_{n_k}+M_{ [n_{k}/2] }}f-f\right\Vert _{weak-L_p}\nrightarrow
	0,\,\,\,\text{as\thinspace \thinspace \thinspace }k\rightarrow \infty .
	\end{equation*}
\end{corollary}

In the next corollary we state theorem for Walsh system only to clearly see difference of divergence rates for the various subsequences:
\begin{corollary}
	\label{corollary3} \textit{a)} Let $0<p<1/2,$ \  $f\in H_{p}$ and \begin{equation*}
	\omega _{H_{p}}\left( 1/2^{ n},f\right) =o\left(\frac{1}{2^{ n(1/p-2)}}\right) ,\text{ as \ }	n\rightarrow \infty . 
	\end{equation*}
	Then 
	\begin{equation*}
	\left\Vert \sigma^w _{2^{n}+1}f-f\right\Vert_{H_p}\rightarrow 0,\text{ \ \ as \ }n\rightarrow \infty . 
	\end{equation*}
	
	b) There exists a
	martingale $f\in H_{p}(G)$ $\left( 0<p<1/2\right) ,$\ \ for which 
	\begin{equation*}
	\omega _{H_{p}}\left( 1/2^{n},f\right) =O\left(\frac{1}{2^{ n(1/p-2)}}\right) ,
	\end{equation*}%
	and 
	\begin{equation*}
	\left\Vert \sigma^w _{2^n+1}f-f\right\Vert _{weak-L_p}\nrightarrow
	0,\,\,\,\text{as\thinspace \thinspace \thinspace }k\rightarrow \infty .
	\end{equation*}
\end{corollary}

\begin{corollary}
	\label{corollary3aa} \textit{a)} Let $0<p<1/2,$ \  $f\in H_{p}$ and \begin{equation*}
	\omega _{H_p}\left( 1/2^{ n},f\right) =o\left(\frac{1}{2^{ n(1/p-2)/2 }}\right) ,\text{ as \ }	n\rightarrow \infty.  
	\end{equation*}
	Then 
	\begin{equation*}
	\left\Vert \sigma^w _{2^{n}+2^{ [n/2] }}f-f\right\Vert_{H_p}\rightarrow 0,\text{ \ \ as \ }n\rightarrow \infty ,
	\end{equation*}
	where $[n/2]$ denotes integer part of $n/2$.
	
	b) There exists a
	martingale $f\in H_{p}(G_m)$ $\left( 0<p<1/2\right) ,$\ \ for which 
	\begin{equation*}
	\omega _{H_{p}}\left( 1/M_{ n },f\right) =O\left(\frac{1}{2^{ n(1/p-2)/2 }}\right)
	\end{equation*}%
	and 
	\begin{equation*}
	\left\Vert \sigma^w _{2^{n}+2^{ [n/2] }}f-f\right\Vert _{weak-L_p}\nrightarrow
	0,\,\,\,\text{as\thinspace \thinspace \thinspace }n\rightarrow \infty .
	\end{equation*}
\end{corollary}

\newpage

\section[$T$ means of Vilenkin-Fourier series on $H_p$ spaces]{$T$ means of Vilenkin-Fourier series on  martingale Hardy spaces}

\vspace{0.5cm}

\subsection{Some classical results on $T$ means of Vilenkin-Fourier series}

\quad \quad It is well-known in the literature that the so-called $T$ means are generalizations of the Fej\'er, Ces\`aro and logarithmic means. The $T$ summation is a general summability method. Therefore it is of prior
interest to study the behavior of operators related to N\"orlund means of
Fourier series with respect to orthonormal systems.

Since $T$ means are inverse of N\"orlund means we state some interesting results concerning  N\"orlund summability, which has high influence on the new results for $T$ means of Vilenkin-Fourier series.

In \cite{gog8} Goginava investigated the behavior of Ces\`aro means of Walsh-Fourier series in detail. In the two-dimensional case approximation properties of N\"orlund and Ces\`aro means were considered by Nagy (see \cite{na}, \cite{n} and \cite{nagy}). The maximal operator $\sigma ^{\alpha ,\ast }$ $\left( 0<\alpha <1\right) $ of the $\left( C,\alpha \right) $ means of Vilenkin systems was investigated by Weisz \cite{We6}. In this paper Weisz proved that $\sigma ^{\alpha ,\ast }$ is bounded from the martingale space $H_{p}$ to the Lebesgue space $L_{p}$ for $p>1/\left( 1+\alpha \right) .$ Goginava \cite{gog4} gave a counterexample which shows that boundedness does not hold for $0<p\leq 1/\left( 1+\alpha \right) .$ Weisz and Simon \cite{sw} showed that the maximal operator $\sigma ^{\alpha ,\ast }$ is bounded from the Hardy space $H_{1/\left( 1+\alpha \right) }$ to the space $weak-L_{1/\left( 1+\alpha \right)}$.

Strong convergence theorems and boundedness of weighted maximal operators of of the $\left( C,\alpha \right) $ means of
Vilenkin systems on the Hardy spaces, when $0<p\leq 1/( 1+\alpha) $ were considered by Blahota, Tephnadze \cite{bt3} and Blahota, Tephnadze, Toledo \cite{btt}. Summability of some general methods were considered by Blahota, Nagy and Tephnadze \cite{bnt1}, Weisz,

In Persson,Tephnadze and Wall \cite{ptw} (see also \cite{tep8}) considered 
maximal operator of the N\"orlund summation method (see (\ref{1.2})). In particular, the maximal operator $t^{\ast }$ of  the
summability method (\ref{1.2}) with non-decreasing sequence $\{q_{k}:k\in
\mathbb{N}\}$ is bounded from the Hardy space $H_{1/2}$ to the space $%
weak-L_{1/2}.$

Moreover,for any $0<p<1/2$ and non-decreasing sequence $\{q_{k}:k\in \mathbb{N}\}$ 
satisfying the condition
\begin{equation}
\frac{q_{0}}{Q_{n}}\geq \frac{c}{n},\text{ \ \ }\left( c>0\right) ,
\label{cond1nor}
\end{equation}%
there exists a martingale $f\in H_{p},$ such that
\begin{equation*}
\underset{n\in
	\mathbb{N}}{\sup }\left\Vert t_{n}f\right\Vert _{weak-L_{p}}=\infty .
\end{equation*}

In Persson, Tephnadze and Wall \cite{ptw2} was proved that if $0<p<1/2$ and the sequence $\left\{ q_{k}:k\in
\mathbb{N}\right\} $ be non-decreasing, then the maximal operator%
\begin{equation*}
\widetilde{t}_{p,1}^{\ast }f:=\sup_{n\in \mathbb{N}}\frac{\left\vert
	t_{n}f\right\vert }{\left( n+1\right) ^{1/p-2}}
\end{equation*}%
is bounded from the Hardy martingale space $H_{p}$ to the Lebesgue space $L_{p}.$

On the other hand, according the fact that Fej\'er means are examples of N\"orlund means with
non-decreasing sequence $\left\{ q_{k}:k\in \mathbb{N}\right\} $ we
immediately obtain that the
asymptotic behaviour of the sequence of weights
\begin{equation*}
\left\{ 1/\left( k+1\right) ^{1/p-2}:k\in \mathbb{N}\right\}
\end{equation*}
in N\"orlund means can not be improved.

Let the sequence $\left\{ q_{k}:k\in \mathbb{N}%
\right\} $ be non-decreasing. Then the maximal operator%
\begin{equation*}
\overset{\sim }{t}_{1}^{\ast }f:=\sup_{n\in \mathbb{N}}\frac{\left\vert
	t_{n}f\right\vert }{\log ^{2}\left( n+1\right) }
\end{equation*}%
is bounded from the Hardy space  $H_{1/2}$ to the Lebesgue space $L_{1/2}.$
On the other hand, according the fact that Fej\'er means are examples of N\"orlund means with non-decreasing sequence $\left\{ q_{k}:k\in \mathbb{N}\right\} $ we
immediately obtain  that the
asymptotic behaviour of the sequence of  weights
\begin{equation*}
\left\{ 1/\log ^{2}\left( n+1\right) :n\in \mathbb{N}\right\}
\end{equation*}
in N\"orlund means can not be improved.

Persson, Tephnade and Wall \cite{ptw} proved that for all N\"orlund means with
non-increasing sequence $\left\{ q_{k}:k\in \mathbb{N}\right\} $ there
exists a martingale $f\in H_{p}$ such that%
\begin{equation*}
\underset{n\in
	\mathbb{N}}{\sup }\left\Vert t_{n}f\right\Vert _{weak-L_{p}}=\infty .
\end{equation*}

It follows that for any $0<p<1/2\ \ $ and  N\"orlund means $t_{n}$ with
non-increasing sequence $\{q_{k}:k\in \mathbb{N}\},$  the maximal
operator $t^{\ast }$ is not bounded from the martingale Hardy space $H_{p}$
to the space $weak-L_{p},$ that is there exists a martingale $f\in H_{p},$
such that
\begin{equation*}
\underset{n\in\mathbb{N}}{\sup }\left\Vert t^{\ast }f\right\Vert _{weak-L_{p}}=\infty .
\end{equation*}

Persson, Tephnade and
Wall \cite{ptw} find necessary condition for the N\"orlund means with
non-increasing sequence $\left\{ q_{k}:k\in
\mathbb{N}
\right\} $, when $1/2\leq p<1.$ In particular, for  $0<p<1/\left(1+\alpha \right),$ $0<\alpha
\leq 1,$ and non-increasing sequence $\{q_{k}:k\in \mathbb{N}\}$ 
satisfying the condition
\begin{equation} \label{cond4}
\overline{\lim_{n\rightarrow \infty }}\frac{n^{\alpha }}{Q_{n}}=c>0,\text{ }%
0<\alpha \leq 1,  
\end{equation}%
there exists a martingale $f\in H_{p}$ such that
\begin{equation}
\underset{n\in\mathbb{N}}{\sup }\left\Vert t_{n}f\right\Vert _{weak-L_{p}}=\infty . \label{cond444}
\end{equation}

Moreover, for any non-increasing sequence $\{q_{k}:k\in
\mathbb{N}\}$ satisfying the condition 
\begin{equation}
\overline{\lim_{n\rightarrow \infty }}\frac{n^{\alpha }}{Q_{n}}=\infty ,%
\text{ \ \ \ }\left( 0<\alpha \leq 1\right),  \label{cond22}
\end{equation}
there exists an martingale $f\in H_{1/\left( 1+\alpha \right) },$ such
that
\begin{equation}
\underset{n\in
	\mathbb{N}
}{\sup }\left\Vert t_{n}f\right\Vert_{weak-L_{1/\left( 1+\alpha \right)
}}=\infty .\label{cond445}
\end{equation}

It follows that for any $0<p<1/\left( 1+\alpha \right) ,$ $0<\alpha \leq 1$ and non-increasing
sequence $\{q_{k}:k\in \mathbb{N}\}$  satisfying the condition (\ref{cond4}). Then there exists a
martingale $f\in H_{p}$ such that
\begin{equation*}
\left\Vert t^{\ast }f\right\Vert _{weak-L_{p}}=\infty .
\end{equation*}

Moreover, if $\{q_{k}:k\in\mathbb{N}\}$ be a non-increasing sequence satisfying the condition (\ref{cond22}), then there exists an martingale $f\in H_{1/\left( 1+\alpha \right) }$ such
that
\begin{equation*}
\left\Vert t^{\ast }f\right\Vert _{weak-L_{1/\left( 1+\alpha \right)
}}=\infty .
\end{equation*}

In \cite{mpt} it was proved that the maximal operator $t^{\ast }$ of the N\"orlund summability method with non-increasing sequence $\{q_{k}:k\in
\mathbb{N}\},$ satisfying the condition 
\begin{equation} \label{6a}
\frac{1}{Q_n}=O\left(\frac{1}{n^{\alpha}}\right),\text{ \ \ when \ \ }
n\rightarrow \infty
\end{equation}
and
\begin{equation} \label{7a}
q_n-q_{n+1}=O\left(\frac{1}{n^{2-\alpha}}\right) ,\text{ \ \ when \ \ }
n\rightarrow \infty,
\end{equation}%
is bounded from the
Hardy space $H_{1/\left( 1+\alpha \right) }$ to the space $weak-L_{1/\left(
	1+\alpha \right) },$ for $0<\alpha \leq 1.$

Moreover, for $0<\alpha \leq 1$ and non-increasing
sequence $\{q_{k}:k\in \mathbb{N}\}$  satisfying the conditions%
\begin{equation}
\overline{\lim_{n\rightarrow \infty }}\frac{n^{\alpha }}{Q_{n}}\geq
c_{\alpha }>0  \label{cond29}
\end{equation}%
and
\begin{equation}
\left\vert q_{n}-q_{n+1}\right\vert \geq c_{\alpha }n^{\alpha -2},\text{ \ }n\in
\mathbb{N}\text{.}  \label{cond30}
\end{equation}
there exists a martingale $f\in H_{1/\left( 1+\alpha \right) }$ such
that
\begin{equation*}
\underset{n\in
	\mathbb{N}
}{\sup }\left\Vert t_{n}f\right\Vert _{{1/\left( 1+\alpha \right)
}}=\infty .
\end{equation*}

In \cite{ptw} (see  also \cite{bt2}) was proved that if $f\in H_{p},$ where $0<p\,<1/\left(1+\alpha \right)$ for some $0<\alpha \leq 1,$ and $\{q_{k}:k\in \mathbb{N}\} $ be a sequence of non-increasing numbers satisfying conditions (\ref{6a}%
) and (\ref{7a})  the maximal operator%
\begin{equation*}
\overset{\sim }{t}_{p,\alpha }^{\ast }:=\frac{\left\vert t_{n}f\right\vert }{\left( n+1\right) ^{1/p-1-\alpha }}
\end{equation*}
is bounded from the martingale Hardy space $H_{p}$ to the Lebesgue space $L_{p}.$

Moreover, let $\left\{ \Phi _{n}:n\in \mathbb{N}_{+}\right\} $ be any
non-decreasing sequence, satisfying the condition
\begin{equation}
\overline{\lim_{n\rightarrow \infty }}\frac{\left( n+1\right) ^{1/p-1-\alpha
}}{\Phi _{n}}=+\infty,  \label{6lbbb}
\end{equation}
then there exists N\"orlund means with non-increasing sequence $\{q_{k}:k\in
\mathbb{N}\}$ satisfying the conditions (\ref{cond29}) and (\ref{cond30})  such
that
\begin{equation*}
\sup_{k\in \mathbb{N}}\frac{\left\Vert \frac{t_{M_{_{2n_{k}}}+1}f_{k}}{\Phi
		_{M_{_{2n_{k}}}+1}}\right\Vert _{weal-L_{p}}}{\left\Vert f_{k}\right\Vert
	_{H_{p}}}=\infty.
\end{equation*}

It follows that if $0<p<1/\left( 1+\alpha \right) $ and $%
f\in H_{p},$ then there exists an absolute constant $c_{p,\alpha },$
depending only on $p$ and $\alpha $, such that
\begin{equation*}
\left\Vert t_{n}f\right\Vert _{p}\leq c_{p,\alpha }\left( n+1\right)
^{1/p-1-\alpha }\left\Vert f\right\Vert _{H_{p}},\text{ \ }n\in \mathbb{N}%
_{+}.
\end{equation*}

On the other hand, let $\left\{ \Phi _{n}:n\in \mathbb{N}%
\right\} $ be any non-decreasing sequence satisfying the condition (\ref
{6lbbb}), then there exists a martingale $f\in H_{p}$ such that
\begin{equation*}
\sup_{n\in \mathbb{N}}\left\Vert \frac{t_{n}f}{\Phi _{n}}\right\Vert
_{weak-L_{p}}=\infty.
\end{equation*}

Moreover, let $\left\{ \Phi _{n}:n\in \mathbb{N}\right\} $ be any non-decreasing sequence satisfying the condition (\ref%
{6lbbb}), then the maximal operator
\begin{equation*}
\sup_{n\in \mathbb{N}}\frac{\left\vert t_{n}f\right\vert }{\Phi _{n}}
\end{equation*}%
is not bounded from the Hardy space  $H_{p}$ to the space $weak-L_{p}.$

In \cite{bpt1} (see  also \cite{bt2}) was proved if $f\in H_{1/(1+\alpha )},$ where $0<\alpha \leq 1$
and $\{q_{k}:k\in \mathbb{N}\}$ be a sequence of non-increasing numbers
satisfying the conditions (\ref{6a}) and (\ref{7a}), then there exists an absolute
constant $c_{\alpha }$ depending only on $\alpha$ such that the maximal
operator%
\begin{equation*}
\overset{\sim }{t}_{\alpha }^{\ast }:=\frac{\left\vert t_{n}f\right\vert }{%
	\log ^{1+\alpha }\left( n+1\right) }
\end{equation*}
is bounded from the martingale Hardy space $H_{1/(1+\alpha )}$ to the Lebesgue space $
L_{1/(1+\alpha )}.$

Moreover, if $\left\{ \Phi _{n}:n\in \mathbb{N}_{+}\right\} $ be any
non-decreasing sequence satisfying the condition%
\begin{equation}
\overline{\lim_{n\rightarrow \infty }}\frac{\log ^{1+\alpha }\left(
	n+1\right) }{\Phi _{n}}=+\infty,  \label{nom1}
\end{equation}%
then there exists N\"orlund means with non-increasing sequence $\{q_{k}:k\in
\mathbb{N}\}$ satisfying the conditions (\ref{cond29}) and (\ref{cond30}) such that
\begin{equation*}
\sup_{k\in \mathbb{N}}\frac{\left\Vert \sup_{n}\left\vert \frac{t_{n}f_{k}}{\Phi _{n}}\right\vert \right\Vert _{1/(1+\alpha )}}{\left\Vert f\right\Vert
	_{H_{1/(1+\alpha )}}}=\infty .
\end{equation*}

In Persson, Tephnadze and Wall \cite{ptw2} was proved that if $0<p<1/2$, $f\in H_{p}$ and the sequence $\left\{ q_{k}:k\in \mathbb{N}\right\} $ be non-decreasing, then there
exists an absolute constant $c_{p}$ depending only on $p$ such that
\begin{equation*}
\overset{\infty }{\underset{k=1}{\sum }}\frac{\left\Vert t_{k}f\right\Vert
	_{p}^{p}}{k^{2-2p}}\leq c_{p}\left\Vert f\right\Vert _{H_{p}}^{p}.
\end{equation*}

On the other hand, according the fact that Fej\'er means are examples of N\"orlund means with
non-decreasing sequence $\left\{ q_{k}:k\in \mathbb{N}\right\} $ we
immediately obtain  that
the asymptotic behaviour of the sequence of weights
\begin{equation*}
\left\{ 1/k^{2-2p}:k\in \mathbb{N}\right\}
\end{equation*}
in N\"orlund means  can not be improved.

In Persson, Tephnadze and Wall \cite{ptw2} was proved if $f\in H_{1/2}$ and the sequence $\left\{q_{k}:k\in \mathbb{N}\right\} $ be non-decreasing satisfying condition (\ref{fn01}), then there exists an absolute constant $c,$ such that
\begin{equation*}
\frac{1}{\log n}\overset{n}{\underset{k=1}{\sum }}\frac{\left\Vert
	t_{k}f\right\Vert _{1/2}^{1/2}}{k}\leq c\left\Vert f\right\Vert
_{H_{1/2}}^{1/2}.
\end{equation*}

In  Blahota and Tephnadze \cite{bt2} was investigated N\"orlund means with non-increasing sequence  $\{q_{k}:k\in \mathbb{N}\}$ in the case $0<p<1/\left( 1+\alpha \right)$ where $0<\alpha <1$. In particular, if $f\in H_{p},$ where $0<p\,<1/\left(1+\alpha \right) ,$ $0<\alpha \leq 1$ and $\{q_{k}:k\in \mathbb{N}\},$ be a
sequence of non-increasing numbers satisfying the conditions (\ref{6a}) and (\ref{7a}), then there exists an absolute constant $c_{\alpha ,p},$ depending
only on $\alpha $ and $p$ such that
\begin{equation*}
\overset{\infty }{\underset{k=1}{\sum }}\frac{\left\Vert t_{k}f\right\Vert
	_{H_{p}}^{p}}{k^{2-\left( 1+\alpha \right) p}}\leq c_{\alpha ,p}\left\Vert
f\right\Vert _{H_{p}}^{p}.
\end{equation*}

In Blahota, Persson and Tephnadze \cite{bpt1} was proved that if $f\in H_{1/(1+\alpha )}$ where $0<\alpha \leq 1$ and
$\{q_{k}:k\in \mathbb{N}\}$ be a sequence of non-increasing numbers
satisfying the conditions (\ref{6a}) and (\ref{7a}), then there exists an
absolute constant $c_{\alpha }$ depending only on $\alpha$ such that
\begin{equation*}
\frac{1}{\log n}\overset{n}{\underset{m=1}{\sum }}\frac{\left\Vert
	t_{m}f\right\Vert _{H_{1/\left( 1+\alpha \right) }}^{1/\left( 1+\alpha
		\right) }}{m}\leq c_{\alpha }\left\Vert f\right\Vert _{H_{1/\left( 1+\alpha
		\right) }}^{1/\left( 1+\alpha \right)}.
\end{equation*}

In \cite{tut3}  we investigate the maximal operator $T^{\ast }$ of \ the summability method (\ref{nor})
with non-increasing sequence $\{q_{k}:k\geq 0\},$ is bounded from the Hardy space $H_{1/2}$ to the space $weak-L_{1/2}.$

Moreover, for any $0<p<1/2$ and non-increasing sequence $\{q_{k}:k\geq 0\}$ satisfying the condition
\begin{equation*}
\frac{q_{n+1}}{Q_{n+2}}\geq \frac{c}{n},\text{ \ \ }\left( c\geq 1\right) .
\end{equation*}%
then there exists a martingale $f\in H_{p},$ such that
\begin{equation*}
\underset{n\in	\mathbb{N}}{\sup }\left\Vert T_{n}f\right\Vert _{weak-L_{p}}=\infty .
\end{equation*}

We also proved that the maximal operator $T^{\ast }$ of  the summability method (\ref{nor})
with non-decreasing sequence $\{q_{k}:k\geq 0\}$ satisfying the condition 
\begin{equation*}
\frac{q_{n-1}}{Q_n}= O\left(\frac{1}{n}\right)
\end{equation*}
is bounded from the Hardy space $H_{1/2}$ to the space $weak-L_{1/2}.$

Moreover, for any $0<p<1/2$ and non-decreasing sequence $\{q_{k}:k\geq 0\},$ there exists a martingale $f\in H_{p},$ such that
\begin{equation*}
\underset{n\in
	\mathbb{N}
}{\sup }\left\Vert T_{n}f\right\Vert _{weak-L_{p}}=\infty .
\end{equation*}

In \cite{tut4} we proved that for $0<p\leq 1/2,$ $f\in H_{p}$ and non-decreasing sequence  $\{q_{k}:k\geq
0\}$  the maximal operator 
\begin{equation} \label{Tn100}
\widetilde{T}_{p}^{\ast }f:=\sup_{n\in \mathbb{N}_{+}}\frac{\left\vert
	T_{n}f\right\vert }{\left( n+1\right) ^{1/p-2}\log ^{2\left[ 1/2+p\right]
	}\left( n+1\right) }
\end{equation}%
is bounded from the Hardy space $H_{p}$ to the space $L_{p}.$

On the other hand if $0<p\leq 1/2,$ $f\in H_{p}$ and $\{q_{k}:k\geq
0\}$ be a sequence of non-increasing numbers,satisfying the condition 
\begin{equation}
\frac{q_{n-1}}{Q_{n}}=O\left( \frac{1}{n}\right) ,\text{ \ \ as \ \ }\
n\rightarrow \infty .  \label{fn01}
\end{equation},
then the maximal operator 
\begin{equation} \label{Tn1000}
\widetilde{T}_{p}^{\ast }f:=\sup_{n\in \mathbb{N}_{+}}\frac{\left\vert
	T_{n}f\right\vert }{\left( n+1\right) ^{1/p-2}\log ^{2\left[ 1/2+p\right]
	}\left( n+1\right) }
\end{equation}%
is bounded from the Hardy space $H_{p}$ to the space $L_{p}.$

Since the maximal operator
\begin{equation*}
\widetilde{\sigma }_{p}^{\ast }f:=\sup_{n\in \mathbb{N}_{+}}\frac{\left\vert
	\sigma _{n}f\right\vert }{\left( n+1\right) ^{1/p-2}\log ^{2\left[ 1/2+p\right] }\left( n+1\right) }
\end{equation*}%
is bounded from the martingale Hardy space $H_{p}$ to the space $L_{p}$ and the rate of denominator $\left( n+1\right) ^{1/p-2}\log ^{2\left[ 1/2+p	\right] }$ is in the sense sharp and Fejer means is example of $T$ means non-decreasing sequence we obtain that this weights are also sharp in \eqref{Tn100} and \eqref{Tn1000}.

In \cite{tut4} we also investigate strong convergence of $T$ means with respect to Vilenkin systems. In particular, if $0<p<1/2,$ $f\in H_{p}$ and $\{q_{k}:k\geq0\}$ be a sequence of non-decreasing numbers, then there exists an absolute constant $c_{p},$ depending only on $p,$ such that the inequality%
\begin{equation*}
\overset{\infty }{\underset{k=1}{\sum }}\frac{\left\Vert T_{k}f\right\Vert
	_{p}^{p}}{k^{2-2p}}\leq c_{p}\left\Vert f\right\Vert _{H_{p}}^{p}
\end{equation*}%
holds.

Moreover, if  $f\in H_{1/2}$ and $\{q_{k}:k\geq 0\}$ be a sequence of non-increasing
numbers, satisfying the condition \eqref{fn01}, then there exists an absolute constant $c,$ such that the inequality%
\begin{equation*}
\frac{1}{\log n}\overset{n}{\underset{k=1}{\sum }}\frac{\left\Vert
	T_{k}f\right\Vert _{1/2}^{1/2}}{k}\leq c\left\Vert f\right\Vert
_{H_{1/2}}^{1/2}  \label{7nor}
\end{equation*}%
holds.

\subsection{Auxiliary lemmas}

\begin{lemma}\label{T1} 
	Let $n\in\mathbb{N}.$ Then
	
	\begin{eqnarray} \label{2b}
	Q_n&:=&\overset{n-1}{\underset{j=0}{\sum}}q_j =\overset{n-2}{\underset{j=0}{\sum}}\left(q_{j}-q_{j+1}\right) j+q_{n-1}{(n-1)}
	\end{eqnarray}
	
	\begin{equation} 	\label{2c}
	F_n=\frac{1}{Q_n}\left(\overset{n-2}{\underset{j=0}{\sum}}\left(q_j-q_{j+1}\right) jK_j+q_{n-1}(n-1)K_{n-1}\right).
	\end{equation}	
	
	\begin{equation} 	\label{2d}
	t_n=\frac{1}{Q_n}\left(\overset{n-2}{\underset{j=0}{\sum}}\left(q_j-q_{j+1}\right) j\sigma_j+q_{n-1}(n-1)\sigma_{n-1}\right).
	\end{equation}	
\end{lemma}
\textbf{Proof:}
If we invoke Abel transformation we immediately get  identities \eqref{2b}, \eqref{2c} and \eqref{2d}. The proof is complete.
\QED

\begin{lemma}\label{T2} 
	Let $n\in\mathbb{N}$ and $\{q_k:k\in\mathbb{N}\}$ be a sequence of non-increasing numbers, or non-decreasing function satisfying condition \eqref{fn01}. Then
	
	\begin{equation} 	\label{2e}
	\Vert F_n\Vert_1<c<\infty.
	\end{equation}	
\end{lemma}
\textbf{Proof:}
Let $n\in\mathbb{N}$ and $\{q_k:k\in\mathbb{N}\}$ be a sequence of non-increasing numbers. By combining \eqref{knbounded} with  \eqref{2b} and \eqref{2d} we can conclude that 

\begin{eqnarray*} 	
	\Vert T_n\Vert_1&\leq& \frac{1}{Q_n} \left(\overset{n-2}{\underset{j=0}{\sum}}\left\vert q_j-q_{j+1}\right\vert j\Vert \sigma_j\Vert_1+q_{n-1}(n-1)\Vert \sigma_{n-1}\Vert_1\right)\\
	&\leq &\frac{c}{Q_n} \left(\overset{n-2}{\underset{j=0}{\sum}}\left(q_j-q_{j+1}\right) j+q_{n-1}(n-1)\right)\leq c<\infty.
\end{eqnarray*}	

Let $n\in\mathbb{N}$ and $\{q_k:k\in\mathbb{N}\}$ be a sequence non-decreasing function satisfying condition \eqref{fn01}. Then By using again  \eqref{knbounded} with  \eqref{2b} and \eqref{2d} we find that

\begin{eqnarray*} 	
	\Vert T_n\Vert_1&\leq& \frac{1}{Q_n} \left(\overset{n-2}{\underset{j=0}{\sum}}\left\vert q_j-q_{j+1}\right\vert j\Vert \sigma_j\Vert_1+q_{n-1}(n-1)\Vert \sigma_{n-1}\Vert_1\right)\\
	&\leq &\frac{c}{Q_n} \left(\overset{n-2}{\underset{j=0}{\sum}}\left(q_{j+1}-q_j\right) j+q_{n-1}(n-1)\right)\\
	&=&\frac{c}{Q_n} \left(2q_{n-1}(n-1)-\left(\overset{n-2}{\underset{j=0}{\sum}}\left(q_{j}-q_{j+1}\right) j+q_{n-1}(n-1)\right)\right)
	\\
	&=&\frac{c}{Q_n}(2q_{n-1}(n-1)-Q_n) \leq c<\infty.
\end{eqnarray*}

The proof is complete.
\QED

\begin{lemma}\label{lemma0nnT0}
	Let $\{q_k:k\in\mathbb{N}\}$ be a sequence of non-increasing numbers and $n>M_N$. Then
	\begin{equation*}
	\left\vert\frac{1}{Q_n}\overset{n-1}{\underset{j=M_{N}}{\sum }}q_{j}D_j\left( x\right)\right\vert
	\leq\frac{c}{M_N}\left\{\sum_{j=0}^{\left\vert n\right\vert }M_j\left\vert K_{M_j}\right\vert\right\},
	\end{equation*}
	where $c$ is an absolute constant.
\end{lemma}

\textbf{Proof:}
Since sequence is non-increasing number  we get that
\begin{eqnarray*}
	&&\frac{1}{Q_n}\left(q_{M_N}+\overset{n-2}{\underset{j=M_N}{\sum }}\left\vert q_j-q_{j+1}\right\vert+q_{n-1}\right)\\ &\leq&\frac{1}{Q_n}\left(q_{M_N}+\overset{n-2}{\underset{j=M_N}{\sum}}\left(q_{j}-q_{j+1} \right)+q_{n-1}\right)\\
	&\leq &
	\frac{2q_{M_N}}{Q_n}\leq\frac{2q_{M_N}}{Q_{M_N+1}}\leq \frac{c}{M_N}.
\end{eqnarray*}

If we apply \eqref{2b} and \eqref{2c} in Lemma \ref{T1} and Abel transformation we immediately get that
\begin{eqnarray*}
	&&\left\vert\frac{1}{Q_n}\overset{n-1}{\underset{j=M_{N}}{\sum }}q_{j}D_j\left( x\right)\right\vert \\
	&=&  \frac{1}{Q_n}\left( q_{M_n}\sigma_{{M_n-1}}+\overset{n-2} {\underset{j=M_N}{\sum}}\left( q_j-q_{j+1}\right)\sigma_j+q_{n-1}\sigma_{n-1}\right)\\
	&\leq &  \frac{1}{Q_n}\left( q_{M_n}+\overset{n-2} {\underset{j=M_N}{\sum}}\left\vert q_j-q_{j+1}\right \vert+q_{n-1}\right)\sum_{i=0}^{\left\vert n\right\vert }M_i\left\vert K_{M_i}\right\vert \\
	&\leq &\frac{c}{M_N}\sum_{i=0}^{\left\vert n\right\vert }M_i\left\vert K_{M_i}\right\vert.
\end{eqnarray*}
The proof is complete.
\QED

\begin{lemma} \label{lemma5aaTin}
	Let $x\in I_N^{k,l}, \ \ k=0,\dots,N-1, \ \ l=k+1,\dots ,N$
	and $\{q_k:k\in \mathbb{N}\}$ be a sequence of non-increasing numbers. Then
	\begin{equation*}
	\int_{I_N}\left\vert\frac{1}{Q_n}\overset{n-1}{\underset{j=M_{N}}{\sum }}q_{j}D_j\left( x-t\right)\right\vert d\mu\left(t\right) \leq\frac{cM_lM_k}{M^2_N}.
	\end{equation*}
	Here $c$ is an absolute constant.
\end{lemma}

{\bf Proof}:
Let $x\in I_N^{k,l},$ for $0\leq k<l\leq N-1$ and $t\in I_N.$ First, we observe that $x-t\in $ $I_N^{k,l}.$ Next, we apply equality \ref{kn8}  and  Lemmas \ref{lemma5aaTin} to obtain that
\begin{eqnarray}  \label{88811}
&&\int_{I_N}\left\vert\frac{1}{Q_n}\overset{n-1}{\underset{j=M_{N}}{\sum }}q_{j}D_j\left( x-t\right)\right\vert d\mu\left(t\right)  \\ \notag
&\leq & \frac{c}{M_N}\underset{i=0}{\overset{\left\vert n\right\vert }{\sum}}M_i\int_{I_N}\left\vert K_{M_i}\left(x-t\right)\right\vert d\mu
\left(t\right) \\ \notag
&\leq& \frac{c}{M_N}\int_{I_N}\overset{l}{\underset{i=0}{\sum }}M_iM_kd\mu\left(t\right)  \\ \notag
&\leq &\frac{cM_kM_l}{M^2_N}
\end{eqnarray}
and the first estimate is proved.

Now, let $x\in I_N^{k,N}$. Since $x-t\in I_{N}^{k,N}$ for $t\in I_N,$ by
combining \eqref{9dn} and  \eqref{2dna} we have that
$$\left\vert D_{i}\left(x-t\right)\right\vert\leq M_k$$
and
\begin{eqnarray} \label{1111001}
&&\int_{I_N}\left\vert\frac{1}{Q_n}\overset{n-1}{\underset{j=M_{N}}{\sum }}q_{j}D_j\left( x-t\right)\right\vert d\mu\left(t\right)   \\ \notag
&\leq&\frac{c}{Q_n}\underset{i=0}{\overset{\left\vert n\right\vert } {\sum}}q_i\int_{I_N}\left\vert D_{i}\left(x-t\right)\right\vert d\mu\left(t\right) \\ \notag
&\leq&\frac{c}{Q_n}\overset{\left\vert n\right\vert-1} {\underset{i=0}{\sum }}q_i\int_{I_N}M_kd\mu\left(t\right) \\ \notag
&\leq&\frac{cM_k}{M_N}.
\end{eqnarray}

According to (\ref{88811}) and (\ref{1111001}) the proof is complete.

\begin{lemma}\label{lemma0nnT}
	Let $n>M_N$ and $\{q_k:k\in\mathbb{N}\}$ be a sequence of non-increasing numbers, satisfying condition \eqref{fn0111}. Then
	\begin{equation*}
	\left\vert\frac{1}{Q_n}\overset{n-1}{\underset{j=M_{N}}{\sum }}q_{j}D_j\left( x\right)\right\vert
	\leq\frac{c}{n}\left\{\sum_{j=0}^{\left\vert n\right\vert }M_j\left\vert K_{M_j}(x)\right\vert\right\},
	\end{equation*}
	where $c$ is an absolute constant.
\end{lemma}

\textbf{Proof:}
Since sequence is non-increasing number  we get that
\begin{eqnarray*}
	&&\frac{1}{Q_n}\left(\overset{n-2}{\underset{j=M_N}{\sum }}\left\vert q_j-q_{j+1}\right\vert+q_{n-1}\right)\\ &\leq&\frac{1}{Q_n}\left(\overset{n-2}{\underset{j=M_N}{\sum}}\left(q_{j}-q_{j+1} \right)+q_{n-1}\right)\\
	&\leq &
	\frac{q_{M_N}}{Q_n}\leq \frac{c}{Q_n}\leq \frac{c}{n}.
\end{eqnarray*}

If we apply \eqref{2b} and \eqref{2c} in Lemma \ref{T1} we immediately get that
\begin{eqnarray*}
	&&\left\vert\frac{1}{Q_n}\overset{n-1}{\underset{j=M_{N}}{\sum }}q_{j}D_j\left( x\right)\right\vert \\
	&\leq & \left( \frac{1}{Q_n}\left( \overset{n-2} {\underset{j=M_N+1}{\sum}}\left\vert q_j-q_{j+1}\right \vert+q_{n-1}\right)\right)\sum_{i=0}^{\left\vert n\right\vert }M_i\left\vert K_{M_i}(x)\right\vert \\
	&\leq &\frac{c}{n}\sum_{i=0}^{\left\vert n\right\vert }M_i\left\vert K_{M_i}(x)\right\vert.
\end{eqnarray*}
The proof is complete.
\QED

\begin{lemma} \label{lemma5a}
	Let $x\in I_N^{k,l}, \ \ k=0,\dots,N-2, \ \ l=k+1,\dots ,N-1$
	and $\{q_k:k\in \mathbb{N}\}$ be a sequence of non-increasing numbers,
	satisfying condition \eqref{fn0111}. Then
	\begin{equation*}
	\int_{I_N}\left\vert\frac{1}{Q_n}\overset{n-1}{\underset{j=M_{N}}{\sum }}q_{j}D_j\left( x-t\right)\right\vert d\mu\left(t\right) \leq\frac{cM_lM_k}{nM_N}.
	\end{equation*}
	Let $x\in I_N^{k,N},$ \ \ $k=0,\dots,N-1.$ Then
	\begin{equation*}
	\int_{I_N} \left\vert\frac{1}{Q_n}\overset{n-1}{\underset{j=M_{N}}{\sum }}q_{j}D_j\left( x-t\right)\right\vert d\mu\left(t\right) \leq\frac{cM_k}{M_N}.
	\end{equation*}
	Here $c$ is an absolute constant.
\end{lemma}

{\bf Proof}:
Let $x\in I_N^{k,l},$ for $0\leq k<l\leq N-1$ and $t\in I_N.$ First, we observe that $x-t\in $ $I_N^{k,l}.$ Next, we apply equality \ref{kn8}  and Lemma \ref{lemma0nnT}  to obtain that
\begin{eqnarray}  \label{8881}
&&\int_{I_N}\left\vert\frac{1}{Q_n}\overset{n-1}{\underset{j=M_{N}}{\sum }}q_{j}D_j\left( x-t\right)\right\vert d\mu\left(t\right)  \\ \notag
&\leq & \frac{c}{n}\underset{i=0}{\overset{\left\vert n\right\vert }{\sum}}M_i\int_{I_N}\left\vert K_{M_i}\left(x-t\right)\right\vert d\mu
\left(t\right) \\ \notag
&\leq& \frac{c}{n}\int_{I_N}\overset{l}{\underset{i=0}{\sum }}M_iM_kd\mu\left(t\right)  \\ \notag
&\leq &\frac{cM_kM_l}{nM_N}
\end{eqnarray}
and the first estimate is proved.

Now, let $x\in I_N^{k,N}$. Since $x-t\in I_{N}^{k,N}$ for $t\in I_N,$ by
combining again equality \eqref{kn8}  and Lemma \ref{lemma0nnT} we have that
\begin{eqnarray} \label{111100}
&&\int_{I_N}\left\vert\frac{1}{Q_n}\overset{n-1}{\underset{j=M_{N}}{\sum }}q_{j}D_j\left( x-t\right)\right\vert d\mu\left(t\right)   \\ \notag
&\leq&\frac{c}{n}\underset{i=0}{\overset{\left\vert n\right\vert } {\sum}}M_i\int_{I_N}\left\vert K_{M_i}\left(x-t\right)\right\vert d\mu\left(t\right) \\ \notag
&\leq&\frac{c}{n}\overset{\left\vert n\right\vert-1} {\underset{i=0}{\sum }}M_i\int_{I_N}M_kd\mu\left(t\right) \\ \notag
&\leq&\frac{cM_k}{M_N}.
\end{eqnarray}

By combining (\ref{8881}) and (\ref{111100}) we complete the proof.

\begin{lemma} \label{lemma5bT}
	Let $n\geq M_N, \ \ x\in I_N^{k,l}, \ \ k=0,\dots,N-1, \ \
	l=k+1,\dots,N$ and $\{q_k:k\in\mathbb{N}\}$ be a sequence of
	non-increasing sequence, satisfying condition  \eqref{fn0111}. Then
	\begin{equation*}
	\int_{I_N}\left\vert\frac{1}{Q_n}\overset{n-1}{\underset{j=M_{N}}{\sum }}q_{j}D_j\left( x-t\right)\right\vert d\mu\left(
	t\right)\leq\frac{cM_lM_k}{M_N^2},
	\end{equation*}%
	where $c$ is an absolute constant.
\end{lemma}

{\bf Proof}:
Since $n\geq M_N$ if we apply Lemma \ref{lemma5a} we immediately get the  proof.

\begin{lemma}\label{lemma0nnT1}
	Let $\{q_k:k\in\mathbb{N}\}$ be a sequence of non-decreasing numbers satisfying  (\ref{fn01}).
	Then
	\begin{equation*}
	\left\vert F_n\right\vert\leq\frac{c}{n}\left\{\sum_{j=0}^{\left\vert n\right\vert }M_j\left\vert K_{M_j}\right\vert \right\},
	\end{equation*}
	where $c$ is an absolute constant.
\end{lemma}

\textbf{Proof:}
Since sequence $\{q_k:k\in \mathbb{N}\}$ be non-decreasing. If we apply condition \eqref{fn01} we can conclude that
\begin{eqnarray*}
	&&\frac{1}{Q_n}\left(\overset{n-2}{\underset{j=0}{\sum}}\left\vert q_j-q_{j+1}\right\vert+q_{n-1}\right) \\ &&\leq\frac{1}{Q_n}\left(\overset{n-2}{\underset{j=0}{\sum }}\left (q_{j+1}-q_j \right)+q_{n-1}\right)\\
	&\leq& \frac{2q_{n-1}-q_0}{Q_{n}}\leq \frac{q_{n-1}}{Q_{n}}
	\leq \frac{c}{n}.
\end{eqnarray*}

If we apply \eqref{2b} and \eqref{2c} in Lemma \ref{T1} we immediately get that
\begin{eqnarray*}
	\left\vert F_n\right\vert &\leq & \left( \frac{1}{Q_n}\left( \overset{n-1}{\underset{j=1}{\sum }}\left\vert q_{j}-q_{j+1} \right\vert+q_0\right)\right)\sum_{i=0}^{\left\vert n\right\vert } M_i\left\vert K_{M_i}\right\vert \\
	&=&\left(\frac{1}{Q_n}\left(\overset{n-1}{\underset{j=1}{\sum}}\left(q_{j}-q_{j+1}\right)+q_0\right)\right) \sum_{i=0}^{\left\vert n\right\vert}M_i\left\vert K_{M_i}\right\vert \\
	&\leq & \frac{q_{n-1}}{Q_n}\sum_{i=0}^{\left\vert n\right\vert}M_i\left\vert K_{M_i}\right\vert \leq\frac{c}{n}\sum_{i=0}^{\left\vert n\right\vert }M_i\left\vert K_{M_i}\right\vert.
\end{eqnarray*}
The proof is complete.
\QED

\begin{lemma} \label{lemma5aT}
	Let $x\in I_N^{k,l}, \ \ k=0,\dots,N-2, \ \ l=k+1,\dots ,N-1$
	and $\{q_k:k\in \mathbb{N}\}$ be a sequence of non-decreasing numbers,
	satisfying condition (\ref{fn01}). Then
	\begin{equation*}
	\int_{I_N}\left\vert F_n\left(x-t\right)\right\vert d\mu\left(t\right) \leq\frac{cM_lM_k}{nM_N}.
	\end{equation*}
	Let $x\in I_N^{k,N},$ \ \ $k=0,\dots,N-1.$ Then
	\begin{equation*}
	\int_{I_N}\left\vert F_n\left(x-t\right)\right\vert d\mu\left(t\right) \leq\frac{cM_k}{M_N}.
	\end{equation*}%
	Here $c$ is an absolute constant.
\end{lemma}

{\bf Proof}:
The proof is quite analogously to Lemma \ref{lemma5a}. So we leave out the details.

\begin{lemma} \label{lemma5b}
	Let $n\geq M_N, \ \ x\in I_N^{k,l}, \ \ k=0,\dots,N-1, \ \
	l=k+1,\dots,N$ and $\{q_k:k\in\mathbb{N}\}$ be a sequence of
	non-decreasing sequence, satisfying condition (\ref{fn01}). Then
	\begin{equation*}
	\int_{I_N}\left\vert F_n\left(x-t\right)\right\vert d\mu\left(
	t\right)\leq\frac{cM_lM_k}{M_N^2},
	\end{equation*}%
	where $c$ is an absolute constant.
\end{lemma}

{\bf Proof}:
Since $n\geq M_N$ if we apply Lemma \ref{lemma5aT} we immediately get the  proof.

\subsection{On the maximal operators of $T$ means with respect to Vilenkin systems on the martingale Hardy spaces}

\ \ \ First we state our main result concerning the maximal operator of the
summation method (\ref{nor}), which we also show is in a sense sharp.

\begin{theorem}\label{theorem1T}
	a) The maximal operator $T^{\ast }$ of the summability method (\ref{nor})
	with non-increasing sequence $\{q_{k}:k\geq 0\},$ is bounded from the Hardy
	space $H_{1/2}$ to the space $weak-L_{1/2}.$
	
	The statement in a) is sharp in the following sense:
	
	\bigskip b) Let $0<p<1/2$ and\ $\{q_{k}:k\geq 0\}$ be a non-increasing sequence, satisfying the condition
	\begin{equation}
	\frac{q_{n+1}}{Q_{n+2}}\geq \frac{c}{n},\text{ \ \ }\left( c\geq 1\right) .
	\label{cond1}
	\end{equation}%
	Then there exists a martingale $f\in H_{p},$ such that
	\begin{equation*}
	\underset{n\in	\mathbb{N}}{\sup }\left\Vert T_{n}f\right\Vert _{weak-L_{p}}=\infty .
	\end{equation*}
\end{theorem}

\textbf{Proof:} a). Let the sequence $\{q_{k}:k\geq 0\}$ be non-increasing. By combining (\ref{2b}) with (\ref{2c}) and using Abel transformation we get that

\begin{eqnarray*}
	\left\vert T_{n}f\right\vert &\leq &\left\vert \frac{1}{Q_{n}}\overset{n-1}{%
		\underset{j=1}{\sum }}q_{j}S_{j}f\right\vert \\
	&\leq &\frac{1}{Q_{n}}\left( \overset{n-2}{\underset{j=1}{\sum }}\left\vert
	q_{j}-q_{j+1}\right\vert j\left\vert \sigma _{j}f\right\vert
	+q_{n-1}(n-1)\left\vert \sigma _{n}f\right\vert \right) \\
	&\leq &\frac{1}{Q_{n}}\left( \overset{n-2}{\underset{j=1}{\sum }}\left(
	q_{j}-q_{j+1}\right) j+q_{n-1}(n-1)\right) \sigma ^{\ast }f\leq \sigma ^{\ast}f
\end{eqnarray*}
so that
\begin{equation}
T^{\ast }f\leq \sigma ^{\ast }f.  \label{12aaaa}
\end{equation}

If we apply (\ref{12aaaa}), according that $\sigma^*$ is bounded from the Hardy space $H_{1/2}$ to the space $weak-L_{1/2},$ we can conclude that the maximal operators  $T^{\ast }$ of all $T$ means with non-increasing sequence $\{q_{k}:k\geq 0\},$ are bounded from the Hardy space $H_{1/2}$ to the space $weak-L_{1/2}.$ The proof of part a) of Theorem 1 is complete.

b) Let $0<p<1/2$ and $\left\{ \alpha _{k}:k\in\mathbb{N}\right\} $ be an increasing sequence of positive integers such that:\qquad
\begin{equation}
\sum_{k=0}^{\infty }1/\alpha _{k}^{p}<\infty ,  \label{2t}
\end{equation}

\begin{equation}
\lambda \sum_{\eta =0}^{k-1}\frac{M_{\alpha _{\eta }}^{1/p}}{\alpha _{\eta }}%
<\frac{M_{\alpha _{k}}^{1/p}}{\alpha _{k}},  \label{3t}
\end{equation}

\begin{equation}
\frac{32\lambda M_{\alpha _{k-1}}^{1/p}}{\alpha _{k-1}}<\frac{M_{\alpha_{k}}^{1/p-2}}{\alpha _{k}},  \label{4t}
\end{equation}
where $\lambda =\sup_{n}m_{n}.$

We note that such an increasing sequence 
$\left\{ \alpha _{k}:k\in\mathbb{N}\right\} $ 
which satisfies conditions (\ref{2t}), (\ref{3t}) and (\ref{4t}) can be constructed.

Let \qquad
\begin{equation}
f^{\left( A\right) }=\sum_{\left\{ k;\text{ }\lambda _{k}<A\right\} }\lambda
_{k}a_{k},  \label{55}
\end{equation}%
where

\begin{equation*}
\lambda _{k}=\frac{\lambda }{\alpha _{k}}  \ \ \ \  and \ \ \ \
a_{k}=\frac{M_{\alpha _{k}}^{1/p-1}}{\lambda }\left( D_{M_{\alpha
		_{k}+1}}-D_{M_{_{\alpha _{k}}}}\right) .  \label{77}
\end{equation*}

By using Lemma \ref{lemma2.1}, it is easy to show that the martingale $f \in H_{1/2}.$ Moreover, it is easy to show that

\begin{equation}
\widehat{f}(j)=\left\{
\begin{array}{ll}
\frac{M_{\alpha _{k}}^{1/p-1}}{\alpha _{k}},\, & \text{if \thinspace
	\thinspace }j\in \left\{ M_{\alpha _{k}},...,M_{\alpha _{k}+1}-1\right\} ,%
\text{ }k=0,1,2..., \\
0, & \text{if \ \thinspace \thinspace \thinspace }j\notin
\bigcup\limits_{k=1}^{\infty }\left\{ M_{\alpha _{k}},...,M_{\alpha
	_{k}+1}-1\right\} .%
\end{array}%
\right.  \label{6aa}
\end{equation}

We can write

\begin{equation}
T_{M_{\alpha _{k}}+2}f=\frac{1}{Q_{M_{\alpha _{k}}+2}}\sum_{j=0}^{M_{\alpha
		_{k}}}q_{j}S_{j}f+\frac{q_{M_{\alpha
			_k}+1}}{Q_{M_{\alpha _{k}}+2}}S_{M_{\alpha
		_{k}}+1}f:=I+II. \label{155aba}
\end{equation}

Let $M_{\alpha _{s}}\leq $ $j\leq M_{\alpha _{s}+1},$ where $s=0,...,k-1.$
Moreover,
\begin{equation*}
\left\vert D_{j}-D_{M_{_{\alpha _{s}}}}\right\vert \leq j-M_{_{\alpha
		_{s}}}\leq \lambda M_{_{\alpha _{s}}},\text{ \ \ }\left( s\in
\mathbb{N}\right)
\end{equation*}%
so that, according to (\ref{dn21}) and (\ref{6aa}), we have that%
\begin{eqnarray}
&&\left\vert S_{j}f\right\vert =\left\vert \sum_{v=0}^{M_{\alpha _{s-1}+1}-1}%
\widehat{f}(v)\psi _{v}+\sum_{v=M_{\alpha _{s}}}^{j-1}\widehat{f}(v)\psi
_{v}\right\vert  \label{8a} \\
&\leq &\left\vert \sum_{\eta =0}^{s-1}\sum_{v=M_{\alpha _{\eta
}}}^{M_{\alpha _{\eta }+1}-1}\frac{M_{\alpha _{\eta }}^{1/p-1}}{\alpha
	_{\eta }}\psi _{v}\right\vert +\frac{M_{\alpha _{s}}^{1/p-1}}{\alpha _{s}}%
\left\vert \left( D_{j}-D_{M_{_{\alpha _{s}}}}\right) \right\vert  \notag \\
&=&\left\vert \sum_{\eta =0}^{s-1}\frac{M_{\alpha _{\eta }}^{1/p-1}}{\alpha
	_{\eta }}\left( D_{M_{_{\alpha _{\eta }+1}}}-D_{M_{_{\alpha _{\eta
}}}}\right) \right\vert +\frac{M_{\alpha _{s}}^{1/p-1}}{\alpha _{s}}%
\left\vert \left( D_{j}-D_{M_{_{\alpha _{s}}}}\right) \right\vert  \notag \\
&\leq &\lambda \sum_{\eta =0}^{s-1}\frac{M_{\alpha _{\eta }}^{1/p}}{\alpha
	_{\eta }}+\frac{\lambda M_{\alpha _{s}}^{1/p}}{\alpha _{s}}\leq \frac{%
	2\lambda M_{\alpha _{s-1}}^{1/p}}{\alpha _{s-1}}+\frac{\lambda M_{\alpha
		_{s}}^{1/p}}{\alpha _{s}}\leq \frac{4\lambda M_{\alpha _{k-1}}^{1/p}}{\alpha
	_{k-1}}.  \notag
\end{eqnarray}

Let $M_{\alpha _{s-1}+1}+1\leq $ $j\leq M_{\alpha _{s}},$ where $s=1,...,k.$
Analogously to (\ref{8}) we can prove that

\begin{eqnarray*}
	&&\left\vert S_{j}f\right\vert =\left\vert \sum_{v=0}^{M_{\alpha _{s-1}+1}-1}%
	\widehat{f}(v)\psi _{v}\right\vert =\left\vert \sum_{\eta
		=0}^{s-1}\sum_{v=M_{\alpha _{\eta }}}^{M_{\alpha _{\eta }+1}-1}\frac{%
		M_{\alpha _{\eta }}^{1/p-1}}{\alpha _{\eta }}\psi _{v}\right\vert \\
	&=&\left\vert \sum_{\eta =0}^{s-1}\frac{M_{\alpha _{\eta }}^{1/p-1}}{\alpha
		_{\eta }}\left( D_{M_{_{\alpha _{\eta }+1}}}-D_{M_{_{\alpha _{\eta
	}}}}\right) \right\vert \leq \frac{2\lambda M_{\alpha _{s-1}}^{1/p}}{\alpha
		_{s-1}}\leq \frac{4\lambda M_{\alpha _{k-1}}^{1/p}}{\alpha _{k-1}}.
\end{eqnarray*}
Hence,
\begin{equation}
\left\vert I\right\vert \leq \frac{1}{Q_{M_{\alpha _{k}}+2}}
\sum_{j=0}^{M_{\alpha _{k}}}q_{j}\left\vert S_{j}f\right\vert \leq \frac{4\lambda M_{\alpha _{k-1}}^{1/p}}{\alpha _{k-1}}\frac{1}{Q_{M_{\alpha
			_{k}}+2}}\sum_{j=0}^{M_{\alpha _{k}}}q_{j}\leq \frac{4\lambda M_{\alpha
		_{k-1}}^{1/p}}{\alpha _{k-1}}.  \label{10}
\end{equation}
If we now apply (\ref{6aa}) and (\ref{8a}) we get that
\begin{eqnarray}
\left\vert II\right\vert &=&\frac{q_{M_{\alpha _{k}}+1}}{Q_{M_{\alpha_{k}}+2}}\left\vert \frac{%
	M_{\alpha _{k}}^{1/p-1}}{\alpha _{k}}\psi _{M_{\alpha _{k}}}+S_{M_{\alpha
		_{k}}}f\right\vert  \label{100} \\
&=&\frac{q_{M_{\alpha _{k}}+1}}{Q_{M_{\alpha_{k}}+2}}\left\vert \frac{M_{\alpha _{k}}^{1/p-1}}{%
	\alpha _{k}}\psi _{M_{\alpha _{k}}}+S_{M_{\alpha _{k-1}+1}}f\right\vert
\notag \\
&\geq &\frac{q_{M_{\alpha _{k}}+1}}{Q_{M_{\alpha_{k}}+2}}\left( \left\vert \frac{M_{\alpha
		_{k}}^{1/p-1}}{\alpha _{k}}\psi _{M_{\alpha _{k}}}\right\vert -\left\vert
S_{M_{\alpha _{k-1}+2}}f\right\vert \right)  \notag \\
&\geq &\frac{q_{M_{\alpha _{k}}+1}}{Q_{M_{\alpha _{k}}+2}}\left( \frac{M_{\alpha
		_{k}}^{1/p-1}}{\alpha _{k}}-\frac{4\lambda M_{\alpha _{k-1}}^{1/p}}{\alpha
	_{k-1}}\right) \notag\\ 
&\geq& \frac{q_{M_{\alpha _{k}}+1}}{Q_{M_{\alpha _{k}}+2}}\frac{M_{\alpha
		_{k}}^{1/p-1}}{4\alpha _{k}}.  \notag
\end{eqnarray}

Without lost the generality we may assume that $c=1$ in (\ref{cond1}). By
combining (\ref{10}) and (\ref{100}) we get
\begin{eqnarray} \label{777}
\left\vert T_{M_{\alpha _{k}}+2}f\right\vert &\geq &\left\vert II\right\vert
-\left\vert I\right\vert \geq \frac{q_{M_{\alpha _{k}}+1}}{Q_{M_{\alpha _{k}}+2}}\frac{	M_{\alpha _{k}}^{1/p-1}}{4\alpha _{k}}-\frac{4\lambda M_{\alpha _{k-1}}^{1/p}}{\alpha _{k-1}} \\
&\geq &\frac{M_{\alpha _{k}}^{1/p-2}}{4\alpha _{k}}-\frac{4\lambda M_{\alpha
		_{k-1}}^{1/p}}{\alpha _{k-1}}\geq \frac{M_{\alpha _{k}}^{1/p-2}}{16\alpha _{k}}. \notag
\end{eqnarray}
On the other hand,
\begin{equation} \label{88}
\mu \left\{ x\in G_{m}:\left\vert T_{M_{\alpha _{k}}+2}f\left( x\right)
\right\vert \geq \frac{M_{\alpha _{k}}^{1/p-2}}{16\alpha _{k}}\right\} =\mu
\left( G_{m}\right) =1.  
\end{equation}
Let $0<p<1/2.$ Then
\begin{eqnarray}\label{99}
&&\frac{M_{\alpha _{k}}^{1/p-2}}{16\alpha _{k}}\cdot \left(\mu \left\{ x\in
G_{m}:\left\vert T_{M_{\alpha _{k}}+2}f\left( x\right) \right\vert \geq
\frac{M_{\alpha _{k}}^{1/p-2}}{16\alpha _{k}}\right\} \right)^{1/p} \\
&=&\frac{M_{\alpha _{k}}^{1/p-2}}{16\alpha _{k}}\rightarrow \infty ,\text{ \	as }k\rightarrow \infty .  \notag
\end{eqnarray}%
\textbf{\ }
The proof is complete.

\QED

A number of special cases of our results are of particular interest and give both well-known and new information. We just give the following examples of such $T$ means with non-increasing sequence $\{q_{k}:k\geq 0\}:$

\begin{corollary}\label{cor0}
	\textbf{\ }The maximal operators $U^{\alpha, *}$, $V^{\alpha, *}$ and $R^{*}$ are  bounded from the Hardy space $H_{1/2}$ to the space $%
	weak-L_{1/2}$ but are not bounded from $H_{p}$ to the space $weak-L_{p},$
	when $0<p<1/2.$
\end{corollary}

\textbf{Proof:} 
Since $R_{n},	U_{n}^{\alpha} \text{ and } V_{n}^{\alpha}$ are  the $T$ means with non-increasing sequence $\{q_{k}:k\geq 0\},$ then the proof of this corollary is direct consequence of Theorem \ref{theorem1T}.

\begin{corollary}\label{cor2}
	Let $f\in L_{1}$ and $T_{n}$ be the $T$ means with non-increasing sequence $\{q_{k}:k\geq 0\}$. Then
	$\	T_{n}f\rightarrow f,\text{ \ \ a.e., \ \ as \ }n\rightarrow \infty .$
\end{corollary}
\textbf{Proof:} 
According to part a) of Theorem \ref{theorem1T} and Lemma \ref{lemma2.3} we also have weak $(1,1)$ type inequality and by  well-known density argument due to Marcinkiewicz and Zygmund \cite{13z} (see Lemma \ref{lemmaae}) we have $T_{n}f\rightarrow f,$ a.e., for all $f\in L_1.$ Which follows proof of Corollary \ref{cor2}.

\begin{corollary}\label{cor1}
	Let $f\in L_{1}$. Then%
	\begin{eqnarray*}
		R_{n}f &\rightarrow &f,\text{ \ \ \ a.e., \ \ \ \ as \ }n\rightarrow \infty, \\
		U_{n} ^{\alpha}f &\rightarrow &f,\text{ \ \ \ a.e., \ \ \ as \ }n\rightarrow
		\infty ,\text{\ \ \ }\\
		V_{n} ^ {\alpha}f &\rightarrow &f,\text{ \ \ \ a.e., \ \ \ as \ }n\rightarrow
		\infty ,\text{\ \ \ }
	\end{eqnarray*}
\end{corollary}
\textbf{Proof:} 
Since $R_{n},	U_{n}^{\alpha} \text{ and } V_{n}^{\alpha}$ are  the $T$ means with non-increasing sequence $\{q_{k}:k\geq 0\},$ then the proof of this corollary is direct consequence of Corollary \ref{cor2}.

Our next main result reads:
\begin{theorem} \label{theorem2}
	a) The maximal operator $T^{\ast }$ of \ the summability method (\ref{nor})
	with non-decreasing sequence $\{q_{k}:k\geq 0\}$ satisfying the condition 
	\begin{equation} \label{Tmeanscond}
	\frac{q_{n-1}}{Q_n}= O\left(\frac{1}{n}\right)
	\end{equation}
	is bounded from the Hardy space $H_{1/2}$ to the space $weak-L_{1/2}.$
	
	b) Let $0<p<1/2.$ For any non-decreasing sequence $\{q_{k}:k\geq 0\},$ there exists a martingale $f\in H_{p},$ such that
	\begin{equation*}
	\underset{n\in
		\mathbb{N}
	}{\sup }\left\Vert T_{n}f\right\Vert _{weak-L_{p}}=\infty .
	\end{equation*}

\end{theorem}

\textbf{Proof:} 
Let the sequence $\{q_{k}:k\geq 0\}$ be
non-decreasing. By combining (\ref{2b}) with (\ref{2c}) in Lemma \ref{T1} and using Abel transformation we get that
\begin{eqnarray*}
	\left\vert T_{n}f\right\vert &\leq &\left\vert \frac{1}{Q_{n}}\overset{n-1}{%
		\underset{j=1}{\sum }}q_{j}S_{j}f\right\vert \\
	&\leq &\frac{1}{Q_{n}}\left( \overset{n-2}{\underset{j=1}{\sum }}\left\vert
	q_{j}-q_{j+1}\right\vert j\left\vert \sigma _{j}f\right\vert
	+q_{n-1}(n-1)\left\vert \sigma _{n}f\right\vert \right) \\
	&\leq &\frac{1}{Q_{n}}\left( \overset{n-2}{\underset{j=1}{\sum }}-\left(
	q_{j}-q_{j+1}\right) j-q_{n-1}(n-1)+2q_{n-1}(n-1)\right) \sigma ^{\ast }f \\
	&\leq &\frac{1}{Q_{n}}\left( 2q_{n-1}(n-1)-Q_{n}\right) \sigma ^{\ast }f\leq c\sigma ^{\ast	}f
\end{eqnarray*}
so that
\begin{equation}
T^{\ast }f\leq c \sigma ^{\ast }f.  \label{12aaaaa}
\end{equation}

If we apply (\ref{12aaaaa}), according that $\sigma^*$ is bounded from the Hardy space $H_{1/2}$ to the space $weak-L_{1/2},$ we can conclude that the maximal operators  $T^{\ast }$ of all $T$ means	with non-decreasing sequence $\{q_{k}:k\geq 0\}$ satisfying the condition \ref{Tmeanscond}
are bounded from the Hardy space $H_{1/2}$ to the space $weak-L_{1/2}.$ 
The proof of part a) is complete.

To prove part b) of Theorem 2 we use the martingale, defined by (\ref{55}) where $\alpha_k$ satisfy conditions (\ref{2t}), (\ref{3t}) and (\ref{4t}). It is easy to show that for every non-increasing sequence $\{q_{k}:k\geq 0\}$ it automatically holds that
\begin{equation*}
q_{M_{\alpha _{k}+1}}/Q_{M_{\alpha _{k}+2}}\geq c/M_{\alpha _{k}}.
\end{equation*}

According to (\ref{777}), (\ref{88}) and \eqref{99} we can conclude that
\begin{equation*}
\left\vert T_{M_{\alpha _{k}}+2}f\right\vert \geq \left\vert II\right\vert
-\left\vert I\right\vert \geq \frac{M_{\alpha _{k}}^{1/p-2}}{8\alpha _{k}}.
\end{equation*}

Analogously to\ (\ref{88})  we then get that%
\begin{equation*}
\sup_{k\in \mathbb{N}}\left\Vert T_{M_{\alpha _{k}}+2}f\right\Vert _{weak-L_{p}}=\infty .
\end{equation*}

The proof is complete.

\QED
A number of special cases of our results are of particular interest and give both well-known and new information. We just give the following examples of such  $T$ means with non-decreasing sequence $\{q_{k}:k\geq 0\}:$

\begin{corollary} \label{cor4}
	\textbf{\ }The maximal operator $B^{\alpha,\beta, *}$ is  bounded from the Hardy space $H_{1/2}$ to the space $
	weak-L_{1/2}$ but is not bounded from $H_{p}$ to the space $weak-L_{p},$
	when $0<p<1/2.$
\end{corollary}

\textbf{Proof:} Since $B^{\alpha,\beta, *}$ are  the $T$ means with non-decreasing sequence $\{q_{k}:k\geq 0\},$ then the proof of this corollary is direct consequence of Theorem \ref{theorem2}.

\begin{corollary}\label{cor5}
	Let $f\in L_{1}$ and $T_{n}$ be the $T$ means with \textit{non-decreasing sequence }$\{q_{k}:k\geq 0\}$ \textit{and} satisfying condition (\ref{Tmeanscond}). Then
	\begin{equation*}
	T_{n}f\rightarrow f,\text{ \ \ a.e., \ \ as \ }n\rightarrow \infty .
	\end{equation*}
\end{corollary}
\textbf{Proof of Corollary \ref{cor5}.} According to Theorem \ref{theorem2} and Lemma and  \ref{lemma2.3} we can conclude that $T^{\ast }$ has weak type-(1,1) and by well-known density argument due to Marcinkiewicz and Zygmund \cite{13z} (see Lemma \ref{lemmaae}) we also have $T_{n}f\rightarrow f,$ a.e.. Which follows proof of Corollary \ref{cor5}.

\begin{corollary}\label{cor6}
	Let $f\in L_{1}$. Then $ \ B_{n}^{\alpha,\beta}f \rightarrow f,\text{ \ \ \ a.e., \ \ \ as \ }n\rightarrow 	\infty.$
\end{corollary}

\textbf{Proof:} Since $B^{\alpha,\beta, *}$ are  the $T$ means with non-decreasing sequence $\{q_{k}:k\geq 0\},$ then the proof of this corollary is direct consequence of Corollary  \ref{cor5}.

\begin{theorem}
	\label{theorem3fejermax2222}Let $0<p\leq 1/2,$ $f\in H_{p}$ and $\{q_{k}:k\geq
	0\}$ be a sequence of non-increasing numbers. Then the maximal operator 
	\begin{equation*}
	\widetilde{T}_{p}^{\ast }f:=\sup_{n\in \mathbb{N}_{+}}\frac{\left\vert
		T_{n}f\right\vert }{\left( n+1\right) ^{1/p-2}\log ^{2\left[ 1/2+p\right]
		}\left( n+1\right) }
	\end{equation*}%
	is bounded from the Hardy space $H_{p}$ to the space $L_{p}.$
\end{theorem}

\textbf{Proof:}
Let the sequence $\{q_{k}:k\geq 0\}$ be non-increasing. By combining (\ref{2b}) and (\ref{2d})  we get that

\begin{eqnarray*}
	\widetilde{T}_{p}^{\ast }f&:=&\frac{\left\vert T_{n}f\right\vert}{{\left( n+1\right) ^{1/p-2}\log ^{2\left[ 1/2+p\right] }\left( n+1\right) }} \\
	&\leq &\frac{1}{{\left( n+1\right) ^{1/p-2}\log ^{2\left[ 1/2+p%
				\right] }\left( n+1\right) }}\left\vert \frac{1}{Q_{n}}\overset{n-1}{
		\underset{j=1}{\sum }}q_{j}S_{j}f\right\vert \\
	&\leq &\frac{1}{\left( n+1\right) ^{1/p-2}\log ^{2\left[ 1/2+p%
			\right] }\left( n+1\right) }\frac{1}{Q_{n}}\left( \overset{n-2}{\underset{j=1}{\sum }}\left\vert
	q_{j}-q_{j+1}\right\vert j\left\vert \sigma _{j}f\right\vert
	+q_{n-1}(n-1)\left\vert \sigma _{n}f\right\vert \right) \\
	&\leq &\frac{1}{Q_{n}}\left( \overset{n-2}{\underset{j=1}{\sum }} \frac{\left\vert
		q_{j}-q_{j+1}\right\vert j\left\vert \sigma _{j}f\right\vert}{{\left( j+1\right)^{1/p-2}\log ^{2\left[ 1/2+p\right] }\left( j+1\right) }}
	+\frac{q_{n-1}(n-1)\left\vert \sigma _{n}f\right\vert}{\left( n+1\right)^{1/p-2}\log ^{2\left[ 1/2+p\right] }\left( n+1\right) } \right) \\
	&\leq &\frac{1}{Q_{n}}\left( \overset{n-2}{\underset{j=1}{\sum }}\left(
	q_{j}-q_{j+1}\right) j+q_{n-1}(n-1)\right) \sup_{n\in \mathbb{N}_{+}}\frac{\left\vert
		\sigma _{n}f\right\vert }{\left( n+1\right) ^{1/p-2}\log ^{2\left[ 1/2+p%
			\right] }\left( n+1\right) } \\
	&\leq& \sup_{n\in \mathbb{N}_{+}}
	\frac{\left\vert
		\sigma _{n}f\right\vert }{\left( n+1\right) ^{1/p-2}\log ^{2\left[ 1/2+p
			\right] }\left( n+1\right) }:=\widetilde{\sigma }_{p}^{\ast }f,
\end{eqnarray*}
so that
\begin{equation}
\widetilde{T}_{p}^{\ast }f\leq \widetilde{\sigma }_{p}^{\ast }f.  \label{12aaaaaa}
\end{equation}

If we apply (\ref{12aaaaaa}), according (see Tephnadze \cite{tep2,tep3})  that  $\widetilde{\sigma }_{p}^{\ast }f$ is bounded from the Hardy space $H_{p}$ to the space $L_{p}$ for $p\leq 1/2$ we can conclude that the maximal operators  $\widetilde{T}_{p}^{\ast }$ of  $T$ means with non-increasing sequence $\{q_{k}:k\geq 0\},$ are bounded  from the Hardy space $H_{p}$ to the space $L_{p}$.

The  proof of theorem  is complete.
\QED

\begin{corollary}\label{corolary1}
	Let $0<p\leq 1/2$ and $f\in H_{p}$. Then the maximal operator 
	
	\begin{equation*}
	\widetilde{R}_{p}^{\ast }f:=\sup_{n\in \mathbb{N}_{+}}\frac{\left\vert
		R_{n}f\right\vert }{\left( n+1\right) ^{1/p-2}\log ^{2\left[ 1/2+p\right]
		}\left( n+1\right) }
	\end{equation*}%
	is bounded from the Hardy space $H_{p}$ to the space $L_{p}.$
\end{corollary}

\begin{corollary}\label{corolary2}
	Let $0<p\leq 1/2$ and $f\in H_{p}$. Then the maximal operator 
	
	\begin{equation*}
	\widetilde{U}^{\alpha,\ast}_{p}f:=\sup_{n\in \mathbb{N}_{+}}\frac{\left\vert
		U^{\alpha}_{n}f\right\vert }{\left( n+1\right) ^{1/p-2}\log ^{2\left[ 1/2+p\right]
		}\left( n+1\right) }
	\end{equation*}%
	is bounded from the Hardy space $H_{p}$ to the space $L_{p}.$
\end{corollary}

\begin{corollary}\label{corolary3}
	Let $0<p\leq 1/2$ and $f\in H_{p}$. Then the maximal operator 
	
	\begin{equation*}
	\widetilde{V}^{\alpha,\ast}_{p}f:=\sup_{n\in \mathbb{N}_{+}}\frac{\left\vert
		V^{\alpha}_{n}f\right\vert }{\left( n+1\right) ^{1/p-2}\log ^{2\left[ 1/2+p\right]
		}\left( n+1\right) }
	\end{equation*}%
	is bounded from the Hardy space $H_{p}$ to the space $L_{p}.$
\end{corollary}

\begin{theorem}\label{theorem3fejermax22221}
	Let $0<p\leq 1/2,$ $f\in H_{p}$ and $\{q_{k}:k\geq
	0\}$ be a sequence of  non-increasing
	numbers, satisfying the condition 
	
	\begin{equation}
	\frac{q_{n-1}}{Q_{n}}=O\left( \frac{1}{n}\right) ,\text{ \ \ as \ \ }\
	n\rightarrow \infty .  \label{fn01T}
	\end{equation}
	Then the maximal operator 
	
	\begin{equation*} 
	\widetilde{T}_{p}^{\ast }f:=\sup_{n\in \mathbb{N}_{+}}\frac{\left\vert
		T_{n}f\right\vert }{\left( n+1\right) ^{1/p-2}\log ^{2\left[ 1/2+p\right]
		}\left( n+1\right) }
	\end{equation*}%
	is bounded from the Hardy space $H_{p}$ to the space $L_{p}.$
\end{theorem}

\textbf{Proof:}
Let the sequence $\{q_{k}:k\geq 0\}$ be non-decreasing satisfying the condition \eqref{fn01T}. By combining (\ref{2b}) and (\ref{2d}) we get that

\begin{eqnarray*}
	&&\frac{\left\vert T_{n}f\right\vert}{{\left( n+1\right) ^{1/p-2}\log ^{2\left[ 1/2+p\right] }\left( n+1\right) }} \\
	&\leq &\frac{1}{{\left( n+1\right) ^{1/p-2}\log ^{2\left[ 1/2+p%
				\right] }\left( n+1\right) }}\left\vert \frac{1}{Q_{n}}\overset{n-1}{
		\underset{j=1}{\sum }}q_{j}S_{j}f\right\vert \\
	&\leq &\frac{1}{\left( n+1\right) ^{1/p-2}\log ^{2\left[ 1/2+p%
			\right] }\left( n+1\right) }\frac{1}{Q_{n}}\left( \overset{n-2}{\underset{j=1}{\sum }}\left\vert
	q_{j}-q_{j+1}\right\vert j\left\vert \sigma _{j}f\right\vert
	+q_{n-1}(n-1)\left\vert \sigma _{n}f\right\vert \right) \\
	&\leq &\frac{1}{Q_{n}}\left( \overset{n-2}{\underset{j=1}{\sum }} \frac{\left\vert
		q_{j}-q_{j+1}\right\vert j\left\vert \sigma _{j}f\right\vert}{{\left( j+1\right)^{1/p-2}\log ^{2\left[ 1/2+p\right] }\left( j+1\right) }}
	+\frac{q_{n-1}(n-1)\left\vert \sigma _{n}f\right\vert}{\left( n+1\right)^{1/p-2}\log ^{2\left[ 1/2+p\right] }\left( n+1\right) } \right) \\
	&\leq &\frac{1}{Q_{n}}\left( \overset{n-2}{\underset{j=1}{\sum }}\left(
	q_{j+1}-q_{j}\right) j+q_{n-1}(n-1)\right) \sup_{n\in \mathbb{N}_{+}}\frac{\left\vert
		\sigma _{n}f\right\vert }{\left( n+1\right) ^{1/p-2}\log ^{2\left[ 1/2+p%
			\right] }\left( n+1\right) } \\
	&\leq& \frac{2q_{n-1}(n-1)-Q_n}{Q_{n}}\sup_{n\in \mathbb{N}_{+}}
	\frac{\left\vert
		\sigma _{n}f\right\vert }{\left( n+1\right) ^{1/p-2}\log ^{2\left[ 1/2+p
			\right] }\left( n+1\right) }\\
	&\leq& \sup_{n\in \mathbb{N}_{+}}
	\frac{\left\vert
		\sigma _{n}f\right\vert }{\left( n+1\right) ^{1/p-2}\log ^{2\left[ 1/2+p
			\right] }\left( n+1\right) }
\end{eqnarray*}
so that

\begin{equation}
\sup_{n\in \mathbb{N}_{+}}\frac{\left\vert T_{n}f\right\vert}{{\left( n+1\right) ^{1/p-2}\log ^{2\left[ 1/2+p\right] }\left( n+1\right) }}\leq \sup_{n\in \mathbb{N}_{+}}
\frac{\left\vert\sigma _{n}f\right\vert }{\left( n+1\right) ^{1/p-2}\log ^{2\left[ 1/2+p\right] }\left( n+1\right) }.  \label{12aaaaa1}
\end{equation}

If we apply (\ref{12aaaaa1}),   according (see Tephnadze \cite{tep2,tep3})  that  $\widetilde{\sigma }_{p}^{\ast }f$ is bounded from the Hardy space $H_{p}$ to the space $L_{p}$ for $p\leq 1/2$ we can conclude that the maximal operators  $\widetilde{T}_{p}^{\ast }$ of  $T$ means with non-decreasing sequence $\{q_{k}:k\geq 0\},$ are bounded  from the Hardy space $H_{p}$ to the space $L_{p}$.

\QED

\begin{corollary}
	Let $0<p\leq 1/2,$ $f\in H_{p}$ and $\{q_{k}:k\geq 0\}$ be a sequence of non-decreasing numbers, such that
	\begin{equation} \label{condT1}
	\sup_{n\in \mathbb{N}}q_{n}<c<\infty .
	\end{equation}%
	Then all such $T$ means are bounded from the Hardy space $H_{p}$ to the space $L_{p}.$
\end{corollary}
\textbf{Proof:}
By using \eqref{condT1} we get
\begin{equation*}
\frac{q_{n-1}}{Q_{n}}\leq \frac{c}{Q_{n}}\leq \frac{c}{q_{0}n}=\frac{c_{1}}{n%
}=O\left( \frac{1}{n}\right),\text{ as }n\rightarrow 0,
\end{equation*}%
It follows that condition (\ref{fn01T}) is satisfied and for such $T$ means is bounded from the Hardy space $H_{p}$ to the space $L_{p}.$

\QED

\begin{remark} Since (see Tephnadze \cite{tep2,tep3}) the maximal operator
	
	\begin{equation*}
	\widetilde{\sigma }_{p}^{\ast }f:=\sup_{n\in \mathbb{N}_{+}}\frac{\left\vert
		\sigma _{n}f\right\vert }{\left( n+1\right) ^{1/p-2}\log ^{2\left[ 1/2+p\right] }\left( n+1\right) }
	\end{equation*}%
	is bounded from the martingale Hardy space $H_{p}$ to the space $L_{p}$ and the rate of denominator $\left( n+1\right) ^{1/p-2}\log ^{2\left[ 1/2+p	\right] }$ is in the sense sharp and Fejer means is example of $T$ means of as non-increasing as non-decreasing sequences we obtain that this weights are also sharp in Theorems \ref{theorem3fejermax2222} and \ref{theorem3fejermax22221}. 
\end{remark}

\subsection{Strong convergence  of $T$ means with respect to Vilenkin systems on the martingale Hardy spaces}

\begin{theorem}
	\label{theorem2fejerstrong}a) Let $0<p<1/2,$ $f\in H_{p}$ and $\{q_{k}:k\geq
	0\}$ be a sequence of non-increasing numbers. Then there exists an absolute
	constant $c_{p},$ depending only on $p,$ such that the inequality
	
	\begin{equation*}
	\overset{\infty }{\underset{k=1}{\sum }}\frac{\left\Vert T_{k}f\right\Vert
		_{p}^{p}}{k^{2-2p}}\leq c_{p}\left\Vert f\right\Vert _{H_{p}}^{p}
	\end{equation*}%
	holds.
	
	b)Let $f\in H_{1/2}$ and $\{q_{k}:k\geq 0\}$ be a sequence of non-increasing
	numbers, satisfying the condition 
	
	\begin{equation}\label{fn0111}
	\frac{1}{Q_{n}}=O\left(\frac{1}{n}\right),\text{ \ \ as \ \ }\
	n\rightarrow \infty .  
	\end{equation}
	
	Then there exists an absolute constant $c,$ such that the inequality
	
	\begin{equation}
	\frac{1}{\log n}\overset{n}{\underset{k=1}{\sum }}\frac{\left\Vert
		T_{k}f\right\Vert _{1/2}^{1/2}}{k}\leq c\left\Vert f\right\Vert
	_{H_{1/2}}^{1/2}  \label{7nor7}
	\end{equation}%
	holds.
\end{theorem}

\textbf{Proof:} Let the sequence $\{q_{k}:k\geq 0\}$ be non-increasing. By Lemma \ref{lemma2.1} (see also \eqref{lemma2.2}), the proof of part a) will be complete, if we show that

\begin{equation}
\frac{1}{\log ^{\left[ 1/2+p\right] }n}\overset{n}{\underset{m=1}{\sum }}%
\frac{\left\Vert T_{m}a\right\Vert _{H_{p}}^{p}}{m^{2-2p}}\leq c_{p},
\label{14cc}
\end{equation}%
for every $p$-atom $a,$ with support$\ I$, $\mu \left( I\right) =M_{N}^{-1}.$
We may assume that $I=I_{N}.$ It is easy to see that $S_{n}\left( a\right)
=T_{n}\left( a\right) =0,$ when $n\leq M_{N}$. Therefore, we can suppose
that $n>M_{N}$.

Let $x\in I_{N}.$ Since $T_{n}$ is bounded from $L_{\infty }$ to $L_{\infty
} $ (boundedness follows from Lemma \ref{T2}) and $\left\Vert a\right\Vert
_{\infty }\leq M_{N}^{1/p}$ we obtain that 

\begin{equation*}
\int_{I_{N}}\left\vert T_{m}a\right\vert ^{p}d\mu \leq \frac{\left\Vert
	a\right\Vert _{\infty }^{p}}{M_{N}}\leq c<\infty ,\text{ \ }0<p< 1/2.
\end{equation*}%
Hence, 
\begin{eqnarray} \label{14b14}
&&\frac{1}{\log ^{\left[ 1/2+p\right] }n}\overset{n}{\underset{m=1}{\sum }}%
\frac{\int_{I_{N}}\left\vert T_{m}a\right\vert ^{p}d\mu }{m^{2-2p}} \\ \notag
&\leq &
\frac{1}{\log ^{\left[ 1/2+p\right] }n}\overset{n}{\underset{k=1}{\sum }}%
\frac{1}{m^{2-2p}}\leq c<\infty .  
\end{eqnarray}

It is easy to see that 
\begin{eqnarray}\label{14a14}
\left\vert T_{m}a\left( x\right) \right\vert &=&\left\vert\int_{I_{N}} a\left(t\right) F_{n}\left( x-t\right)  d\mu \left( t\right)\right\vert\\ \notag
&=&\left\vert\int_{I_{N}} a\left( t\right) \frac{1}{Q_{n}}\overset{n}{\underset%
	{j=M_{N}}{\sum }}q_{j}D_{j}\left( x-t\right)  d\mu \left(
t\right) \right\vert  \\ \notag
&\leq& \left\Vert a\right\Vert _{\infty }\int_{I_{N}}\left\vert \frac{1}{Q_{n}}%
\overset{n}{\underset{j=M_{N}}{\sum }}q_{j}D_{j}\left( x-t\right)
\right\vert d\mu \left( t\right) \\ \notag
&\leq& M_{N}^{1/p}\int_{I_{N}}\left\vert 
\frac{1}{Q_{n}}\overset{n}{\underset{j=M_{N}}{\sum }}q_{j}D_{j}\left(
x-t\right) \right\vert d\mu \left( t\right) 
\end{eqnarray}

Let $T_{n}$ be $T$ means, with non-decreasing coefficients $%
\{q_{k}:k\geq 0\}$ and $x\in I_{N}^{k,l},\,0\leq k<l\leq N.$ Then, in the view of Lemma \ref{lemma5aa} we get that 

\begin{equation} \label{12q1}
\left\vert T_{m}a\left( x\right) \right\vert \leq cM_{l}M_{k}M_{N}^{1/p-2},%
\text{ for }0<p\leq 1/2.  
\end{equation}

Let $0<p<1/2.$ \ By using (\ref{1.1}), (\ref{14a14}), (\ref{12q1}) we find that

\begin{eqnarray}\label{7aaaa}
\int_{\overline{I_{N}}}\left\vert T_{m}a\right\vert ^{p}d\mu &=&\overset{N-2}{%
	\underset{k=0}{\sum }}\overset{N-1}{\underset{l=k+1}{\sum }}%
\sum\limits_{x_{j}=0,\text{ }j\in \{l+1,\dots
	,N-1\}}^{m_{j-1}}\int_{I_{N}^{k,l}}\left\vert T_{m}a\right\vert ^{p}d\mu  \\ \notag
&+&\overset{N-1}{\underset{k=0}{\sum }}\int_{I_{N}^{k,N}}\left\vert
T_{m}a\right\vert ^{p}d\mu \\ \notag
&\leq & c\overset{N-2}{\underset{k=0}{\sum }}\overset{N-1}{\underset{l=k+1}{%
		\sum }}\frac{m_{l+1}\dotsm m_{N-1}}{M_{N}}\left( M_{l}M_{k}\right)
^{p}M_{N}^{1-2p} \\ \notag
&+&\overset{N-1}{\underset{k=0}{\sum }}\frac{1}{M_{N}}%
M_{k}^{p}M_{N}^{1-p}  \\ \notag
&\leq & cM_{N}^{1-2p}\overset{N-2}{\underset{k=0}{\sum }}\overset{N-1}{\underset%
	{l=k+1}{\sum }}\frac{\left( M_{l}M_{k}\right) ^{p}}{M_{l}} \\ \notag
&+&\overset{N-1}{%
	\underset{k=0}{\sum }}\frac{M_{k}^{p}}{M_{N}^{p}}\leq cM_{N}^{1-2p}.
\end{eqnarray}

Moreover, according to (\ref{7aaaa}), we get that

\begin{equation*}
\overset{\infty }{\underset{m=M_{N}+1}{\sum }}\frac{\int_{\overline{I_{N}}%
	}\left\vert T_{m}a\right\vert ^{p}d\mu }{m^{2-2p}}\leq \overset{\infty }{%
	\underset{m=M_{N}+1}{\sum }}\frac{cM_{N}^{1-2p}}{m^{2-2p}}<c<\infty ,\text{
	\ }\left( 0<p<1/2\right) .
\end{equation*}%
The proof of part a) is complete.

Let $p=1/2$ and $T_{n}$ be $T$ means, with non-increasing coefficients $%
\{q_{k}:k\geq 0\}$, satisfying condition (\ref{fn0111}). By Lemma \ref{lemma2.1}, the proof of part b)  will be complete, if we show that%

\begin{equation}
\frac{1}{\log n}\overset{n}{\underset{m=1}{\sum }}%
\frac{\left\Vert T_{m}a\right\Vert _{H_{1/2}}^{1/2}}{m}\leq c_{p},
\label{14c}
\end{equation}%
for every $p$-atom $a,$ with support$\ I$, $\mu \left( I\right) =M_{N}^{-1}.$
We may assume that $I=I_{N}.$ It is easy to see that $S_{n}\left( a\right)
=T_{n}\left( a\right) =0,$ when $n\leq M_{N}$. Therefore, we can suppose
that $n>M_{N}$.

Let $x\in I_{N}.$ Since $T_{n}$ is bounded from $L_{\infty }$ to $L_{\infty
} $ (boundedness follows from Lemma \ref{T2}) and $\left\Vert a\right\Vert
_{\infty }\leq M_{N}^{2}$ we obtain that 

\begin{equation*}
\int_{I_{N}}\left\vert T_{m}a\right\vert ^{1/2}d\mu \leq \frac{\left\Vert
	a\right\Vert _{\infty }^{1/2}}{M_{N}}\leq c<\infty ,\text{ \ }0<p\leq 1/2.
\end{equation*}%
Hence, 

\begin{equation}\label{14b}
\frac{1}{\log n}\overset{n}{\underset{m=1}{\sum }}%
\frac{\int_{I_{N}}\left\vert T_{m}a\right\vert ^{1/2}d\mu }{m}\leq 
\frac{1}{\log n}\overset{n}{\underset{k=1}{\sum }}%
\frac{1}{m}\leq c<\infty .  
\end{equation}

It is easy to see that

\begin{eqnarray} \label{14aaa} 
\left\vert T_{m}a\left( x\right) \right\vert &=&\left\vert\int_{I_{N}} a\left( t\right) \frac{1}{Q_{n}}\overset{n}{\underset%
	{j=M_{N}}{\sum }}q_{j}D_{j}\left( x-t\right)  d\mu \left(
t\right) \right\vert  \\ \notag
&\leq& \left\Vert a\right\Vert _{\infty }\int_{I_{N}}\left\vert F_{m}\left(
x-t\right) \right\vert d\mu \left( t\right) \\ \notag
&\leq& M_{N}^{2}\int_{I_{N}}\left\vert F_{m}\left( x-t\right) \right\vert d\mu
\left( t\right) .
\end{eqnarray}

Let $x\in I_{N}^{k,l},\,0\leq k<l<N.$ Then, in the view of Lemma \ref{lemma5a}
we get that 
\begin{equation}
\left\vert T_{m}a\left( x\right) \right\vert \leq \frac{cM_{l}M_{k}M_{N}}{m}.
\label{13q1}
\end{equation}

Let $x\in I_{N}^{k,N}.$ Then, according to Lemma \ref{lemma5a} we obtain that 
\begin{equation}
\left\vert T_{m}a\left( x\right) \right\vert \leq cM_{k}M_{N}.  \label{13q2}
\end{equation}

By combining (\ref{1.1}), (\ref{14aaa}), (\ref{13q1}) and (\ref{13q2}) we
obtain that%
\begin{eqnarray*}
	&&	\int_{\overline{I_{N}}}\left\vert T_{m}a\left( x\right) \right\vert
	^{1/2}d\mu \left( x\right)\\
	&\leq &c\overset{N-2}{\underset{k=0}{\sum }}\overset{N-1}{\underset{l=k+1}{%
			\sum }}\frac{m_{l+1}\dotsm m_{N-1}}{M_{N}}\frac{\left( M_{l}M_{k}\right)
		^{1/2}M_{N}^{1/2}}{m^{1/2}} \\
	&+&\overset{N-1}{\underset{k=0}{\sum }}\frac{1}{%
		M_{N}}M_{k}^{1/2}M_{N}^{1/2} \\  \notag
	&\leq & M_{N}^{1/2}\overset{N-2}{\underset{k=0}{\sum }}\overset{N-1}{\underset{%
			l=k+1}{\sum }}\frac{\left( M_{l}M_{k}\right) ^{1/2}}{m^{1/2}M_{l}}+\overset{%
		N-1}{\underset{k=0}{\sum }}\frac{M_{k}^{1/2}}{M_{N}^{1/2}} \\
	&\leq& \frac{%
		cM_{N}^{1/2}N}{m^{1/2}}+c.
\end{eqnarray*}

It follows that%
\begin{eqnarray}
&&\frac{1}{\log n}\overset{n}{\underset{m=M_{N}+1}{\sum }}\frac{\int_{%
		\overline{I_{N}}}\left\vert T_{m}a\left( x\right) \right\vert ^{1/2}d\mu
	\left( x\right) }{m} \\ \notag
&\leq& \frac{1}{\log n}\overset{n}{\underset{m=M_{N}+1}{%
		\sum }}\left( \frac{cM_{N}^{1/2}N}{m^{3/2}}+\frac{c}{m}\right) <c<\infty .
\label{15b}
\end{eqnarray}

The proof of part b) is completed by just combining (\ref{14b}) and (\ref%
{15b}).
\QED

\begin{corollary}
	Let $0<p\leq 1/2$ and $f\in H_{p}.$ Then there exists absolute constant $%
	c_{p},$ depending only on $p,$ such that the following inequality holds: 
	\begin{equation*}
	\frac{1}{\log ^{\left[ 1/2+p\right] }n}\overset{n}{\underset{k=1}{\sum }}%
	\frac{\left\Vert \sigma _{k}f\right\Vert _{p}^{p}}{k^{2-2p}}\leq
	c_{p}\left\Vert f\right\Vert _{H_{p}}^{p}.
	\end{equation*}
\end{corollary}

\begin{corollary}
	Let $0<p\leq 1/2$ and $f\in H_{p}.$ Then there exists an absolute constant $%
	c_{p},$ depending only on $p,$ such that the following inequality holds: 
	\begin{equation*}
	\frac{1}{\log ^{\left[ 1/2+p\right] }n}\overset{n}{\underset{k=1}{\sum }}%
	\frac{\left\Vert U_k^{\alpha}f\right\Vert _{p}^{p}}{k^{2-2p}}\leq
	c_{p}\left\Vert f\right\Vert _{H_{p}}^{p}.
	\end{equation*}
\end{corollary}

\begin{corollary}
	Let $0<p\leq 1/2$ and $f\in H_{p}.$ Then there exists an absolute constant $%
	c_{p},$ depending only on $p,$ such that the following inequality holds: 
	\begin{equation*}
	\frac{1}{\log ^{\left[ 1/2+p\right] }n}\overset{n}{\underset{k=1}{\sum }}%
	\frac{\left\Vert V_k^{\alpha}f\right\Vert _{p}^{p}}{k^{2-2p}}\leq
	c_{p}\left\Vert f\right\Vert _{H_{p}}^{p}.
	\end{equation*}
\end{corollary}

\begin{corollary}
	Let $0<p\leq 1/2$ and $f\in H_{p}.$ Then there exists an absolute constant $%
	c_{p},$ depending only on $p,$ such that the following inequality holds: 
	\begin{equation*}
	\overset{\infty}{\underset{k=1}{\sum }}
	\frac{\left\Vert R_k^{\alpha}f\right\Vert _{p}^{p}}{k^{2-2p}}\leq
	c_{p}\left\Vert f\right\Vert _{H_{p}}^{p}.
	\end{equation*}
\end{corollary}

\begin{theorem}
	\label{theorem2fejerstrong1}a) Let $0<p<1/2,$ $f\in H_{p}$ and $\{q_{k}:k\geq
	0\}$ be a sequence of non-decreasing numbers. Then there exists an absolute
	constant $c_{p},$ depending only on $p,$ such that the inequality%
	\begin{equation*}
	\overset{\infty }{\underset{k=1}{\sum }}\frac{\left\Vert T_{k}f\right\Vert
		_{p}^{p}}{k^{2-2p}}\leq c_{p}\left\Vert f\right\Vert _{H_{p}}^{p}
	\end{equation*}%
	holds.
	
	b)Let $f\in H_{1/2}$ and $\{q_{k}:k\geq 0\}$ be a sequence of non-increasing
	numbers, satisfying the condition \eqref{fn01}.
	Then there exists an absolute constant $c,$ such that the inequality%
	\begin{equation}
	\frac{1}{\log n}\overset{n}{\underset{k=1}{\sum }}\frac{\left\Vert
		T_{k}f\right\Vert _{1/2}^{1/2}}{k}\leq c\left\Vert f\right\Vert
	_{H_{1/2}}^{1/2}  \label{7nor}
	\end{equation}%
	holds.
\end{theorem}

\textbf{Proof:}[Proof of Theorem \protect\ref{theorem2fejerstrong1}]
If we use Lemmas \ref{lemma5aT} and \ref{lemma5b} and follows analogical steps of Theorem \ref{theorem2fejerstrong} we immediately get the proof of Theorem \protect\ref{theorem2fejerstrong1}. So, we leave out the details.
\QED

\begin{corollary}
	Let $0<p\leq 1/2,$ $f\in H_{p}$ and $\{q_{k}:k\geq 0\}$ be a sequence of
	non-decreasing numbers, such that
	\begin{equation*}
	\sup_{n\in \mathbb{N}}q_{n}<c<\infty .
	\end{equation*}%
	Then condition (\ref{fn01}) is satisfied and for such $T$ means there
	exists an absolute constant $c,$ such that the inequality (\ref{7nor}) holds.
\end{corollary}

We have already considered the case when the sequence $\{q_{k}:k\geq 0\}$ is
bounded. Now, we consider some Nörlund means, which are generated by a
unbounded sequence $\{q_{k}:k\geq 0\}.$

\begin{corollary}
	Let $0<p\leq 1/2$ and $f\in H_{p}.$ Then there exists an absolute constant $%
	c_{p},$ depending only on $p,$ such that the following inequality holds: 
	\begin{equation*}
	\frac{1}{\log ^{\left[ 1/2+p\right] }n}\overset{n}{\underset{k=1}{\sum }}%
	\frac{\left\Vert B^{\alpha,\beta}_kf\right\Vert _{p}^{p}}{k^{2-2p}}\leq
	c_{p}\left\Vert f\right\Vert _{H_{p}}^{p}.
	\end{equation*}
\end{corollary}

\newpage

\section{Reisz and N\"orlund logarithmoic means means on $H_p$ spaces}

\subsection{Introduction}

Riesz logarithmic means with respect to the trigonometric system was studied by a lot of authors. We mention, for instance, the paper by Szasz \cite{Sz} and Yabuta \cite{Ya}.
These means with respect to the Walsh and Vilenkin systems were investigated
by Simon \cite{Si11} and G\`{a}t \cite{Ga1}. Blahota  and G{\'a}t \cite{bg}
considered norm summability of N\"orlund logarithmic means and showed that
Riesz logarithmic means $R_{n}$ have better approximation properties on some
unbounded Vilenkin groups than the Fej\'er means. Moreover, in \cite{tep10}
it was proved that the maximal operator of Riesz means is bounded from the
Hardy space $H_{p}$ to the Lebesgue space $L_{p}$ for $p>1/2$ but not
when $0<p\leq 1/2.$ Strong convergence theorems and boundedness of weighted maximal operators of Riesz logarithmic means was considered in Lukkassen, Persson, Tutberidze, Tephnadze \cite{LPTT} and Tephnadze \cite{tep10}.

In \cite{tep11} Tephnadze proved that the maximal operator of Riesz logarithmic means $R^{\ast }$ is bounded from the Hardy space $H_{1/2}$  to the space $weak-L_{1/2}.$ Moreover,  there exists a martingale $f\in H_{p},$ where $0<p\leq 1/2$ such that
\begin{equation*}
\left\Vert R^{\ast}f\right\Vert_p=+\infty.
\end{equation*}

In \cite{tep10} Tephnadze proved that for any $0<p<1/2,$  the maximal operator 
$$\overset{\sim }{R}_{p}^{\ast }:=\sup_{n\in \mathbb{N}}\frac{\log n\left\vert R_n f\right\vert }{\left( n+1\right)^{1/p-2}}$$ 
is bounded from the Hardy space $H_p$ to the space $L_p.$

Moreover, for $0<p<1/2$ and  non-decreasing function $\varphi :\mathbb{N}_+\rightarrow
\lbrack 1,\infty )$ satisfying the condition
\begin{equation} \label{66}
\frac{\left( n+1\right)^{1/p-2}}{\log \left( n+1\right) \varphi \left(n\right) }=\infty, 
\end{equation}%
the maximal operator
\begin{equation*}
\sup_{n\in \mathbb{N}}\frac{\left\vert R_n f\right\vert }{\varphi \left(n\right) }
\end{equation*}
is not bounded from the Hardy space $H_p$ to the space $weak-L_p.$

In the case $p=1/2$ he also proved that the maximal operator 
$$\widetilde{R}^{\ast}f:=\sup_{n\in \mathbb{N}}\frac{\left\vert R_nf\right\vert}{\log\left(n+1\right)}$$ 
is bounded from the Hardy space $H_{1/2}$ to the space $L_{1/2}.$

Moreover, for any non-decreasing function  $\varphi:\mathbb{N}_+\rightarrow\lbrack 1,\infty)$ satisfying the condition
\begin{equation*}
\overline{\lim_{n\rightarrow\infty}}\frac{\log\left(n+1\right) }{\varphi\left(n\right)}=+\infty,  
\end{equation*}
the maximal operator
\begin{equation*}
\sup_{n\in \mathbb{N}}\frac{\left\vert R_nf\right\vert}{\varphi \left(n\right)}
\end{equation*}
is not bounded from the Hardy space $H_{1/2}$ to the space $L_{1/2}.$

In this thesis (see also \cite{LPTT}) we also proved that if $0<p<1/2$ and $f\in H_p(G_m),$ there exists an absolute constant $c_p,$ depending only on $p,$ such that the inequality
\begin{equation*} 
\overset{\infty}{\underset{n=1}{\sum}}\frac{\log^p n\left\Vert R_nf\right\Vert _{H_p}^p}{n^{2-2p}}\leq c_p\left\Vert f\right\Vert_{H_p}^p
\end{equation*}
holds, where $R_nf$ denotes the $n$-th Reisz logarithmic mean with respect to the Vilenkin-Fourier series of $f.$

M\'oricz and Siddiqi \cite{Mor} investigate the approximation properties of
some special N\"orlund means of Walsh-Fourier series of $L_{p}$ functions in
norm. The case when $\left\{ q_{k}=1/k:k\in \mathbb{N}\right\} $ was
excluded, since the methods of M\'oricz and Siddiqi are not applicable to N\"o%
rlund logarithmic means. 
Fridli, Manchanda and Siddiqi \cite{FMS} improved and extended results of M\'oricz and Siddiqi \cite{Mor} to dyadic homogeneous Banach spaces and Martingale Hardy spaces.
In \cite{Ga2} G\'{a}t and Goginava proved some
convergence and divergence properties of the N\"orlund logarithmic means of
functions in the class of continuous functions and in the Lebesgue space $%
L_{1}.$ In particular, they gave a negative answer to the question of M\'oricz
and Siddiqi \cite{Mor}. G\'{a}t and Goginava \cite{Ga3} proved that for each
measurable function satisfying
$
\phi \left( u\right) =o\left( u\log ^{1/2}u\right) , \  as \ 
u\rightarrow \infty ,
$
there exists an integrable function $f$ such that
$$
\int_{G_{m}}\phi \left( \left\vert f\left( x\right) \right\vert \right) d\mu
\,\left( x\right) <\infty
$$
and that there exists a set with positive measure such that the
Walsh-logarithmic means of the function diverges on this set.
It follows that that weak-(1,1) type inequality does not hold for the maximal operator of N\"orlund logarithmic means:
\index{S}{maximal operator of N\"orlund logarithmic means}
\begin{eqnarray*}
	L^{\ast}f&:=&\sup_{n\in\mathbb{N}}\left\vert L_nf\right\vert
\end{eqnarray*}

On the other hand, there exists an absolute constant $c_p$ such that
\begin{equation*}
\left\Vert L^{\ast }f\right\Vert _{p}\leq c_{p}\left\Vert f\right\Vert _{p},%
\text{ \ when \ }f\in L_{p},\text{ \ }p>1.
\end{equation*}

If we consider the following restricted maximal operator \index{S}{restricted maximal operator of partial sums}
\index{N}{$\widetilde{S}_{\#}^{\ast}f$}
\begin{eqnarray*}
	\widetilde{L}_{\#}^{\ast}f&:=&\sup_{n\in\mathbb{N}}\left\vert L_{M_n}f\right\vert, \ \ \ (M_{k}:=m_0...m_{k-1}, \ \ \  k=0,1...)
\end{eqnarray*}
then
\begin{eqnarray*}
	\lambda \mu\left\{ \widetilde{L}_{\#}^{\ast}f>\lambda \right\} \leq c\left\Vert f\right\Vert_{1}, \ \ \ f\in L_1(G_m), \ \ \lambda>0.
\end{eqnarray*}
Hence, if $f\in L_1(G_m)$ then
$$L_{M_n}f\to f, \ \ \text{a.e. on } \ \ G_m.$$

In this thesis we prove that if $f\in L_1(G_m)$ then $L_{M_n}f(x)\to f(x)$ for all Lebesgue points.

In \cite{tep4} it was proved that there exists a martingale $f\in H_{p},$ $(0<p\leq 1),$ such that the maximal operator of N\"orlund logarithmic means $L^{\ast }$ is not bounded in the Lebesgue space $L_{p}.$ In particular, it was proved that  there exists a martingale $f\in
H_{p} $ such that
\begin{equation*}
\left\Vert L^{\ast }f\right\Vert _{p}=+\infty .
\end{equation*}

Boundedness of weighted maximal operators of N\"orlund logarithmic means was considered  Persson, Tephnadze, Wall \cite{ptw3}. In particular,  maximal operator 
\begin{equation*}  \label{wemax}
\overset{\sim }{L}^{\ast }f:=\sup_{n\in \mathbb{N}}\frac{\left\vert
	L_{n}f\right\vert }{\log\left(n+1\right)}
\end{equation*}
is bounded from the Hardy space $H_{1}\left( G_{m}\right) $ to the space $%
L_{1}\left( G_{m}\right) .$

Moreover, if  $\varphi :\mathbb{N}_{+}\rightarrow \lbrack 1,\infty )$ be a non-decreasing function satisfying the condition%
\begin{equation} \label{5}
\overline{\lim_{n\rightarrow \infty }}\frac{\log \left( n+1\right) }{\varphi
	\left( n\right) }=+\infty ,
\end{equation}%
then there exists a martingale $f\in H_{1}\left( G_{m}\right),$ such that the maximal operator 
\begin{equation*}
\sup_{n\in \mathbb{N}}\frac{\left\vert L_{n}f\right\vert }{\varphi \left(
	n\right) }
\end{equation*}%
is not bounded from the Hardy space $H_{1}\left(G_{m}\right) $ to the
Lebesgue space $L_{1}\left( G_{m}\right) .$

In Tephnadze and Tutberidze \cite{tt1} was proved the maximal operator
\begin{equation*}
\overset{\sim}{L}_{p}^{*}f :=\sup_{n\in \mathbb{N}}
\frac{\left| L_n f \right|}{\left(n+1\right)^{1/p-1}}
\end{equation*}
is bounded from the Hardy space $H_{p}\left( G_{m}\right) $ to the space $L_{p}\left( G_{m}\right).$

We also proved that for $0<p<1$ and   a non-decreasing function $\varphi :\mathbb{N}_{+}\rightarrow [1, \infty)$ satisfying the condition
\begin{equation*} 
\overline{\lim_{n\rightarrow\infty}}\frac{n^{1/p-1}}{\log n\varphi\left( n\right)}=+\infty,  
\end{equation*}
then there exists a martingale $f\in H_p\left(G_m\right),$ such that the maximal operator
\begin{equation*}
\sup_{n\in \mathbb{N}}\frac{\left| L_{n}f \right| }{\varphi \left( n+1\right) }
\end{equation*}
is not bounded from the Hardy space $H_{p}\left( G_{m}\right) $ to the space $L_{p}\left( G_{m}\right).$

In this paper also state the following open problem:

\textbf{Open Problem.} For any  $0<p<1$ let find non-decreasing function $\Theta :\mathbb{N}_{+}\rightarrow [1, \infty)$ suth that the following maximal operator

\begin{equation*}
\overset{\sim}{L}_{p}^{*}f :=\sup_{n\in \mathbb{N}}
\frac{\left| L_n f \right|}{\Theta\left(n+1\right)}
\end{equation*}
is bounded from the Hardy space $H_{p}\left( G_{m}\right) $ to the Lebesgue space $L_p\left(G_m\right)$ and the rate of $\Theta :\mathbb{N}_{+}\rightarrow [1,\infty)$ is sharp, that is, for any non-decreasing function
$\varphi :\mathbb{N}_{+}\rightarrow [1, \infty)$ satisfying the condition
\begin{equation*} 
\overline{\lim_{n\rightarrow \infty }}\frac{\Theta\left( n\right)}{ \varphi\left( n\right)}=+\infty,  
\end{equation*}
then there exists a martingale $f\in H_p\left( G_m\right),$ such that the maximal operator
\begin{equation*}
\sup_{n\in\mathbb{N}}\frac{\left|L_nf\right|}{\varphi\left(n+1\right)}
\end{equation*}
is not bounded from the Hardy space $H_p\left(G_m\right)$ to the space
$L_p\left(G_m\right).$

According to Theorems above we can conclude that there exist absolute constants $C_1$ and $C_2$ such that
\begin{equation*}
\frac{C_1 n^{1/p-1}}{\log(n+1)}\leq\Theta\left(n\right) \leq C_2 n^{1/p-1}.
\end{equation*}

Later on, Memic generalized result of Tephnadze and Tutberidze \cite{tt1} and proved  that the maximal operator 
\begin{equation*}
\sup_{ n\in \mathbb{N}}
\frac{\log n\left| L_n f \right|}{\left(n+1\right)^{1/p-1}}
\end{equation*}
is bounded from the Hardy space $H_{p}\left( G_{m}\right) $ to the space $L_{p}\left( G_{m}\right).$

Sharpness of this result immediately follows negative result of Tephnadze and Tutberidze \cite{tt1}, which is already stated above.

\subsection{Auxiliary lemmas}

We need the following lemma of independent interest:

\begin{lemma}\label{reiszkernel}
	Let $n\in \mathbb{N}$. Then
	\begin{equation}\label{25}
	Y_n=\frac{1}{l_n}\overset{n-1}{\underset{j=1}{\sum}}\frac{K_j}{j+1} +\frac{K_n}{l_n}.  
	\end{equation}
	
	Moreover, 
	\begin{equation}\label{normreisz}
	\Vert Y_n \Vert_1 <c <\infty.
	\end{equation}	
\end{lemma}
\textbf{Proof:} To rewrite the kernels of the Riesz
logarithmic means as in equality \eqref{25} we have to just use Abel transformation. On the other hand, equality \eqref {normreisz} immediately follows equality \eqref{25} and \eqref{knbounded}.

\QED

To prove our main results we need the following lemma, which is proved in Tephnadze \cite{tep7}. Here we also give some short proof.

\begin{lemma} \label{l2} 
	Let $x\in I_{N}\left( x_{k}e_{k}+x_{l}e_{l}\right) ,$ $\ 1\leq x_{k}\leq
	m_{k}-1,$ $1\leq x_{l}\leq m_{l}-1,$ $\ k=0,...,N-2,$ $l=k+1,...,N-1.$ Then
	\begin{equation*}
	\int_{I_{N}}\underset{j=M_{N}+1}{\overset{n}{\sum }}\frac{\left|
		K_j\left( x-t\right) \right| }{j+1}d\mu \left( t\right) \leq \frac{
		cM_{k}M_{l}}{M_{N}^{2}}.
	\end{equation*}
	
	Let $x\in I_{N}\left( x_{k}e_{k}\right) ,$ $1\leq x_{k}\leq m_{k}-1,$ $\
	k=0,...,N-1.$ Then
	\begin{equation*}
	\int_{I_{N}}\underset{j=M_{N}+1}{\overset{n}{\sum }}\frac{\left| K_{j}\left(x-t\right) \right| }{j+1}d\mu \left( t\right) \leq \frac{cM_{k}}{M_{N}}l_{n}.
	\end{equation*}
\end{lemma}
\textbf{proof:} Let $x\in I_{N}\left( x_{k}e_{k}+x_{l}e_{l}\right) ,$ $\
1\leq x_{k}\leq m_{k}-1,$ $1\leq x_{l}\leq m_{l}-1,$ $\ k=0,...,N-2,$ $%
l=k+1,...,N-1.$ By using Lemma \ref{lemma5} we have that
\begin{eqnarray}\label{12aa}
&&\int_{I_{N}}\underset{j=M_{N}+1}{\overset{n}{\sum }}\frac{\left\vert
	K_{j}\left( x-t\right) \right\vert }{j+1}d\mu \left( t\right) \\ \notag
&\leq& \underset{j=M_{N}+1}{\overset{n}{\sum }}\frac{cM_{k}M_{l}}{\left( j+1\right) jM_{N}}
\\ \notag
&\leq&\frac{cM_{k}M_{l}}{M_{N}}\underset{j=M_{N}+1}{\overset{\infty }{\sum }
}\left( \frac{1}{j}-\frac{1}{j+1}\right) \\ \notag
&\leq& \frac{cM_{k}M_{l}}{M_{N}^{2}}.
\notag
\end{eqnarray}
Let $x\in I_{N}\left( x_{k}e_{k}\right) ,$ $1\leq x_{k}\leq m_{k}-1,$ $\
k=0,...,N-1.$ Then
\begin{eqnarray} \label{12b}
&&\int_{I_{N}}\underset{j=M_{N}+1}{\overset{n}{\sum }}\frac{\left\vert
	K_{j}\left( x-t\right) \right\vert }{j+1}d\mu \left( t\right) \\ \notag
&\leq& \underset{j=M_{N}+1}{\overset{n}{\sum }}\frac{cM_{k}}{\left( j+1\right) M_{N}} \\ \notag
&\leq&\frac{cM_{k}}{M_{N}}l_{n}.
\end{eqnarray}

Combining (\ref{12aa}) and (\ref{12b}) we complete the proof of Lemma.

\QED

Next we study some special consequences of kernels of N\"orlund logarithmic means:

\begin{lemma}\label{lemma0nnT121}Let $n\in \mathbb{N}$. Then
	\begin{eqnarray} \label{1.71alphaT2j} P_{M_n}(x)=D_{M_n}(x)-\psi_{M_n-1}(x)\overline{Y}_{M_n}(x)
	\end{eqnarray}
	Moreover,
	\begin{eqnarray} \label{1.71alphaT2jj} 
	\Vert P_{M_n}(x)\Vert_1<c<\infty.
	\end{eqnarray}
\end{lemma}
\textbf{Proof:} By using \ref{dn22} we get that	
\begin{eqnarray*}
	P_{M_n}(x)&=&\frac{1}{Q_{M_n}}\overset{M_n}{\underset{k=1}{\sum}}\frac{D_k(x)}{M_n-k}=\frac{1}{Q_{M_n}}\overset{M_n-1}{\underset{k=0}{\sum }}\frac{D_{M_n-k}(x)}{k}\\
	&=&\frac{1}{Q_{M_n}}\overset{M_n-1}{\underset{k=0}{\sum }}\frac{1}{k}\left(D_{M_n}(x)-\psi_{M_n-1}(x)\overline{D}_k(x)\right)\\
	&=&D_{M_n}(x)-\psi_{M_n-1}(x)\overline{Y}_{M_n}(x).
\end{eqnarray*}
which compete proof of identity \eqref{1.71alphaT2j}.

On the other hand, if we combine  \eqref{5aa} and \eqref{reisz} we also get proof of \eqref{1.71alphaT2jj}. The proof is complete.

\QED

The proof of next Lemma can be found in Tephnadze \cite{tep6}:
\begin{lemma} \label{dn2.7}
	Let $x\in I_s\backslash I_{s+1}, \ \  s=0,...,N-1.$  Then
	\begin{equation*}
	\int_{I_N}\left\vert P_n\left(x-t\right)\right\vert d\mu\left(
	t\right)\leq\frac{cM_s}{M_N},
	\end{equation*}%
	where $c$ is an absolute constant.
\end{lemma}

\textbf{Proof:} The proof is direct consequence of Lemma \ref{dn2.6}. So, we leave out the details.

\QED

\subsection{strong convergence of Reisz means means with respect to Vilenkin systems on the martingale Hardy spaces}

Our first main result reads:
\begin{theorem} \label{threisz_2} Let $0<p<1/2$ and $f\in H_p(G_m).$ Then there exists an absolute constant $c_{p},$ depending only on $p,$ such that the inequality
	\begin{equation} \label{star}
	\overset{\infty}{\underset{n=1}{\sum}}
	\frac{\log^p n \left\Vert R_nf\right\Vert _{H_p(G_m)}^p}{n^{2-2p}}\leq
	c_p\left\Vert f\right\Vert_{H_p(G_m)}^p
	\end{equation}
	holds, where $Rf$ denotes the $n$-th Reisz logarithmic mean with respect to the Vilenkin-Fourier series of $f.$ 
\end{theorem}

\textbf{Proof:} According to \eqref{normreisz} in Lemma \ref{reiszkernel} we get that
\begin{equation*}
\sup_{n\in \mathbb{N}}\int_{G_m}\left\vert Y_n a\right\vert d\mu \leq c<\infty.
\end{equation*}
and it follows that $R_n$ is bounded from $L_{\infty}$  to $L_{\infty}.$
By Proposition \ref{lemma2.2}, the proof of theorem will be complete, if we show that
\begin{equation} \label{reiszbounded}
\overset{\infty}{\underset{n=1}{\sum}}\frac{\log^pn \int\limits_{\overset{-}{I}}\left\vert R_na\right\vert^pd\mu}{n^{2-2p}}\leq c_p<\infty,\ \ \text{for} \ \  0<p<1/2,
\end{equation}
for every $p$-atom $a,$ where $I$ denotes the support of the atom.

Let $a$ be an arbitrary p-atom with support $I$ and $\mu \left( I\right)
=M_{N}^{-1}.$ We may assume that $I=I_{N}.$ It is easy to see that $%
R_n a=\sigma_n\left(a\right)=0,$ when $n\leq M_{N}$.
Therefore we suppose that $n>M_{N}.$

Since $\left\Vert a\right\Vert _{\infty }\leq cM_{N}^{2}$ if we apply (\ref{25}) in Lemma \ref{reiszkernel}, we can conclude that
\begin{eqnarray} \label{26}
&&\left\vert R_na\left(x\right)\right\vert \\ \notag
&=&\int_{I_N}\left\vert a\left(t\right)\right\vert\left\vert Y_n\left(x-t\right)\right\vert d\mu(t)   \\
&\leq &\left\Vert a\right\Vert _{\infty}\int_{I_{N}}\left\vert Y_n\left( x-t\right) \right\vert d\mu \left(t\right)  \notag \\
&\leq &\frac{cM_{N}^{1/p}}{l_n}\underset{I_{N}}{\int }
\underset{j=M_{N}+1}{\overset{n-1}{\sum }}\frac{\left\vert
	K_{j}\left( x-t\right) \right\vert }{j+1}d\mu \left( t\right)  \notag \\
&+&\frac{cM_N^{1/p}}{l_n}\int_{I_N}\left\vert
K_n\left(x-t\right)\right\vert d\mu\left(t\right).  \notag
\end{eqnarray}

Let $x\in I_{N}\left( x_{k}e_{k}+x_{l}e_{l}\right) ,$ $\ 1\leq x_{k}\leq
m_{k}-1,$ $1\leq x_{l}\leq m_{l}-1,$ $k=0,...,N-2,$ $l=k+1,...,N-1.$ From
Lemma  \ref{l2} it follows that
\begin{equation} \label{29}
\left| R_{n}a\left( x\right)\right| \leq\frac{cM_{l}M_{k}M_{N}^{1/p-2}}{\log(n+1)}.  
\end{equation}

Let $x\in I_{N}\left( x_{k}e_{k}\right) ,$ $1\leq x_{k}\leq m_{k}-1,$ $%
k=0,...,N-1.$ Applying Lemma \ref{l2} we can conclude that
\begin{equation} \label{32}
\left\vert R_{n}a\left(x\right)\right\vert\leq{M_{N}^{1/p-1}M_{k}}.  
\end{equation}

By combining (\ref{1.1}) and (\ref{26}-\ref{32}) we obtain that
\begin{eqnarray} \label{7aaa}
&&\int_{\overline{I_{N}}}\left\vert R_na\left( x\right) \right\vert
^{p}d\mu \left( x\right) \\ \notag
&=&\overset{N-2}{\underset{k=0}{\sum}}\overset{N-1}{\underset{l=k+1}{\sum}}\sum\limits_{x_j=0,j\in\{l+1,...,N-1}^{m_{j-1}}\int_{I_N^{k,l}}\left\vert R_na \right\vert^{p}d\mu \\ \notag
&+&\overset{N-1}{\underset{k=0}{\sum}} \int_{I_N^{k,N}}\left\vert R_na\right\vert^pd\mu   \\ \notag
&\leq&c\overset{N-2}{\underset{k=0}{\sum }}\overset{N-1} {\underset{l=k+1}{\sum }}\frac{m_{l+1}\dots m_{N-1}}{M_{N}}\frac{\left( M_{l}M_{k}\right)^{p}M_{N}^{1-2p}}{{\log}^{p}(n+1)}\\ \notag
&+&\overset{N-1}{\underset{k=0}{\sum }}\frac{1}{M_{N}}M_{k}^{p}M_{N}^{1-p}  \\
&\leq&\frac{cM_N^{1-2p}}{{\log^p(n+1)}}\overset{N-2}{\underset{k=0}{\sum}}\overset{N-1}{\underset{l=k+1}{\sum}}\frac{\left(M_lM_k\right)^p}{M_l}+\overset{N-1}{\underset{k=0}{\sum}}\frac{M_{k}^{p}}{M_{N}^p}  \notag \\
&=&\frac{cM_N^{1-2p}}{\log^p(n+1)}\overset{N-2}{\underset{k=0}{\sum}}\frac{1}{M_k^{1-2p}}\overset{N-1}{\underset{l=k+1}{\sum}}\frac{M_k^{1-p}}{M_l^{1-p}}+\overset{N-1}{\underset{k=0}{\sum}}\frac{M_k^p}{M_{N}^p} \notag \\ &\leq&\frac{cM_N^{1-2p}}{\log^p(n+1)}+c_p.  \notag
\end{eqnarray}

It is easy to see that
\begin{equation} \label{6aaa}
\overset{\infty}{\underset{n=M_{N}+1}{\sum }}\frac{1}{n^{2-2p}}\leq \frac{c}{M_{N}^{1-2p}}, \ \ \text{ for } \ \ 0<p< 1/2.
\end{equation}

By combining (\ref{7aaa}) and (\ref{6aaa}) we get that
\begin{eqnarray*}
	&&\overset{\infty}{\underset{n=M_N+1}{\sum}}\frac{\log^p n\int_{\overline{I_N}}\left\vert R_na\right\vert ^{p}d\mu}{n^{2-2p}}  \\
	&\leq &\overset{\infty}{\underset{n=M_N+1}{\sum}}\left( \frac{c_pM_N^{1-2p}}{n^{2-p}}+\frac{c_p}{n^{2-p}}\right)+c_p \\
	&\leq&c_pM_N^{1-2p}\overset{\infty}{\underset{n=M_{N}+1}{\sum}}\frac{1}{n^{2-2p}}\\
	&+&\overset{\infty}{\underset{n=M_N+1}{\sum }}\frac{1}{n^{2-p}}+c_p\leq C_p<\infty.
\end{eqnarray*}
It means that (\ref{reiszbounded}) holds true and the proof is complete.

\QED
Our next main result shows in particular that the inequality in Theorem \ref{threisz_2} is in a special sense sharp at least in the case of Walsh-Fourier series.
\begin{theorem} \label{reisz_negative_th_2}
	Let $0<p<1/2$ and $\Phi :\mathbb{N}\rightarrow \lbrack 1,\infty )$ be any non-decreasing function, satisfying the condition
	\begin{equation}\label{dir222}
	\underset{n\rightarrow \infty }{\lim }\Phi \left( n\right) =+\infty .
	\end{equation}
	Then there exists a martingale $f\in H_p\left(G_2\right) $ such that
	\begin{equation}\label{theoremjemala}
	\sum_{n=1}^{\infty}\frac{\log^p n
		\left\Vert R^w_nf\right\Vert_p^p\Phi\left(n\right)}{n^{2-2p}} =\infty,
	\end{equation}
	where $R^w_{n}f$ denotes the $n$-th Reisz logarithmic means with respect to Walsh-Fourier series of
	$f.$ 
\end{theorem}

\textbf{Proof:}
It is evident that if we assume that $\Phi\left( n\right)\geq cn,$ where $c$ is some positive constant then 
\begin{equation*}
\frac{\log^p n\Phi\left(n\right)}{n^{2-2p}}\geq n^{1-2p}\log^p n \rightarrow \infty, \ \text{as} \ n\rightarrow \infty,
\end{equation*}
and also \eqref{theoremjemala} holds. So, without lost the generality we may assume that 	there exists an increasing sequence of positive integers
$\left\{ \alpha^{\prime} _{k}:k\in \mathbb{N} \right\} $ such that
\begin{equation} \label{3jemala}
\Phi\left( \alpha^{\prime} _{k}\right)=o(\alpha^{\prime} _{k}), \ \text{as} \ k\rightarrow \infty.
\end{equation}

Let $\left\{ \alpha _{k}:k\in \mathbb{N} \right\}\subseteq \left\{ \alpha^{\prime} _{k}:k\in \mathbb{N} \right\}$ be an increasing sequence of positive integers such that $\alpha _{0}\geq 2$ and
\begin{equation}
\sum_{k=0}^{\infty }\frac{1}{\Phi^{1/2}(2^{2\alpha_k})}<\infty,  \label{3a}
\end{equation}

\begin{equation} \label{4a}
\sum_{\eta =0}^{k-1}\frac{2^{2\alpha _{\eta }/p}} {\Phi^{1/2p}(2^{2\alpha_\eta})}\leq \frac{2^{2\alpha _{k-1}/p+1}}{\Phi^{1/2p}(2^{2\alpha_{k-1}})}, 
\end{equation}

\begin{equation} \label{5a}
\frac{2^{2\alpha _{k-1}/p+1}}{\Phi^{1/2p}(2^{2\alpha_{k-1}})}\leq \frac{1}{128\alpha_k}\frac{2^{2\alpha _k(1/p-2)}}{\Phi^{1/2p}(2^{2\alpha_k})}.  
\end{equation}

We note that under condition \eqref{3jemala} we can conclude that
$$\frac{2^{2\alpha _{\eta }/p}} {\Phi^{1/2p}(2^{2\alpha_\eta})}
\geq
{\left(\frac{2^{2\alpha _{\eta }}} {\Phi(2^{2\alpha_\eta})}\right)}^{1/2p} \rightarrow \infty, \ \text{as} \ \eta\rightarrow \infty$$
and it immediately follows that such an increasing sequence $\left\{ \alpha_k:k\in \mathbb{N}\right\} ,$ which satisfies conditions (\ref{3a})-(\ref{5a}), can be constructed.

Let 
\begin{equation*}
f^{\left( n\right) }\left( x\right) :=\sum_{\left\{ k;\text{ }2\alpha_{k}<n\right\} }\lambda _{k}a_{k},
\end{equation*}
where 
\begin{equation*}
\lambda_k=\frac{1}{\Phi^{1/2p}(2^{2\alpha_k})}
\end{equation*}
and
\begin{equation*}
a_{k}={2^{2\alpha _{k}(1/p-1)}}\left( D_{2^{2\alpha _{k}+1}}-D_{2^{2\alpha_{k}}}\right) .
\end{equation*}

From (\ref{3a}) and Lemma \ref{lemma2.1} we can conclude that $f=\left( f^{\left(n\right) },n\in \mathbb{N} \right) \in H_p(G_2).$

It is easy to show that
\begin{equation} \label{8aa}
\widehat{f}^w(j)=\left\{ 
\begin{array}{l}
\frac{2^{2\alpha _{k}(1/p-1)}}{\Phi^{1/2p}(2^{2\alpha_k})},\,\,\text{ if \thinspace\thinspace }j\in \left\{ 2^{2\alpha _{k}},...,\text{ ~}2^{2\alpha_{k}+1}-1\right\} ,\text{ }k\in \mathbb{N}, \\ 
0,\text{ \thinspace \thinspace \thinspace if \thinspace \thinspace
	\thinspace }j\notin \bigcup\limits_{k=1}^{\infty }\left\{ 2^{2\alpha
	_k},...,\text{ ~}2^{2\alpha _{k}+1}-1\right\}.
\end{array}\right.  
\end{equation}

For $n=\sum_{i=1}^s 2^{n_i},$ $n_1<n_2<...<n_s$ we denote
\begin{equation*}
\mathbb{A}_{0,2}:=\left\{ n\in \mathbb{N}:\text{ }n=2^0+2^2+
\sum_{i=3}^{s_n}2^{n_i}\right\}.
\end{equation*}

Let $2^{2\alpha _{k}}\leq j\leq 2^{2\alpha_{k}+1}-1$ and $j\in \mathbb{A}_{0,2}.$ Then 
\begin{eqnarray} \label{10a}
R^w_jf &=&\frac{1}{l_j}\sum_{n=1}^{2^{2\alpha_k}-1}\frac{S_nf}{n} +\frac{1}{l_j}\sum_{n=2^{2\alpha_k}}^{j}\frac{S_{n}f}{n}:=I+II. 
\end{eqnarray}

Let $n<2^{2\alpha _{k}}.$ Then from (\ref{4a}), (\ref{5a}) and (\ref{8aa}) we have that
\begin{eqnarray*}
	\left\vert S^w_{n}f\left( x\right) \right\vert 
	&\leq &\sum_{\eta =0}^{k-1}\sum_{v=2^{2\alpha _{\eta }}}^{2^{2\alpha _{\eta}+1}-1}\left\vert \widehat{f}^w(v)\right\vert \\
	&\leq& \sum_{\eta=0}^{k-1}\sum_{v=2^{2\alpha _{\eta }}}^{2^{2\alpha _{\eta }+1}-1}\frac{2^{2\alpha _{\eta }(1/p-1)}}{\Phi^{1/2p}(2^{2\alpha_\eta})} \\
	&\leq &\sum_{\eta =0}^{k-1}\frac{2^{2\alpha _{\eta }/p}} {\Phi^{1/2p}(2^{2\alpha_\eta})}\\
	&\leq& \frac{2^{2\alpha _{k-1}/p+1}}{\Phi^{1/2p}(2^{2\alpha_{k-1}})}\\
	&\leq &\frac{1}{128\alpha_k}\frac{2^{2\alpha _k(1/p-2)}}{\Phi^{1/2p}(2^{2\alpha_k})}.
\end{eqnarray*}

Consequently,	
\begin{eqnarray} \label{11a}
\left\vert I\right\vert &\leq &\frac{1}{l_j}\underset{n=1}
{\overset{2^{2\alpha_k}-1}{\sum}}\frac{\left\vert S^w_n f\left( x\right)\right\vert}{n} \\ \notag
&\leq&\frac{1}{l_{2^{2\alpha_k}}}\frac{1}{128\alpha_k}\frac{2^{2\alpha _k(1/p-2)}}{\Phi^{1/2p}(2^{2\alpha_k})}
\sum_{n=1}^{2^{2\alpha_k}-1} \frac{1}{n}\\
&\leq&\frac{1}{128\alpha_k}\frac{2^{2\alpha _k(1/p-2)}}{\Phi^{1/2p}(2^{2\alpha_k})}.  \notag
\end{eqnarray}

Let $2^{2\alpha _{k}}\leq n\leq 2^{2\alpha_{k}+1}-1.$ Then we have the
following
\begin{eqnarray*}
	S^w_{n}f &=&\sum_{\eta =0}^{k-1}\sum_{v=2^{2\alpha _{\eta }}}^{2^{2\alpha
			_{\eta }+1}-1}\widehat{f}^w(v)w_{v}+\sum_{v=2^{2\alpha _{k}}}^{n-1}\widehat{f}^w(v)w_{v} \\
	&=&\sum_{\eta=0}^{k-1}\frac{2^{{2\alpha _{\eta }}\left(1/p-1\right)}} {\Phi^{1/2p}(2^{2\alpha_\eta})}\left( D^w_{2^{2\alpha _{\eta}+1}} -D^w_{2^{2\alpha_{\eta}}}\right)\\
	&+&\frac{2^{{2\alpha_k}\left( 1/p-1\right)}} {\Phi^{1/2p}(2^{2\alpha_k})}\left(D^w_{n}-D^w_{2^{2\alpha_k}}\right) .
\end{eqnarray*}

This gives that

\begin{eqnarray}\label{12a}
II&=&\frac{1}{l_j}\underset{n=2^{2\alpha_k}}{\overset{2^{2\alpha _{k}+1}}{\sum }}\ \frac{1}{n}\left(\sum_{\eta=0}^{k-1}\frac{2^{2\alpha _{\eta}\left(1/p-1\right)}}{\Phi^{1/2p}(2^{2\alpha_\eta})}\left( D^w_{2^{2\alpha_{\eta}+1}}-D^w_{2^{2\alpha _{\eta }}}\right)\right)  \\ \notag
&+&\frac{1}{l_{j}}\frac{2^{2\alpha_k\left(1/p-1\right)}}{\Phi^{1/2p}(2^{2\alpha_k})}\sum_{n=2^{2\alpha_k}}^{j}\frac{\left( D^w_n-D^w_{2^{2\alpha _k}}\right) }{n}\\ \notag
&:=&II_{1}+II_{2}.
\end{eqnarray}

Let $x\in I_{2}(e_{0}+e_{1})\in I_{0}\backslash I_{1}.$ According to 
\eqref{1dn} and \eqref{2dn} we conclude that
\begin{equation*}
D^w_n\left( x\right) =\left\{ 
\begin{array}{ll}
w_{n}, & \,\text{if\thinspace \thinspace \thinspace }n\ \ \text{is odd
	number,} \\ 
0, & \text{if}\,\,n\ \ \text{is even number.}
\end{array}%
\right.
\end{equation*}

Since $\alpha _{0}\geq 2,\ \ k\in \mathbb{N}$ we obtain that $2\alpha_{k}\geq 4,$ for all $k\in \mathbb{N}$ and if we apply (\ref{1dn}) we get that
\begin{equation}
II_{1}=0  \label{13a}
\end{equation}
and 
\begin{eqnarray*}
	II_2&=&\frac{1}{l_j}\frac{2^{2\alpha _k(1/p-1)}}{\Phi^{1/2p}(2^{2\alpha_k})}\sum_{n=2^{2\alpha _k-1}}^{(j-1)/2}\frac{w_{2n+1}}{2n+1} \\
	&=&\frac{1}{l_j}\frac{2^{2\alpha _k(1/p-1)}r_{1}}{\Phi^{1/2p}(2^{2\alpha_k})}\sum_{n=2^{2\alpha _k-1}}^{(j-1)/2}\frac{w_{2n}}{2n+1}.
\end{eqnarray*}

Let $x\in I_{2}(e_{0}+e_{1}).$ Then, by the definition of Walsh functions, we get that
\begin{equation*}
w_{4n+2}=r_{1}w_{4n}=-w_{4n}
\end{equation*}
and
\begin{eqnarray} \label{14aa}
&&\left\vert II_2\right\vert=\frac{1}{l_j}\frac{2^{2\alpha _k(1/p-1)}}{\Phi^{1/2p}(2^{2\alpha_k})}\left\vert \sum_{n=2^{2\alpha _k-1}}^{(j-1)/2}\frac{w_{2n}}{2n+1}\right\vert \\ \notag
&=&\frac{1}{l_j}\frac{2^{2\alpha_k(1/p-1)}}{\Phi^{1/2p}(2^{2\alpha_k})}\left\vert\frac{w_{j-1}}{j}+\sum_{n=2^{2\alpha_k-2}+1}^{(j-1)/4}\left(\frac{w_{4n-4}}{4n-3}+\frac{w_{4n-2}}{4n-1}\right)\right\vert  \\ \notag
&=&\frac{1}{l_j}\frac{2^{2\alpha _{k}(1/p-1)}}{\Phi^{1/2p}(2^{2\alpha_k})}\left\vert \frac{w_{j-1}}{j}+\sum_{n=2^{2\alpha_{k}-2}+1}^{(j-1)/4}\left( \frac{w_{4n-4}}{4n-3}-\frac{w_{4n-2}}{4n-1}\right) \right\vert  
\end{eqnarray}
\begin{eqnarray*}
	&\geq&\frac{c}{\log({2^{2\alpha_k+1}})}\frac{2^{2\alpha_k(1/p-1)}}{\Phi^{1/2p}(2^{2\alpha_k})}\left( \left\vert\frac{w_{j-1}}{j}\right\vert-\sum_{n=2^{2\alpha _{k}-2}+1}^{(j-1)/4} \left\vert{w_{4n-4}}\right\vert\left( \frac{1}{4n-3}-\frac{1}{4n-1}\right) \right)  \\ \notag
	&&\geq\frac{1}{4\alpha_k}\frac{2^{2\alpha_k(1/p-1)}}{\Phi^{1/2p}(2^{2\alpha_k})}\left( \frac{1}{j}-\sum_{n=2^{2\alpha_k-2}+1}^{(j-1)/4}\left( \frac{1}{4n-3}-\frac{1}{4n-1}\right)\right).
\end{eqnarray*}

Hence, by simple calculation we can conclude that
\begin{eqnarray*}
	&&\sum_{n=2^{2\alpha_k-2}+1}^{(j-1)/4}\left(\frac{1}{4n-3}- \frac{1} {4n-1}\right)\\
	&=&\sum_{n=2^{2\alpha_k-2}+1}^{(j-1)/4}\frac{2}{(4n-3)(4n-1)} \\
	&\leq& \sum_{n=2^{2\alpha _k-2}+1}^{(j-1)/4} \frac{2}{(4n-4)(4n-2)}\\
	&=&\frac{1}{2}\sum_{n=2^{2\alpha _k-2}+1}^{(j-1)/4} \frac{1}{(2n-2)(2n-1)} \\
	&\leq&\frac{1}{2}\sum_{n=2^{2\alpha _k-2}+1}^{(j-1)/4} \frac{1}{(2n-2)(2n-2)}\\
	&=&\frac{1}{8}\sum_{n=2^{2\alpha _k-2}+1}^{(j-1)/4} \frac{1}{(n-1)(n-1)} \\
	&\leq& \frac{1}{8}\sum_{n=2^{2\alpha _k-2}+1}^{(j-1)/4} \frac{1}{(n-1)(n-2)}\\
	&=&\frac{1}{8}\sum_{l=2^{2\alpha _k-2}+1}^{(j-1)/4} \left(\frac{1}{n-2}-\frac{1}{n-1} \right)\\
	&\leq &\frac{1}{8}\left(\frac{1}{2^{2\alpha_k-2}-1}-\frac{4}{j-5}\right)\\
	&\leq& \frac{1}{8}\left(\frac{1}{2^{2\alpha_k-2}-1}-\frac{4}{j}\right).
\end{eqnarray*}

Since $2^{2\alpha _{k}}\leq j\leq 2^{2\alpha_{k}+1}-1,$ where $\alpha _{k}\geq 2,$ we obtain that 
$$\frac{2}{2^{2\alpha_k}-4}\leq \frac{2}{2^{4}-4}=\frac{1}{6}$$
and

\begin{eqnarray} \label{14aR}
\left\vert II_2\right\vert&\geq&\frac{1}{4\alpha_k}\frac{2^{2\alpha _k(1/p-1)}}{\Phi^{1/2p}(2^{2\alpha_k})}\left( \frac{1}{j}-\frac{1}{8}\left(\frac{1}{2^{2\alpha_k-2}-1}-\frac{4}{j}
\right)\right) \\ \notag
&\geq &\frac{1}{4\alpha_k}\frac{2^{2\alpha _k(1/p-1)}}{\Phi^{1/2p}(2^{2\alpha_k})}\left( \frac{3}{2j}-\frac{1}{2^{2\alpha_k+1}-8}\right) \\ \notag
&\geq &\frac{1}{4\alpha_k}\frac{2^{2\alpha _k(1/p-1)}}{\Phi^{1/2p}(2^{2\alpha_k})}\left( \frac{3}{4}\frac{1}{2^{2\alpha_{k}}}-\frac{1}{2}\frac{1}{2^{2\alpha_k}-4}\right) \\ \notag
&\geq &\frac{1}{4\alpha_k}\frac{2^{2\alpha _k(1/p-1)}}{\Phi^{1/2p}(2^{2\alpha_k})}\left( \frac{1}{4}\frac{1}{2^{2\alpha_{k}}}+\frac{1}{2}\frac{1}{2^{2\alpha_{k}}}-\frac{1}{2}\frac{1}{2^{2\alpha_k}-4}\right) \\ \notag
&=&\frac{1}{4\alpha_k}\frac{2^{2\alpha _k(1/p-1)}}{\Phi^{1/2p}(2^{2\alpha_k})}\left( \frac{1}{4}\frac{1}{2^{2\alpha_{k}}}-\frac{2}{2^{2\alpha_k}(2^{2\alpha_k}-4)}\right)
\\ \notag
&\geq &\frac{1}{4\alpha_k}\frac{2^{2\alpha _k(1/p-1)}}{\Phi^{1/2p}(2^{2\alpha_k})}\left( \frac{1}{4}\frac{1}{2^{2\alpha_{k}}}-\frac{1}{6}\frac{1}{2^{2\alpha_{k}}}\right)
\\ \notag
&\geq &\frac{1}{48\alpha_k}\frac{2^{2\alpha _k(1/p-2)}}{\Phi^{1/2p}(2^{2\alpha_k})}\\ \notag
&\geq &\frac{1}{64\alpha_k}\frac{2^{2\alpha _k(1/p-2)}}{\Phi^{1/2p}(2^{2\alpha_k})}.
\end{eqnarray}

By combining (\ref{5}), (\ref{10a})-(\ref{14aR}) for $\in I_{2}(e_{0}+e_{1})$ and $0<p<1/2$ we find that

\begin{eqnarray*}
	&&\left\vert R^w_j f\left( x\right) \right\vert \\
	&\geq&
	\left\vert II_{2}\right\vert-\left\vert II_{1}\right\vert -\left\vert I \right\vert\\
	&\geq &\frac{1}{64\alpha_k}\frac{2^{2\alpha _k(1/p-2)}}{\Phi^{1/2p}(2^{2\alpha_k})}
	-\frac{1}{128\alpha_k}\frac{2^{2\alpha _k(1/p-2)}}{\Phi^{1/2p}(2^{2\alpha_k})}\\
	&=&\frac{1}{128\alpha_k}\frac{2^{2\alpha _k(1/p-2)}}{\Phi^{1/2p}(2^{2\alpha_k})}.
\end{eqnarray*}
and
\begin{eqnarray} \label{16}
&&\left\Vert R^w_{j}f\right\Vert_{weak-L_p(G_2)}^p  \\ \notag
&\geq &\frac{1}{128\alpha_k^p}\frac{2^{2\alpha _k(1-2p)}}{\Phi^{1/2}(2^{2\alpha_k})}\mu\left\{ x\in G_2:\left\vert R^w_{j}f\right\vert \geq \frac{1}{128\alpha_k}\frac{2^{2\alpha _k(1/p-2)}}{\Phi^{1/2p}(2^{2\alpha_k})}\right\} ^{1/p} \\ \notag
&\geq &\frac{1}{128\alpha_k^p}\frac{2^{2\alpha _k(1-2p)}}{\Phi^{1/2}(2^{2\alpha_k})}\mu \left\{ x\in I_{2}(e_{0}+e_{1}):\left\vert
R^w_{j}f\right\vert \geq \frac{1}{128\alpha_k}\frac{2^{2\alpha _k(1/p-2)}}{\Phi^{1/2p}(2^{2\alpha_k})}\right\} \\ \notag
&\geq &\frac{1}{128\alpha_k^p}\frac{2^{2\alpha _k(1-2p)}}{\Phi^{1/2}(2^{2\alpha_k})}(\mu \left( x\in I_{2}(e_{0}+e_{1})\right))\\ \notag
&>&\frac{1}{516\alpha_k^p}\frac{2^{2\alpha _k(1-2p)}}{\Phi^{1/2}(2^{2\alpha_k})}.
\end{eqnarray}
Moreover,

\begin{eqnarray*}
	&&\underset{j=1}{\overset{\infty }{\sum }}\frac{\left\Vert R^w
		_{j}f\right\Vert _{weak-L_p(G_2)}^{p}\log^p{(j)}\Phi(j)}{j^{2-2p} } \\
	&\geq& \underset{\left\{ j\in \mathbb{A}_{0,2}:\text{ }2^{2\alpha_k}< j\leq 2^{2\alpha_{k}+1}-1\right\} }{\sum }\frac{\left\Vert R^w
		_{j}f\right\Vert _{weak-L_p}^{p}\log^p{(j)}\Phi(j)}{j^{2-2p} } \\
	&\geq &\frac{c}{\alpha_k^p}\frac{2^{2\alpha _k(1-2p)}}{\Phi^{p/2}(2^{2\alpha_k)}}\underset{\left\{ j\in \mathbb{A}_{0,2}:\text{ }2^{2\alpha_k}< j\leq 2^{2\alpha_{k}+1}-1\right\} }{\sum }\frac{\log^p{(j)}\Phi(j)}{j^{2-2p}}\\
	&\geq &\frac{c\Phi(2^{2\alpha_{k}})\log^p{(2^{2\alpha_k})}}{\alpha_k^p}\frac{2^{2\alpha _k(1-2p)}}{\Phi^{1/2}(2^{2\alpha_k})}\underset{\left\{ j\in \mathbb{A}_{0,2}:\text{ }2^{2\alpha_k}< j\leq 2^{2\alpha_{k}+1}-1\right\} }{\sum }\frac{1}{j^{2-2p}}\\ 
	&\geq& \Phi^{1/2}(2^{2\alpha_{k}})\rightarrow \infty ,\text{ \ as \ \ }k\rightarrow \infty.
\end{eqnarray*}

The proof is complete.
\QED

\begin{theorem}\label{Corollaryconv5} Let $p\geq 1$ and $f\in L_p$. Then
	\begin{equation}\label{nata}
	\Vert R_{n} f-f\Vert_p \to 0 \ \ \text{as}\ \ n\to \infty
	\end{equation}
	for all  $f\in L_p(G_m)$. Moreover,
	\begin{equation*}
	\underset{n\rightarrow \infty }{\lim }R_{n}f(x)=f(x)
	\end{equation*}
	for all Vilenkin-Lebesgue points of $f\in L_p(G_m)$.
\end{theorem}

\textbf{Proof:} By using Abel transformation we get 

\begin{equation}\label{25a}
R_n=\frac{1}{l_n}\overset{n-1}{\underset{j=1}{\sum}}\frac{\sigma_j}{j+1} +\frac{\sigma_n}{l_n}.  
\end{equation}
Since 	(for details see Yano \cite{Yano})
$$\Vert \sigma_{n} f-f\Vert_p \to 0 \ \ \text{as}\ \ n\to \infty$$
if we apply \eqref{25a} we immediately get \eqref{nata}.

On the other hand, since (for details see \cite{gg})
\begin{equation*}
\underset{n\rightarrow \infty }{\lim }\sigma_{n}f(x)=f(x)
\end{equation*}
for all Vilenkin-Lebesgue points of $f\in L_p(G_m)$, if we invoke again equality \eqref{25a} we also get 
\begin{equation*}
\underset{n\rightarrow \infty }{\lim }R_{n}f(x)=f(x)
\end{equation*}
for all Vilenkin-Lebesgue points of $f\in L_p(G_m)$.
\QED

\subsection{N\"orlund logarithmic means of Vilenkin-Fourier series on Martingale Hardy spaces}

\begin{theorem}
	a) Let $0<p<1.$ Then the maximal operator
	\begin{equation*}
	\overset{\sim}{L}_{p}^{*}f :=\sup_{n\in \mathbb{N}}
	\frac{\left| L_n f \right|}{\left(n+1\right)^{1/p-1}}
	\end{equation*}
	is bounded from the Hardy space $H_{p}\left( G_{m}\right) $ to the space $%
	L_{p}\left( G_{m}\right).$
	
	b) Let $0<p<1$ and $\varphi :\mathbb{N}_{+}\rightarrow [1, \infty)$ be a
	non-decreasing function satisfying the condition
	\begin{equation} \label{6}
	\overline{\lim_{n\rightarrow\infty}}\frac{n^{1/p-1}}{\log n\varphi\left( n\right)}=+\infty.  
	\end{equation}
	Then there exists a martingale $f\in H_p\left(G_m\right),$ such that the maximal operator
	\begin{equation*}
	\text{ }\sup_{n\in \mathbb{N}}\frac{\left| L_{n}f 
		\right| }{\varphi \left( n+1\right) }
	\end{equation*}
	is not bounded from the Hardy space $H_{p}\left( G_{m}\right) $ to the space
	$L_{p}\left( G_{m}\right).$
\end{theorem}

\textbf{Proof:} Since 
\begin{eqnarray*}
	&&\frac{\left|L_{n} f \right|}{\left(n+1\right)^{1/p-1}}\\
	&\leq& \frac{1}{\left(n+1\right)^{1/p-1}}\underset{1\leq k\leq n}{\sup }\left|S_k f  \right| \\
	&\leq& \underset{1\leq k\leq n}{\sup }\frac{\left|S_k f \right|}{\left(k+1\right)^{1/p-1}} \\
	&\leq& \underset{n\in \mathbb{N}}{\sup }\frac{\left|S_n f  \right|}{\left(n+1\right)^{1/p-1}}
\end{eqnarray*}
and
\begin{equation}\label{LnSn}
\underset{n\in \mathbb{N}}{\sup }\frac{\left|L_n f \right|}{\left(n+1\right)^{1/p-1}}
\leq \underset{n\in \mathbb{N}}{\sup }\frac{\left|S_n f  \right|}{\left(n+1\right)^{1/p-1}}
\end{equation}
On the other hand, Tephnadze \cite{tep7} (see also \cite{tep9} and \cite{tep12}) proved that for every $0<p<1,$  the  maximal operator
\begin{equation*}
\widetilde{S}_{p}^{\ast }f:=\sup_{n\in \mathbb{N}}\frac{\left\vert
	S_{n}f\right\vert }{\left( n+1\right) ^{1/p-1}}
\end{equation*}%
is bounded from the Hardy space $H_{p}$\ to the Lebesgue space $L_{p}.$
Hence, according to inequality \eqref{LnSn} we can conclude that
\begin{equation*}
\left\Vert \underset{n\in \mathbb{N}}{\sup }\frac{\left|L_{n} f \right|}{\left(n+1\right)^{1/p-1}}\right\Vert_{p}\leq \left\Vert \underset{n\in \mathbb{N}}{\sup }\frac{\left|S_{n} f \right|}{\left(n+1\right)^{1/p-1}}\right\Vert _{p}\leq c_{p}\left\Vert f\right\Vert_{H_p}.
\end{equation*}
Now, prove part b) of the Theorem. Let

\begin{equation*}
f_{n_k}=D_{M_{2n_k+1}}-D_{M_{2n_k}}.
\end{equation*}

It is evident

\begin{equation}
\widehat{f}_{n_{k}}\left( i\right) =\left\{
\begin{array}{l}
\text{ }1,\qquad \text{ if }\qquad i=M_{2n_{k}},...,M_{2n_{k}+1}-1, \\
\text{ }0,\qquad \text{otherwise.}
\end{array}
\right.  \label{14a}
\end{equation}
Then we can write
\begin{equation} \label{14}
S_{i}f_{n_{k}} =\left\{
\begin{array}{l}
D_{i}-D_{M_{2n_{k}}}, \quad \text{ if } \quad
i=M_{2n_{k}}+1,...,M_{2n_{k}+1}-1, \\
f_{n_{k}},\qquad \qquad \quad \text{ if } \quad i\geq M_{2n_{k}+1}, \\
0,\qquad \qquad \qquad \text{ otherwise.}
\end{array}
\right.  
\end{equation}

From (\ref{3aa}) we get
\begin{eqnarray} \label{14c}
&&\left\| f_{n_{k}}\right\| _{H_{p}}\\  \notag
&=&\left\|
\sup\limits_{n\in \mathbb{N}}S_{M_{n}}f_{n_{k}}
\right\|_{p}  \\ \notag
&=&\left\| D_{M_{2n_{k}+1}}-D_{M_{_{2n_k}}} \right\|_p \\ \notag
&\leq& \left\| D_{M_{2n_{k}+1}} \right\|_p+ \left\|D_{M_{_{2n_{k}}}}\right\| _p \\ \notag
&\leq& cM_{_{2n_k}}^{1-1/p}<c<\infty.  \notag
\end{eqnarray}

Let $0<p<1$ and $\left\{ \lambda_k:k\in \mathbb{N}_{+}\right\}$ be an increasing sequence of the positive integers such that

\begin{equation*}
\lim_{k\rightarrow\infty}\frac{\lambda_{k}^{1/p-1}}{\varphi \left(
	\lambda _{k}\right)}=\infty.
\end{equation*}

Let$\ \left\{ n_{k}:k\in \mathbb{N}_{+}\right\} \subset \left\{ \lambda_k:k\in
\mathbb{N}_{+}\right\} $ such that

\begin{eqnarray*}
	\underset{k\rightarrow\infty}{\lim}\frac{\left(
		M_{_{2n_{k}}}+2\right)^{1/p-1}}{\log{\left( M_{2n_{k}}+2\right)}\varphi\left(
		M_{2n_k+2}\right)}\geq c\lim_{k\rightarrow \infty }\frac{\lambda _{k}^{1/p-1}}{\varphi \left(\lambda _{k}\right) }=\infty.
\end{eqnarray*}

According to (\ref{14}) we can conclude that

\begin{eqnarray*}
	&&\left| \frac{L_{M_{2n_{k}}+2}f_{n_{k}}}{\varphi\left(
		M_{2n_k+2}\right)} \right| \\
	&=&\frac{\left| D_{M_{_{2n_{k}}}+1} -D_{M_{_{2n_{k}}}}\right|}{l_{M_{2n_{k}}+1}\varphi\left(M_{2n_k+1}\right)}\\
	&=&\frac{\left|\psi _{M_{_{2n_{k}}}} \right|}{
		l_{M_{2n_{k}}+2}\varphi\left(
		M_{2n_k+1}\right)}\\
	&=&\frac{1}{l_{M_{2n_{k}}+1}\varphi\left(
		M_{2n_k+2}\right)}.
\end{eqnarray*}

Hence,
\begin{equation} \label{17}
\mu \left\{ x\in G_{m}:\left| L_{M_{2n_{k}}+2}f_{n_{k}}\right|\geq \frac{1}{l_{M_{2n_{k}}+2}\varphi\left(
	M_{2n_k+2}\right)}\right\}  
=\mu \left( G_{m}\right) =1.
\end{equation}

By combining (\ref{14c}) and (\ref{17}) we get that

\begin{eqnarray*}
	&&\frac{\frac{1}{l_{M_{2n_{k}}+2}\varphi\left(M_{2n_k+2}\right)} \left( \mu \left\{ x\in G_{m}:\quad\left| L_{M_{2n_{k}}+2} f_{n_{k}}  \right| 
		\geq\frac{1}{l_{M_{2n_{k}}+2}\varphi\left(
			M_{2n_k+2}\right)}\right\}\right)^{1/p}}{\left\| f_{n_{k}}\right\| _p} \\
	&\geq&\frac{M_{_{2n_{k}}}^{1/p-1}}{l_{M_{2n_{k}}+2}\varphi\left(	M_{2n_k+2}\right)} \\
	&\geq& \frac{c\left(M_{_{2n_{k}}}+2\right)^{1/p-1}} {\log{\left(M_{2n_{k}}+2\right)}\varphi\left(M_{2n_k+2}\right)}\rightarrow \infty, \text{\quad as \quad }k\rightarrow \infty.
\end{eqnarray*}

Theorem is proved.

\QED

\begin{corollary} Let $0<p<1\ $and $f\in H_{p}\left( G_{m}\right) .$ Then
	there is an absolute constant $c_{p},$ depends only $p,$ such that
	
	\begin{equation*}
	\overset{\infty }{\underset{k=1}{\sum }}\frac{\left\Vert L_{k}f\right\Vert
		_{p}^{p}}{k^{2-p}}\leq c_{p}\left\Vert f\right\Vert _{H_{p}}^{p}.
	\end{equation*}
\end{corollary}
\textbf{Proof:} If we apply Lemma \ref{lemma2.2} we only have to prove that

\begin{equation*}
\overset{\infty }{\underset{k=1}{\sum }}\frac{\left\Vert L_{k}a\right\Vert
	_{p}^{p}}{k^{2-p}}\leq c_{p}<\infty.
\end{equation*}
for every $p$ atom $a$, with support $I_{N}$ and $\mu \left( I_{N}\right)
=M_{N}^{-1}$. Since $S_{k}a=0,$ for $k\leq M_{N},$ we get that $L_{n}a=0,$ for $n\leq M_{N}.$ So, we may assume that $n>M_{N}.$

Let $x\in \overline{I_{N}}$ and $0<p<1.$ Applying (\ref{1.1}) and Lemma \ref{dn2.7} we have that

\begin{eqnarray*}
	&&\overset{\infty }{\underset{k=M_{N}}{\sum }}\frac{\left\Vert
		L_{k}a\right\Vert _{p}^{p}}{k^{2-p}} \\
	&\leq & \overset{\infty }{\underset{k=M_{N}}%
		{\sum }}\frac{1}{k}\int_{\overline{I_{N}}}\left\vert \frac{L_{k}a}{k^{1/p-1}}\right\vert ^{p}d\mu  \\
	&=&\overset{\infty }{\underset{k=M_{N}}{\sum }}\frac{1}{k}\overset{N-1}{
		\underset{s=0}{\sum }}\int_{I_{s}\backslash I_{s+1}}\left\vert \frac{
		L_{k}a }{k^{1/p-1}}\right\vert ^{p}d\mu  \\
	&\leq & c\overset{\infty }{\underset{k=M_{N}}{\sum }}\frac{1}{k}\overset{N-1}{%
		\underset{s=0}{\sum }}\int_{I_{s}\backslash I_{s+1}}\left\vert \frac{%
		M_{N}^{1/p-1}M_{s}}{k^{1/p-1}}\right\vert ^{p}d\mu \\
	&\leq & c_{p}M_{N}^{1-p}\overset{\infty }{\underset{k=M_{N}}{\sum }}\frac{1}{%
		k^{2-p}}\overset{N-1}{\underset{s=0}{\sum }}\int_{I_{s}\backslash
		I_{s+1}}M_{s}^{p} \\
	&\leq & c_{p}M_{N}^{1-p}\overset{\infty }{\underset{k=M_{N}}{\sum }}\frac{1}{%
		k^{2-p}}\overset{N-1}{\underset{s=0}{\sum }}M_{s}^{p-1}d\mu \\
	&+&c_{p}M_{N}^{1-p}\overset{\infty }{\underset{k=M_{N}}{\sum }}\frac{1}{k^{2-p}%
	}\leq c_{p}<\infty .
\end{eqnarray*}

Which complete the proof.
\QED

\begin{theorem}\label{Corollaryconv4} Let $p\geq 1$ and $f\in L_p$. Then
	$$\Vert L_{M_n} f-f\Vert_p \to 0 \ \ \text{as}\ \ n\to \infty$$
	for all  $f\in L_p(G_m)$. Moreover,
	\begin{equation*}
	\underset{n\rightarrow \infty }{\lim }L_{M_n}f(x)=f(x)
	\end{equation*}
	for all Lebesgue points of $f\in L_p(G_m)$.
\end{theorem}

\textbf{Proof:} To prove a.e convergence we use first identity in Lemma \ref{lemma0nnT121}  to write that

\begin{eqnarray*}
	L_{M_n}f\left(x\right)&=&\underset{G_m}{\int}f\left(t\right)P_n\left(x-t\right) d\mu\left(t\right)\\
	&=&\underset{G_m}{\int}f\left(t\right)D_{M_n}\left(x-t\right)d\mu\left(t\right)\\
	&-&\underset{G_m}{\int}f\left(t\right)\psi_{M_n-1}(x-t)\overline{Y}_{M_n}(x-t)\\
	&=&I-II.
\end{eqnarray*}

By applying  \eqref{covstrong} we can conclude that
$$I=S_{M_n}f(x)\to f(x)$$
for all Lebesgue points of $f\in L_p(G_m)$.
By using \eqref{vilenkin} we can conclude that
$$II=\psi_{M_n-1}(x)\underset{G_m}{\int}f\left(t\right)\overline{Y}_{M_n}(x-t)\overline{\psi}_{M_n-1}(t)d(t)$$
By combining \ref{1.71alphaT2jj} in Lemma \ref{lemma0nnT121} and  \ref{covstrong} we find that  $$f\left(t\right)\overline{Y}_{M_n}(x-t)\in L_p \ \ \text{ where} \ \  p\geq  1 \ \ \text{for any } \ \ x\in G_m, $$
and $II$ is Fourier coefficients of integrable function. According to Riemann-Lebesgue Lemma it vanish as $n\to \infty,$
$$II\to 0 \ \ \text{for any } \ \ x\in G_m, \ \  n\to \infty.$$
The proof is complete.

\QED
\newpage




\appendix

\printindex{A}{Author Index}
\printindex{N}{Notation Index}
\printindex{S}{Subject Index}

\newpage

\end{document}